\documentclass[10pt,reqno]{NumPDEsArticle}

\usepackage{NumPDEsMacros}

\usepackage{csquotes}
\usepackage{enumitem}
\usepackage{amsmath,amssymb}
\usepackage{empheq}
\usepackage{graphicx}
\usepackage{tikz}
\usepackage{caption}
\usepackage{subcaption}
\usepackage[capitalize,nameinlink]{cleveref}
\usepackage{float}
\usepackage{booktabs}
\usepackage[width=\textwidth]{caption}
\usepackage{cleveref} 

\DeclareMathAlphabet\mathbfcal{OMS}{cmsy}{b}{n}

\title[BD2-type integrator for LLG, part I]{BDF2-type integrator for Landau-Lifshitz-Gilbert equation in micromagnetics: unconditional weak convergence to weak solutions}

\author{Michele Ald\'e}
\email{Michele.Alde@asc.tuwien.ac.at \quad\rm (corresponding author)}

\author{Michael Feischl}
\email{Michael.Feischl@asc.tuwien.ac.at}

\author{Dirk Praetorius}
\email{Dirk.Praetorius@asc.tuwien.ac.at}

\address{TU Wien, Institute of Analysis and Scientific Computing, Wiedner Hauptstra\ss{}e 8--10, 1040 Wien, Austria}

\thanks{This research was funded by the Austrian Science Fund (FWF) projects \href{https://www.doi.org/10.55776/F65}{10.55776/F65} (SFB F65 ``Taming complexity in PDE systems''), \href{https://www.doi.org/10.55776/I6802}{10.55776/I6802} (international project I6802 ``Functional error estimates for PDEs on unbounded domains''), and \href{https://www.doi.org/10.55776/P33216}{10.55776/P33216} (stand-alone project P33216 ``Computational nonlinear PDEs'') and through the ERC grant \href{https://www.doi.org/10.3030/101125225}{10.3030/101125225} ``New frontiers in optimal adaptivity'' of Michael Feischl. Additionally, the \href{https://www.vsmath.at}{Vienna School of Mathematics} supports Michele Ald\'e.}

\newtheorem*{proposition*}{Proposition}

\newcommand{\mm}{\boldsymbol{m}}
\newcommand{\heff}{\boldsymbol{h}_{\rm eff}}
\newcommand{\ff}{\boldsymbol{f}}
\newcommand{\vv}{\boldsymbol{v}}
\newcommand{\pphi}{\vvarphi}
\newcommand{\bphi}{\pphi}
\newcommand{\bpsi}{\boldsymbol{\psi}}
\newcommand{\lamex}{\lambda_{\rm ex}}
\newcommand{\bT}{\boldsymbol{T}}
\newcommand{\bL}{\mathbf{L}}
\newcommand{\bH}{\mathbf{H}}
\newcommand{\OT}{\Omega_T}
\newcommand{\mhatp}{\widehat{\mm}_{h\tau}^+}
\newcommand{\vht}{\vv_{h\tau}^-}
\newcommand{\scalarproductO}[2]{\langle #1, #2 \rangle_{\Omega}}
\newcommand{\scalarproductOT}[2]{\langle #1, #2\rangle_{\OT}}
\newcommand{\normO}[2][]{\norm[#1]{#2}_{\Omega}}
\newcommand{\normOT}[2][]{\norm[#1]{#2}_{\Omega_T}}

\renewcommand{\set}[3][]{#1\{#2 \,:\, #3#1\}}
\newcommand{\xx}{\boldsymbol{x}}
\newcommand{\uu}{\boldsymbol{u}}
\def\III{\boldsymbol{\II}}
\def\SSS{\boldsymbol{\SS}}
\def\pphi{\boldsymbol{\phi}}
\def\vvarphi{\boldsymbol{\varphi}}

\begin{document}
	
\maketitle

\begin{abstract}
We consider the Landau--Lifshitz--Gilbert equation (LLG) that models time-dependent micromagnetic phenomena.
We propose a full discretization that employs first-order finite elements in space and a BDF2-type two-step method in time. In each time step, only one linear system of equations has to be solved.
We employ linear interpolation in time to reconstruct  the discrete space-time magnetization.
We prove that the integrator is unconditionally stable and thus guarantees that a subsequence of the reconstructed magnetization converges weakly in $H^1$ towards a weak solution of LLG in the space-time domain.
Numerical experiments verify that the proposed integrator is indeed first-order in space and second-order in time.
\end{abstract}
\section{Introduction}
Time-dependent micromagnetic phenomena are usually modeled via the phenomenological Landau--Lifshitz--Gilbert equation (LLG), which is a time-dependent nonlinear partial differential equation (PDE).
LLG describes the evolution of the magnetization $\mm(t,\xx) \in \R^3$ modeling the ferromagnetic spin at time $t$ and position $\xx \in \Omega$ of a ferromagnetic body $\Omega \subset \R^{d}$, $d=2,3$.
It is physically known and inherent part of the PDE model that the modulus $|\mm(t,\xx)| = M_s$ is constant for constant temperatures below the so-called Curie temperature, where $M_s > 0$ is the so-called saturation magnetization.

The understanding of LLG is of utmost importance for technical applications like sensors, actuators, or recording devices~\cite{LandauA1935,BrownB1962,BrownB1963,hubert2008magnetic}.
In the last decades, the mathematical understanding of LLG has matured.
We refer to~\cite{Alouges_Soyeur} for the global-in-time existence of weak solutions, \cite{csg1998, cf2001, cimrak2005, cimrak2007, ft2017b,MR4646547} for local-in-time existence of smooth strong solutions, and~\cite{ds2014, dip2020} for weak-strong uniqueness, i.e., if a strong solution exists until time $t'$, then it is also the unique weak solution until time $t'$.

As far as the recent development of numerical integrators is concerned, these can roughly be classified into integrators that converge to strong solutions at certain rates \cite{prohl2001, bs2006, MR3273326, MR3542789, feischl, MR4454924, MR4362832,MR4673464,MR4617914} and integrators that guarantee weak convergence (along a subsequence) to weak solutions of LLG~\cite{aj2006, bkp2008, bp2006, Alouges_2008, ruggeri2022}.
In particular, certain attention has been paid to effective integrators of coupled LLG systems~\cite{bppr2014, bpp2015, dpprs2020} that use a split implicit-explicit time integration (IMEX) to reduce the computational cost~\cite{akt2012, bsffgppr2014, mrs2018, hpprss2019, dppr2023}.

From all these works, we emphasize the two seminal works~\cite{bp2006, Alouges_2008} that firstly proved unconditional weak convergence to weak solutions of LLG.
Both integrators rely on first-order finite elements in space.
While~\cite{bp2006} employs the implicit midpoint-rule in time, which is formally of second-order in time but requires to solve one nonlinear system of equations per time step, the work~\cite{Alouges_2008} proposes a so-called tangent-space integrator, which is of first-order in time but requires only the solution of one linear equation per time step.

From the viewpoint of computational mathematics, the tangent-plane method is particularly attractive.
Therefore,~\cite{Alouges_2008} triggered several generalizations.
The work \cite{akst2014} proposed a variant of the original algorithm that is still unconditionally convergent to weak solutions, but formally almost second-order in time, and~\cite{dpprs2020} formulated and analyzed an IMEX variant of this.
In both works, the proof of \textsl{a~priori} error estimates in the presence of a smooth strong solution remains open.
Instead, the work~\cite{feischl} employs a modified discrete tangent-plane and provides \textsl{a~priori} error estimates for various BDF-type integrators in the presence of smooth strong solutions of LLG.
However, the unconditional weak convergence to weak solutions of LLG remains open in~\cite{feischl}.

In a series of three papers, our work aims to formulate and analyze a numerical integrator that is second-order in time and first-order in space having the following properties: First, the integrator guarantees unconditional weak convergence to weak solutions of LLG.
Second~\cite{part2}, the integrator converges with given rates to the unique strong solution of LLG, if the latter exists and is sufficiently smooth.
Third~\cite{part3}, the integrator can be extended from mathematical model problems to the full effective field via an efficient IMEX strategy.
Similar results are so far only available for the first-order tangent-plane scheme, with weak convergence analyzed in \cite{Alouges_2008} and \textsl{a~priori} error estimates proved in \cite{AnLiSun2024}. Additionally, the same results were proved also for the projection-free variant of this scheme \cite{Spin-polarized,bartels2016}, where unconditional convergence and IMEX time-stepping is analyzed in \cite{Spin-polarized}, while \textsl{a~priori} error estimates are shown in \cite{ft2017}.

The analysis in the present paper adapts recent ideas from~\cite{abp2024} for a BDF2-type minimization strategy for harmonic maps and the tangent-plane BDF2 integrator of~\cite{feischl}.
While~\cite{feischl} employs a discrete tangent space defined in a variational sense but arbitrary-order finite elements in space, our algorithm employs a node-wise discrete tangent space as in~\cite{Alouges_2008} that is tailored to first-order finite elements in space. Other related work includes~\cite{MR4727106}, which provides error estimates for a discretization scheme of the harmonic map heat flow into spheres and also uses constraints to enforce the orthogonality condition in the discrete tangent space.
Instead of the first-order backward Euler method used in~\cite{Alouges_2008}, our integrator relies on a predictor-corrector approach as in~\cite{feischl}, where the corrector is based on a BDF2-step. 
As in~\cite{Alouges_2008} and the succeeding works, the proposed time-marching scheme leads to lowest-order finite element approximations $\mm_h^j \approx \mm(t_j)$ of a weak solution $\mm$ of LLG at times $t_j$, which are then interpolated in time to provide a discrete space-time approximation $\mm_{h\tau} \approx \mm$ with $\mm_{h\tau}(t_j) = \mm_h^j$.
Then, our analysis adapts the classical energy method for parabolic PDEs to show unconditional weak convergence of a subsequence of $\mm_{h\tau}$ towards $\mm$. As usual in this context, the proof provides a constructive argument for the existence of weak solutions to LLG.


\textbf{Outline.}
%
The outline of the present work reads as follows: 
In Section~\ref{section:model}, we formulate LLG~\eqref{eq: LLG system} in the presence of a simplified effective field consisting only of an applied exterior field and the so-called exchange contribution, which is the leading-order term.
In Section~\ref{section:discretization}--\ref{section:integrator}, we give the precise statement of the proposed numerical integrator (Algorithm~\ref{alg: full discr}).
Section~\ref{section:theorem} recalls the definition of weak solutions to LLG from~\cite{Alouges_Soyeur} and then formulates the main result of this work (Theorem~\ref{thm: main result}) stating unconditional weak convergence of the computed approximations towards a weak solution of LLG.
Its proof is given in Section~\ref{section:proof}, which takes the main part of this paper.
In Section~\ref{section:numerics}, some numerical experiments underline that the proposed integrator is indeed of first-order in space and second-order in time.

\textbf{General notation.}
%
Throughout this paper, we adhere to the standard notation for Lebesgue, Sobolev, and Bochner spaces, using bold letters to indicate spaces of vector valued functions, e.g., $\bL^2(\Omega) \coloneqq L^2(\Omega,\R^3)$ and $\bH^1(\Omega)\coloneqq H^1(\Omega,\R^3)$.
Given $\uu, \vv \in \bL^2(\Omega)$, we will denote the $\bL^2$-scalar product in $\Omega\subset\R^d$, $d=2,3,$ by 
\begin{align*}
 \langle\uu,\vv\rangle_\Omega 
 \coloneqq 
 \langle\uu,\vv\rangle _{\bL^2(\Omega)}= \int_\Omega\uu(\xx)\cdot\vv(\xx)\,\d\xx,
\end{align*}
and we will use $\normO{\,\cdot\,}$ to indicate the induced norm, i.e., $\normO{\uu} \coloneqq \langle\uu,\uu\rangle_\Omega^{1/2}$.
Similar notation is applied for the space-time domain $\OT\coloneqq \Omega \times [0,T]$ and the corresponding scalar product $\langle\cdot,\cdot\rangle_{\OT}$ and norm $\normOT{\,\cdot\,}$.
All other norms will be explicitly stated.

By $C > 0$, we will denote a generic constant, which is always independent of the discretization parameters, but not necessarily the same at each occurrence.
We will use the notation $A\lesssim B$ to indicate that $A\leq C \, B$, where the hidden constant $C > 0$ is clear from the context.
Moreover, $A \simeq B$ abbreviates $A \lesssim B \lesssim A$.

\section{Numerical BDF2 integrator and main result}

\subsection{Mathematical model}\label{section:model}

Let $\Omega \subset \R^{d}$, $d=2,3$, be a bounded Lipschitz domain with boundary $\Gamma \coloneqq \partial\Omega$ describing a ferromagnetic body. We note that the 2-dimensional setting is related to so-called thin-film limits, which investigate the behavior of micromagnetic phenomena in samples where one dimension (the thickness) is much smaller than the other two; see, e.g., \cite{davoli2020} for details and the mathematical derivation of a thin-film model with $d = 2$ from the general physical setting with $d = 3$.
The dynamics of the normalized magnetization $\mm \in \S^2 \coloneqq \set{\xx \in \R^3}{|\xx|=1}$ is governed by the Landau--Lifshitz--Gilbert equation (LLG):
Given the final time $T > 0$, find $\mm\colon\Omega\times[0,T]\to\S^2$ such that
\begin{subequations}\label{eq: LLG system}
\begin{empheq}[left=\empheqlbrace]{align}
 &\partial_t\mm = -\mm\times\heff + \alpha \mm\times\partial_t\mm 
  &&\hspace{-2em} \text{in} \quad \Omega \times [0,T] \label{eq: LLG}, \\
 &\partial_{\vec{n}}\mm = 0 
  &&\hspace{-2em} \text{in} \quad \partial\Omega \times [0,T] \label{eq:subeq_b}, \\
 &\mm(0) = \mm^0 
  &&\hspace{-2em} \text{in} \quad \Omega \label{eq:subeq_c},
\end{empheq}
where~\eqref{eq: LLG} is the so-called Gilbert form of LLG, $\mm^0\colon \Omega\to\S^2$ is the given initial state and $\alpha>0$ is the non-dimensional Gilbert damping parameter.
The dynamics of~\eqref{eq: LLG system} is driven by the effective magnetic field 
\begin{align}
 \heff\coloneqq\heff(\mm,\ff)=-\frac{\delta\EE(\mm,\ff)}{\delta\mm},
\end{align}
where $\EE(\mm,\ff)$ is the total magnetic Gibbs free energy of the system.
In this work, the energy consists of exchange contribution (also called Heisenberg exchange) and an applied external field contribution (also called Zeeman's energy), i.e.,
\begin{equation}\label{eq: energy}
 \EE(\mm,\ff)
 =
 \frac{\lamex^2}{2} \int_\Omega|\nabla\mm|^2\d \xx - \int_\Omega\ff\cdot\mm \d \xx
\end{equation} 
\end{subequations}
so that $\heff(\mm,\ff)=\lamex^2\Delta\mm+\ff$.
The formal application of the chain rule yields 
\begin{equation}\label{eq: energy law}
 \frac{\d}{\d t}\EE(\mm,\ff)
 =
 -\alpha \int_\Omega|\partial_t\mm(t)|^2\d \xx - \int_\Omega \partial_t\ff(t) \cdot\mm(t) \d \xx,
\end{equation}
revealing the dissipative behavior of LLG~\eqref{eq: LLG system} if $\ff$ is constant.

\subsection{Discretization}\label{section:discretization}

Given the final time $T>0$ and $N \in \N$, we consider a uniform partition of the time interval $[0,T]$ with time-step size $\tau\coloneqq T/N$ and time steps $t_j\coloneqq j\tau$ for all $j = 0,\dots,N$.
For a sequence of functions $\{\mm^j\}_{j=0,\dots,N}$ with $\mm^j$ being associated with the time step $t_j$, we define the first-order discrete time-derivative
\begin{align*}
 \d _t{\mm}^{j+1} 
 \coloneqq 
 \dfrac{1}{\tau}\big(\mm^{j+1}-\mm^j\big)
\end{align*}
as well as the second-order discrete time derivative
\begin{equation}\label{eq: second order discrete time derivative}
 \d_t^2{\mm}^{j+1}
 \coloneqq 
 \dfrac{1}{\tau}\big(\d _t \mm^{j+1}-\d _t \mm^j\big)
 =
 \frac{1}{\tau^2}\big(\mm^{j+1}-2\mm^j+\mm^{j-1}\big).
\end{equation}
For the space discretization, let $\TT_h$ be a $\gamma$-quasi-uniform family of regular tetrahedral meshes of the domain $\Omega$ with mesh size $h>0$, i.e., 
\begin{align*}
 \gamma^{-1} h \le |K|^{1/d} \le \text{diam}(K) \le h
 \quad \text{for all } K \in \TT_h \text{ and all } h > 0,\, d=2,3.
\end{align*} 
Let $\SSS^1(\TT_h) \coloneqq \set{\vv_h \in C(\Omega; \R^3)}{\forall K \in \TT_h: \vv_h|_{K} \text{ is affine}}$
be the space of globally continuous and $\TT_h$-piecewise affine functions.
Let $\NN_h$ be the set of vertices of $\TT_h$.

Multiplying~\eqref{eq: LLG} by $\mm$, we note that~\eqref{eq: LLG} formally ensures that 
\begin{equation}\label{eq: orthogonality property}
	\dfrac{1}{2}\partial_t|\mm|^2=\mm\cdot\partial_t\mm\stackrel{\eqref{eq: LLG}}{=}0
	\quad \text{in }\OT.
\end{equation}
With $|\mm^0| = 1$, the magnetization $\mm$ thus satisfies the pointwise modulus constraint $|\mm| = 1$ in $\OT$. 
Moreover, the time-derivative $\partial_t \mm$ is pointwise orthogonal to $\mm$.

To mimic the constraint $|\mm| = 1$ at a discrete level, we consider the set
\begin{align*}
 \boldsymbol{\mathcal{M}}_h 
 \coloneqq 
 \set{\bphi_h\in\SSS^1(\TT_h) }{\forall \vec{z}\in\NN_h: |\bphi_h(\vec{z})|=1}
\end{align*}
of admissible discrete magnetizations satisfying the constraint at every $\vec{z}\in \NN_h$.
To mimic the orthogonality property $\mm \cdot \partial_t\mm = 0$ at a discrete level, we will enforce it at every $\vec{z}\in\NN_h$.
To this end, we define the discrete tangent space at $\pphi_h\in\SSS^1(\TT_h)$ via 
\begin{equation*}
 \bT_h(\pphi_h) 
 \coloneqq 
 \set[\big]{\bpsi_h\in\SSS^1(\TT_h)}{\forall \vec{z} \in \NN_h: \bphi_h(\vec{z}) \cdot \bpsi_h(\vec{z}) = 0}.
\end{equation*}

\subsection{Numerical integrator}\label{section:integrator}

We make the following assumptions on the given data of LLG~\eqref{eq: LLG system} and their discretization:
\begin{enumerate}[label=(A\arabic*), ref=A\arabic*]
 \item\label{enum: punto A2} 
 $\mm^0 \in \bH^1(\Omega)$, and $\mm_h^0 \in \boldsymbol{\mathcal{M}}_h$ satisfies $\mm_h^0 \to \mm^0$ in $\bH^1(\Omega)$ as $h\to 0$.
 \item\label{enum: punto A1}  
 $\ff\in C^1([0,T], \bL^2(\Omega))$ and the discretization $\ff_{h\tau}^+ \in L^2([0,T], \boldsymbol{\SS^1}(\TT_h))$ is given by $\ff_{h\tau}^+(t,\vec{z})\coloneqq \ff_h^{j+1}(\vec{z}) \coloneqq \ff(t_{j+1},\vec{z})$ for all $t\in[t_j,t_{j+1})$ and $\vec{z} \in \NN_h$, i.e., $\ff_{h\tau}^+$ is piecewise constant in time and the nodal interpoland of $\ff(t_{j+1})$ in space so that, in particular, $\ff_{h\tau}^+ \to \ff$ in $\bL^2(\Omega)$ as $(h,\tau)\to(0,0)$.
\end{enumerate}
Clearly, the assumptions~\eqref{enum: punto A2}--\eqref{enum: punto A1} yield uniform boundedness $\normO{\nabla \mm_h^0}\le C_0$ and $\normOT{\ff_{h\tau}^+}\le C_f$ for all $h > 0$ and $\tau > 0$.

Under the constraint $|\mm|=1$, one can prove that~\eqref{eq: LLG} is equivalent to 
\begin{equation}\label{eq: alternative LLG}
	\alpha\partial_t\mm + \mm\times\partial_t \mm = \heff-(\mm\cdot \heff)\mm;
\end{equation}
see, e.g.,~\cite{goldenits2012konvergente}.
For the numerical integration of~\eqref{eq: LLG system}, the proposed time-marching scheme thus exploits a finite element discretization of~\eqref{eq: alternative LLG} for the time-derivative $\vv \coloneqq \partial_t\mm$ in the tangent space of $\mm$ and uses second-order backward differences (BDF2) in time.
Since BDF2 requires an additional initial value $\mm_h^1 \approx \mm(t_1)$ beyond $\mm_h^0 \approx \mm^0$, we employ one step of the non-projected first-order tangent-plane integrator for LLG~\cite{Spin-polarized}.
Overall, the proposed numerical integrator reads as follows:

\begin{algorithm}\label{alg: full discr}
\textbf{Input:} Conforming mesh $\TT_h$ of $\Omega$, $\mm(0) \approx \mm_h^0\in \boldsymbol{\MM}_h$, $T>0$, $N\in\N$, $\tau\coloneqq T/N$, and $t_j\coloneqq j\tau$ for all $j=0,\dots,N$.
\begin{enumerate}[label=\rm (\roman*), ref=\rm \roman*]
\item\label{eq: first step full discr item} 
For $j=0,$ 
compute $\vv_h^0\in \bT_h(\mm_h^0)$ such that, for all $\pphi_h\in \bT_h(\mm_h^0)$,
\begin{equation}\label{eq: first step full discr}
 \alpha\scalarproductO{\vv_h^0}{\pphi_h} + \scalarproductO{\mm_h^0\times\vv_h^0}{\pphi_h} + \lamex^2\tau\scalarproductO{\nabla\vv_h^0}{\nabla\pphi_h}
 =
 \scalarproductO{\ff_h^1}{\pphi_h} - \lamex^2\scalarproductO{\nabla\mm_h^0}{\nabla\pphi_h}.
\end{equation}
\item\label{item: def mm1}
Set $\mm_h^1\coloneqq \mm_h^0+\tau\vv_h^0$.
\item\label{item: next steps}  
For $j=1,\dots,N-1$, repeat the following steps~\eqref{item: def mhat}--\eqref{item: def mmhj}:
\begin{enumerate}[label=\rm(\alph*), ref=\rm \alph*]
 \item\label{item: def mhat}
 Set $\widehat{\mm}_h^{j+1}=2\mm_h^{j}-\mm_h^{j-1}$.
 \item\label{item: var form}
 Compute $\vv_h^j\in \bT_h(\widehat{\mm}_h^{j+1})$ such that, for all $\pphi_h\in \bT_h(\widehat{\mm}_h^{j+1})$,
 \begin{equation}\label{eq: full discr alg B}
 \begin{split}
 \alpha\scalarproductO{\vv_h^j}{\pphi_h} &+ \scalarproductO{\widehat{\mm}_h^{j+1}\times\vv_h^j}{\pphi_h} + \frac{2}{3}\lamex^2\tau\scalarproductO{\nabla\vv_h^j}{\nabla\pphi_h}
  \\&
  =
  \scalarproductO{\ff_h^{j+1}}{\pphi_h} - \frac{1}{3}\lamex^2\scalarproductO{\nabla[4\mm_h^{j}-\mm_h^{j-1}]}{\nabla\pphi_h}.
 \end{split}
 \end{equation}
 \item\label{item: def mmhj}
 Set $\mm_h^{j+1}\coloneqq \dfrac{4}{3}\mm_h^{j}-\dfrac{1}{3}\mm_h^{j-1}+\dfrac{2}{3}\tau \vv_h^j$.
\end{enumerate}
\end{enumerate}
\textbf{Output:} Sequences $\vv_h^j\approx\partial_t\mm(t_j)$ and $\mm_h^{j+1}\approx\mm(t_{j+1})$ for all $j=0,\dots,N-1$.
\end{algorithm}

Thanks to the Lax--Milgram lemma, it follows that~\eqref{eq: first step full discr}--\eqref{eq: full discr alg B} admit unique solutions $\vv_h^j$ so that 
Algorithm~\ref{alg: full discr} is indeed well-defined.

\subsection{Convergence theorem}\label{section:theorem}

Throughout this work, we consider weak solutions of LLG~\eqref{eq: LLG system} in the sense of \cite[Definition~1.2]{Alouges_Soyeur}:

\begin{definition}(Weak solution)\label{def: weak solution}
Let $\mm^0\in\bH^1(\Omega)$ satisfy $|\mm^0|=1$ a.e.\ in $\Omega$.
A function $\mm: \OT \rightarrow \mathbb{R}^3$ is a \text{\rm weak solution} of~\eqref{eq: LLG system} if the following properties {\rm (i)-(iv)} are satisfied:
\begin{enumerate}[label={\rm(\roman*)}, ref={\rm\roman*}]
\item\label{def: weak solution i} 
$\mm \in  \bH^1(\OT) \cap L^{\infty}(0, T ; \bH^1(\Omega))$ and $|\mm|=1$ a.e.\ in $\OT$;
\item\label{def: weak solution ii}
$\mm(0)=\mm^0$ in the sense of traces;
\item\label{def: weak solution iii}
for all $\vvarphi \in \bH^1(\OT)$, it holds that 
$$
\begin{aligned}
& \int_0^T\scalarproductO{\partial_t \mm(t)}{\vvarphi(t)} \d t-\alpha \int_0^T\scalarproductO{\mm(t) \times \partial_t \mm(t)}{\vvarphi(t)} \d t \\
&= \lambda_{\mathrm{ex}}^2 \int_0^T\scalarproductO{\mm(t) \times {\nabla} \mm(t)}{{\nabla} \vvarphi(t)} \d t-\int_0^T\scalarproductO{\mm(t) \times \boldsymbol{f}(t)}{\vvarphi(t)} \d t;
\end{aligned}
$$
\item\label{def: weak solution iv} 
for a.e.\ $t' \in(0, T)$, $\mm$ satisfies the energy inequality 
$$
\begin{aligned}\label{eq: energy inequality}
\EE(\mm(t'), \boldsymbol{f}(t'))+\alpha \int_0^{t'}\normO[\big]{\partial_t \mm(t)}^2 \mathrm{~d} t+\int_0^{t'}\langle\partial_t \boldsymbol{f}(t), \mm(t)\rangle _\Omega \d t\leq \EE(\mm^0, \boldsymbol{f}(0)).
\end{aligned}
$$
\end{enumerate}
\end{definition}

Definition~\ref{def: weak solution}\eqref{def: weak solution i} formalizes the modulus constraint~\eqref{eq: orthogonality property}, while Definition~\ref{def: weak solution}\eqref{def: weak solution ii} incorporates the initial condition~\eqref{eq:subeq_c}.
The variational formulation of Definition~\ref{def: weak solution}\eqref{def: weak solution iii} comes from a weak formulation of~\eqref{eq: LLG}--\eqref{eq:subeq_b} in the space-time domain $\OT$ after integrating by parts.
The energy inequality of Definition~\ref{def: weak solution}\eqref{def: weak solution iv} is the weak counterpart of the energy law~\eqref{eq: energy law}.

Using the approximations $\{\mm_h^j\}_{0\le j\le N}$ and $\{\vv_h^j\}_{0\le j \le N-1}$ obtained by Algorithm~\ref{alg: full discr}, we define approximations $\mm_{h\tau}\in\bH^1(\OT)$, $\mm_{h\tau}^{\pm},\mhatp,\vv_{h\tau}^-\in L^2(0,T;\bH^1(\Omega))$ of $\mm$ and $\vv=\partial_t\mm$ as follows:
For all $0\le j\le N-1$ and all $t\in[t_j,t_{j+1}),$ let
	\begin{align}
		\mm_{h\tau}(t) & \coloneqq \frac{t-t_j}{\tau}\mm_h^{j+1} + \frac{t_{j+1}-t}{\tau}\mm_h^{j},\label{eq:mhtau}\\
		\mm_{h\tau}^-(t) & \coloneqq \mm_h^j, \label{eq:mhtau-minus}\\
		\mm_{h\tau}^+(t) & \coloneqq \mm_h^{j+1}, \label{eq:mhtau-plus}\\
		\widehat{\mm}_{h\tau}^+(t) & \coloneqq
			\begin{cases}
				\mm_h^0 & \text{for } j=0,\\
				\widehat{\mm}_h^{j+1} & \text{for } 1 \leq j \leq N-1,
			\end{cases} \label{eq:mhtau-hat-plus}\\
		\vv_{h\tau}^-(t) & \coloneqq \vv_h^j.\label{eq:vhtau-minus}
	\end{align}
The following theorem, which is the main result of this work, states that these interpolands converge weakly towards a weak solution of LLG~\eqref{eq: LLG system} as $(h,\tau)\to(0,0)$.

\begin{theorem}\label{thm: main result}
For $h > 0$, let $\TT_h$ be a family of $\gamma$-quasi-uniform meshes of $\Omega$.
Under the assumptions {\rm~\eqref{enum: punto A2}} and {\rm~\eqref{enum: punto A1}}, the sequences of discrete functions $\mm_{h\tau}$, $\mm_{h\tau}^\pm$ and $\mhatp$, admit subsequences (not relabeled) that unconditionally converge to a function $\mm$ which is a weak solution of LLG~\eqref{eq: LLG system} according to Definition~\ref{def: weak solution}\eqref{def: weak solution i}--\eqref{def: weak solution iii}.
More precisely, it holds that $\mm_{h\tau}\rightharpoonup\mm$ weakly in~$\bH^1(\OT)$ (and hence $\mm_{h\tau} \to \mm$ strongly in $\bL^2(\OT)$) and $\mm_{h\tau},\mm_{h\tau}^\pm,\mhatp\stackrel{*}{\rightharpoonup}\mm$ in $L^\infty(0,T;\bH^1(\Omega))$ as $(h,\tau)\to(0,0)$,
where all convergences hold for the same subsequence.
Moreover, define 
\begin{equation}\label{eq: eta0n}
	\eta_0 \coloneqq \normO{\nabla\mm_h^1}^2-\normO{\nabla\mm_h^0}^2\qquad\text{and}\qquad\eta_n\coloneqq\normO{\nabla\mm_h^{n-1}}^2-\normO{\nabla\mm_h^{n}}^2.
\end{equation}
Then, the limit $\mm$ satisfies also Definition~\ref{def: weak solution}\eqref{def: weak solution iv}, provided that 
\begin{equation}\label{eq: condition eta0n}
\eta_0+\eta_n\xrightarrow{(h,\tau)\to(0,0)} 0.
\end{equation}
In particular, \eqref{eq: condition eta0n} is guaranteed under the CFL-type condition $\tau = o(h^2)$, i.e., $\sqrt{\tau}h^{-1}\to0$ as $(h,\tau)\to(0,0)$, which is, however,
only needed for the first and last time step.
\end{theorem}

\begin{remark}
We note that the limit function $\mm$ of Theorem~\ref{thm: main result} exists and satisfies Definition~\ref{def: weak solution}\eqref{def: weak solution i}--\eqref{def: weak solution iii} already under weaker assumptions on $\mm_h^0$ and $\ff$, namely weak convergence $\mm_h^0 \rightharpoonup \mm^0$ in $\bH^1(\Omega)$ and $\ff \in \bL^2(\OT)$ with $\ff_{h\tau}^+ \rightharpoonup \ff$ in $\bL^2(\OT)$, while the stronger assumptions from~\eqref{enum: punto A2}--\eqref{enum: punto A1} are only needed to verify Definition~\ref{def: weak solution}\eqref{def: weak solution iv}.
\end{remark}

\begin{corollary}\label{cor: main result}
Under the assumptions of Theorem~\ref{thm: main result} and provided that the weak solution $\mm$ to LLG~\eqref{eq: LLG system} is unique (e.g., strong and sufficiently smooth~\cite{dip2020}), the limits of Theorem~\ref{thm: main result} do not only hold for subsequences but also for the full sequences $\mm_{h\tau}$, $\mm_{h\tau}^\pm$ and $\mhatp$ as $(h,\tau)\to(0,0)$.
\end{corollary}

\section{Proof of theorem~\ref{thm: main result}}\label{section:proof}

\subsection{Prerequisites}

We recall some auxiliary results that will be needed throughout this section.

\begin{lemma}\label{lem: identities}
	For all $j=1,\dots,N-1$, the output of Algorithm~\ref{alg: full discr} satisfies
	\begin{align}
	    \vv_h^0 &= \d_t \mm_h^1,\label{eq: identity0}\\
		2\vv_h^j &= 3\d_t\mm_h^{j+1}-\d_t\mm_h^j;\label{eq: identity1}\\
		2(\vv_h^j-\d_t\mm_h^{j+1}) &= \tau \d_t^2\mm_h^{j+1};\label{eq: identity2}\\
		\vv_h^j &= \tau\d_t^2\mm_h^{j+1}+\dfrac{1}{2}[\d _t\mm_h^{j+1}+\d _t\mm_h^j];\label{eq: identity3}\\
		\dfrac{\tau}{3}(\vv_h^j-\d _t\mm_h^j) &= \dfrac{\tau^2}{2}\d _t^2\mm_h^{j+1};\label{eq: identity4}\\
	\mm_h^{j+1}-\widehat{\mm}_h^{j+1}&= \tau(\d _t\mm_h^{j+1}-\d _t\mm_h^j)=\tau^2\d _t^2\mm_h^{j+1}.\label{eq: identity5}
	\end{align}
\end{lemma}

\begin{proof}
The identity~\eqref{eq: identity0} is obvious from Algorithm~\ref{alg: full discr}\eqref{item: def mm1}.
The identities~\eqref{eq: identity1}--\eqref{eq: identity4} follow from the use of BDF2 in Algorithm~\ref{alg: full discr}(\ref{item: next steps}.\ref{item: def mmhj}), i.e.,
\begin{align*}
 \vv_h^j=\dfrac{1}{\tau}\bigg(\dfrac{3}{2}\mm_h^{j+1}-2\mm_h^j+\dfrac{1}{2}\mm_h^{j-1}\bigg)
 \quad \text{for all } j = 1,\dots,N.
\end{align*}
The final identity~\eqref{eq: identity5} follows from the definition of $\widehat{\mm}_h^{j+1}$ in Algorithm~\ref{alg: full discr}(\ref{item: next steps}.\ref{item: def mhat}).
\end{proof}
\begin{lemma}[Discrete $\boldsymbol{L^r}$-norm]\label{lem: norm equivalence 1}
 Let $\mathcal T _h$ be a $\gamma$-quasi-uniform mesh, $d=2,3$, and $~r~\in~[1,+\infty)$.
 Then, it holds that
 \begin{equation}\label{eq: norm equivalence 1}
  C^{-1}\norm{{v}_h}^r_{L^r(\Omega)}
  \leq 
  h^{d} \sum_{\vec{z}\in\mathcal N _h}|v_h(\vec{z})|^r
  \le 
  C\norm{v_h}^r_{L^r(\Omega)}
  \quad \text{for all }\,\, v_h\,\in \SSS^1(\mathcal T _h),
 \end{equation}
 where $C>0$ depends only on $r$ and the $\gamma$-quasi-uniformity of $\TT_h$.\qed
\end{lemma}
	
Moreover, we recall two results from \cite{abp2024} adapted to our notation.
We point out that the proof of~\eqref{eq: norm equivalence 2} in \cite[Lemma~2.1]{abp2024} relies only on the algebraic identity~\eqref{eq: identity1}, while the proof of~\eqref{eq: inverse estimate} in~\cite[Lemma~2.2]{abp2024} relies only on~\eqref{eq: norm equivalence 2}.
This allows to state the same result for both the $\bL^2$ and the $\bH^1$-norms.
For the convenience of the reader, we include the proofs of the following lemmas in the appendix.

\begin{lemma}[{Norm equivalence \cite[Lemma~2.1]{abp2024}}]\label{lem: norm equivalence 2}
For every $n = 1,\dots, N$ it holds 
\begin{equation}\label{eq: norm equivalence 2}
 C_1 \Big( \tau\sum_{j=0}^{n-1}\norm{\vv_h^j}^2_\Omega \Big)^{1/2}
 \le 
 \Big( \tau\sum_{j=0}^{n-1}\norm{\d_t\mm_h^{j+1}}^2_\Omega \Big)^{1/2}
 \le 
 C_2 \Big( \tau\sum_{j=0}^{n-1}\norm{\vv_h^j}^2_\Omega \Big)^{1/2}
\end{equation} 
and
\begin{equation}\label{eq: norm equivalence grad 2}
 C_1 \Big( \tau\sum_{j=0}^{n-1}\norm{\nabla\vv_h^j}^2_\Omega \Big)^{1/2}
 \le \Big( \tau\sum_{j=0}^{n-1}\norm[\big]{\nabla\d_t\mm_h^{j+1}}^2_\Omega \Big)^{1/2}
 \le C_2 \Big( \tau\sum_{j=0}^{n-1}\norm{\nabla\vv_h^j}^2_\Omega \Big)^{1/2},
\end{equation} 
where $C_1=\sqrt{1/5}$ and $C_2 = \sqrt{9/7}$.
\end{lemma}

\begin{lemma}[{Inverse estimate \cite[Lemma~2.2]{abp2024}}]\label{lem: inverse estimate}
For every $n = 1,\dots, N$, it holds 
\begin{equation}\label{eq: inverse estimate}
 \Big( \tau\sum_{j=2}^n\normO{\d_t^2\mm_h^j}^2 \Big) ^{1/2}
 \le
 \frac{1}{{\tau}}C_3 \Big( \tau\sum_{j=0}^{n-1}\normO{\vv_h^j}^2 \Big)^{1/2}
\end{equation}
and
\begin{equation}\label{eq: inverse estimate grad}
 \Big( \tau\sum_{j=2}^n\normO[\big]{\nabla\d_t^2\mm_h^j}^2 \Big) ^{1/2}
 \le
 \frac{1}{\tau} C_3 \Big( \tau\sum_{j=0}^{n-1}\normO{\nabla\vv_h^j}^2 \Big)^{1/2},
\end{equation}
where $C_3 = \sqrt{18/7}$.
\end{lemma}

\subsection{Discrete energy identity}

In this section, we derive a discrete energy identity and a discrete energy estimate for the approximations obtained by Algorithm~\ref{alg: full discr}.
The next result follows ideas from~\cite[Proposition~3.3]{feischl}.
However, we stress that
\cite{feischl} employs a different (less local) discrete tangent space and
that our proof also allows to control the second-order discrete time-derivative $\sum_{j=1}^{n-1}\normO{\nabla\d_t^2\mm_h^{j+1}}^2$.
While our proof relies only on direct algebraic computations, we refer to \cite{Dahlquist, baiocchi_crouzeix, nevanlinna_odeh} or \cite[Section V.6]{hairer} for more details on such energy estimates.


\begin{lemma}[Discrete energy identity]\label{lem: discrete energy identity}
	Let $1 \le n \le N$ and recall $\eta_0$ and $\eta_n$ from \eqref{eq: eta0n}. Then Algorithm \ref{alg: full discr} guarantees that
	\begin{equation}\label{eq: energy inequality 2}
		\begin{split}
		\alpha \tau \sum_{j=0}^{n-1}\normO{\vv_h^j}^2&+\frac{\lamex^2}{2}\normO{\nabla\mm_h^n}^2+\frac{\lamex^2}{2}\tau^2\normO{\d _t\nabla\mm_h^n}^2 + \frac{\lamex^2}{4}\tau^4 \sum_{j=1}^{n-1}\normO{\nabla\d _t^2\mm_h^{j+1}}^2 \\
		& =\frac{\lamex^2}{2}\normO{\nabla\mm_h^0}^2 +\tau\sum_{j=0}^{n-1}\scalarproductO{\ff_h^{j+1}}{\vv_h^j}+\frac{\lamex^2}{4}(\eta_0+\eta_n).
		\end{split}
	\end{equation}
Moreover, this implies the following discrete energy estimates:
	\begin{align}
	\label{eq: bound_3}
		\lamex^2\normO{\nabla \mm_h^1}^2 + \tau \alpha \normO{\vv_h^0}^2 +\tau^2\lamex^2\normO{\nabla\vv_h^0}^2
		&\leq \lamex^2 \normO{\nabla \mm_h^0}^2 + \dfrac{\tau}{\alpha} \normO{\ff_h^1}^2 \quad \text{for } n=1;
	\end{align}
	and
	\begin{align}\label{eq: energy_inequality}
\begin{split}
 \lamex^2 \gamma^- \big( \normO{\nabla \mm_h^n}^2 + \normO{\nabla \mm_h^{n-1}}^2 \big) 
 + \frac{\alpha\tau }{2} \sum_{j=1}^{n-1} \normO{\vv_h^{j}}^2 
 +\frac{\lamex^2}{4} \tau^4 \sum_{j=1}^{n-1} \normO{\nabla\d_t^2\mm_h^{j+1}}^2
 \\\leq
 \lamex^2 \gamma^+ \big( \normO{\nabla \mm_h^0}^2 + \normO{\nabla \mm_h^1}^2 \big) 
 + \frac{\tau}{2 \alpha} \sum_{j=1}^{n-1} \normO{\ff_h^{j+1}}^2 \quad \text{for } n\ge2,
\end{split}
\end{align}
where $\gamma^\pm=(3\pm2 \sqrt{2})/{4}$ are the two eigenvalues of the matrix $G\coloneqq \dfrac{1}{4}\bigg({\begin{array}{cc}
		1 & -2 \\
		-2 & 5 \\
	\end{array} } \bigg)$.
	
\end{lemma}

\begin{proof}
We will split the proof in four steps. Steps 1--3 concern the initialization step of Algorithm~\ref{alg: full discr}\eqref{eq: first step full discr item}--\eqref{item: def mm1}, while Steps 2--4 concern the subsequent BDF2 steps of Algorithm~\ref{alg: full discr}\eqref{item: next steps}.

{\bf Step 1 (Proof of (\ref{eq: energy inequality 2}) for $\boldsymbol{n=1}$).}
Choosing $\pphi_h=\vv_h^0$ in the variational formulation~\eqref{eq: first step full discr} of Algorithm~\ref{alg: full discr} leads to
\begin{equation*}
	\alpha \normO{\vv_h^0}^2 + \lamex^2 \scalarproductO{\nabla (\mm_h^0 + \tau \vv_h^0)}{\nabla \vv_h^0} = \scalarproductO{\ff_h^1}{\vv_h^0}.
\end{equation*}
Multiplying this by $\tau$ and using that $\vv_h^0=(\mm_h^1-\mm_h^0)/\tau$ by Algorithm~\ref{alg: full discr}\eqref{item: def mm1}, we see
\begin{equation*}
	\begin{split}
		&\tau \alpha \normO{\vv_h^0}^2 + \tau \lamex^2 \Big\langle{\nabla \Big( \frac{\mm_h^1 + \mm_h^0}{2} + \frac{\mm_h^1 - \mm_h^0}{2} \Big)},{\nabla \vv_h^0}\Big\rangle _\Omega = \tau \scalarproductO{\ff_h^1}{\vv_h^0}.
	\end{split}
\end{equation*}
This is equivalent to
\begin{align}\label{eq: bound_1}
		\dfrac{\lamex^2}{2}\normO{\nabla \mm_h^1}^2 + \dfrac{\tau^2}{2} \lamex^2\normO{\nabla \vv_h^0}^2 + \tau \alpha \normO{\vv_h^0}^2 
		&= \dfrac{\lamex^2}{2}\normO{\nabla \mm_h^0}^2 + \tau\scalarproductO{\ff_h^1}{\vv_h^0},
\end{align}
which is \eqref{eq: energy inequality 2} for $n=1$.

{\bf Step 2 (Proof of (\ref{eq: bound_3})).}
To prove~\eqref{eq: bound_3},
we bound $\scalarproductO{\ff_h^1}{\vv_h^0}$ by the Young inequality
\begin{equation*}
	\begin{split}
		\scalarproductO{\ff_h^1}{\vv_h^0} &\le \normO{\ff_h^1} \normO{\vv_h^0} \le \frac{1}{2\alpha} \normO{\ff_h^1}^2 +  \frac{\alpha}{2} \normO{\vv_h^0}^2
	\end{split}
\end{equation*}
and insert it into~\eqref{eq: bound_1}. 

{\bf Step 3 (Proof of (\ref{eq: energy inequality 2}) for $\boldsymbol{n \ge 2}$).} Due to the definition of $\mm_h^{j+1}$ in Algorithm~\ref{alg: full discr}(\ref{item: next steps}.\ref{item: def mmhj}), the discrete variational formulation~\eqref{eq: full discr alg B} is equivalent to
\[
 \alpha\scalarproductO{\vv_h^j}{\pphi_h} 
 + \scalarproductO{\widehat{\mm}_h^{j+1}\times\vv_h^j}{\pphi_h}
 + \lamex^2 \scalarproductO{\nabla\mm_h^{j+1}}{\nabla\pphi_h}
 =
 \scalarproductO{\ff_h^{j+1}}{\pphi_h}
 \,\,\text{for all}\,\,\pphi_h\in \bT_h(\widehat{\mm}_h^{j+1}).
\]	
We test this equation with $\pphi_h=\tau\vv_h^j\in \bT_h(\widehat{\mm}_h^{j+1})$ to obtain 
\begin{equation}\label{eq: energy_inequality_halfproof}
	\alpha\tau\normO{\vv_h^j}^2 +\lamex^2\tau\scalarproductO{\nabla\mm_h^{j+1}}{\nabla\vv_h^j} = \tau\scalarproductO{\ff_h^{j+1}}{\vv_h^j}.
\end{equation}
Let $1\le j \le N$. Recall that $\tau\vv^j_h \coloneqq \dfrac{3}{2}\mm^{j+1}_h-2\mm^{j}_h+\dfrac{1}{2}\mm^{j-1}_h$ by definition of $\mm_h^{j+1}$ in Algorithm~\ref{alg: full discr}(\ref{item: next steps}.\ref{item: def mmhj}).
Therefore, it holds that
\begin{align}
 \nonumber
 &\tau\langle\nabla \mm_h^{j+1}, \nabla \vv_h^j\rangle_\Omega 
 = 
 \frac{1}{4} \Big[ 6\normO{\nabla\mm_h^{j+1}}^2 - 8\scalarproductO{\nabla\mm_h^{j+1}}{\nabla\mm_h^{j}} 
 + 2\scalarproductO{\nabla\mm_h^{j+1}}{\nabla\mm_h^{j-1}} \Big] 
 \\& \quad
 \nonumber
 = 
 \frac{1}{4} \Big[ \normO{\nabla\mm_h^{j+1}}^2 - \normO{\nabla\mm_h^{j}}^2 + \normO{2\nabla\mm_h^{j+1}}^2
  + \normO{\nabla\mm_h^{j}}^2  - 2\scalarproductO{2\nabla\mm_h^{j+1}}{\nabla\mm_h^{j}}
 \\& \qquad \qquad  
 \nonumber
 - \normO{2\nabla\mm_h^{j}}^2 - \normO{\nabla\mm_h^{j-1}}^2 + 2\scalarproductO{2\nabla\mm_h^{j}}{\nabla\mm_h^{j-1}}\Big]
 \\& \qquad \quad
 \nonumber
 + \frac{1}{4} \Big[ \normO{\nabla\mm_h^{j+1}}^2 + \normO{2\nabla\mm_h^{j}}^2 + \normO{\nabla\mm_h^{j-1}}^2 + 2\scalarproductO{\nabla\mm_h^{j+1}}{\nabla\mm_h^{j-1}}
 \\& \qquad \qquad \qquad
 \nonumber
 - 2\scalarproductO{2\nabla\mm_h^{j}}{\nabla\mm_h^{j-1}} - 2\scalarproductO{2\nabla\mm_h^{j+1}}{\nabla\mm_h^{j}} \Big]\\
 \nonumber
 &\quad
 = 
 \frac{1}{4} \Big[ \normO{\nabla\mm_h^{j+1}}^2 - \normO{\nabla\mm_h^{j}}^2 + \normO{2\nabla\mm_h^{j+1} - \nabla\mm_h^{j}}^2 - \normO{2\nabla\mm_h^{j} - \nabla\mm_h^{j-1}}^2 \Big]
 \\& \qquad 
 \nonumber
 + \frac{1}{4} \normO{\nabla\big(\mm_h^{j+1} - 2\mm_h^{j} + \mm_h^{j-1})}^2
 \\& \quad
 \nonumber
 = \Big[\frac{1}{4}\normO{\nabla \mm_h^{j}}^2 - \scalarproductO{\nabla \mm_h^{j+1}}{\nabla \mm_h^{j}} + \frac{5}{4}\normO{\nabla \mm_h^{j+1}}^2\Big]
 \\& \qquad 
 - \Big[ \frac{1}{4}\normO{\nabla \mm_h^{j-1}}^2 - \scalarproductO{\nabla \mm_h^{j}}{\nabla \mm_h^{j-1}} + \frac{5}{4}\normO{\nabla \mm_h^{j}}^2 \Big]
 \label{eq: dot_relation_grad}
 + \frac{1}{4}\normO{\tau^2\nabla\d_t^2\mm_h^{j+1}}^2.
\end{align}
For any $1\le j \le n$, there holds the identity
\begin{align*}
  \frac{1}{4} \normO{\nabla \mm_h^{j}}^2-\scalarproductO{\nabla \mm_h^{j+1}}{\nabla \mm_h^{j}}+\frac{5}{4}\normO{\nabla \mm_h^{j+1}}^2
 =&
 \frac{1}{2} \normO{\nabla \mm_h^{j+1}-\nabla\mm_h^j}^2 + \frac{3}{4}\normO{\nabla\mm_h^{j+1}}^2- \frac{1}{4}\normO{\nabla\mm_h^{j}}^2
  \\= & \frac{1}{2} \normO{\nabla \mm_h^{j+1}}^2+
 \frac{\tau}{4} \d{}_t \normO{\nabla\mm_h^{j+1}}^2+\frac{\tau^2}{2} \normO{\d{}_t \nabla\mm_h^{j+1}}^2 .
\end{align*}
Combining this with \eqref{eq: energy_inequality_halfproof}--\eqref{eq: dot_relation_grad}, we obtain
\begin{equation*}
\begin{aligned}
	\alpha \tau\normO{\vv_h^{j}}^2&+\frac{\lamex^2}{2}\normO{\nabla\mm_h^{j+1}}^2+\frac{\lamex^2}{4}\tau \d _t\normO{\nabla\mm_h^{j+1}}^2 + \frac{\lamex^2}{2}\tau^2\normO{\d _t\nabla\mm_h^{j+1}}^2+\frac{\lamex^2}{4}\tau^4\normO{\nabla\d_t^2\mm_h^{j+1}}^2\\
	& = \frac{\lamex^2}{2}\normO{\nabla\mm_h^{j}}^2+\frac{\lamex^2}{4}\tau \d _t\normO{\nabla\mm_h^j}^2+\frac{\lamex^2}{2}\tau^2\normO{\d _t\nabla\mm_h^{j}}^2+\tau \scalarproductO{\ff_h^{j+1}}{\vv_h^{j}}.
\end{aligned}
\end{equation*}
Summing this over $j=1,\dots,n-1$, we get
	\begin{equation*}\label{eq: second energy inequality first step}
		\begin{split}
		\alpha\tau\sum_{j=1}^{n-1}\normO{\vv_h^j}^2&+\frac{\lamex^2}{2}\normO{\nabla\mm_h^n}^2+\frac{\lamex^2}{4}\tau \d _t\normO{\nabla\mm_h^n}^2 + \frac{\lamex^2}{2}\tau^2\normO{\d _t\nabla\mm_h^n}^2+\frac{\lamex^2}{4}\tau^4\sum_{j=1}^{n-1}\normO{\nabla\d_t^2\mm_h^{j+1}}^2\\
		& = \frac{\lamex^2}{2}\normO{\nabla\mm_h^1}^2+\frac{\lamex^2}{4}\tau\d_t\normO{\nabla\mm_h^1}^2+\frac{\lamex^2}{2}\tau^2\normO{\d_t\nabla\mm_h^1}^2+\tau\sum_{j=1}^{n-1}\scalarproductO{\ff_h^{j+1}}{\vv_h^j}.
		\end{split}
	\end{equation*}
	Adding~\eqref{eq: bound_1} to include also $j=0$, the definition of $\eta_0$ and $\eta_n$ in \eqref{eq: eta0n} concludes the proof.

{\bf Step 4: (Proof of (\ref{eq: energy_inequality})).}
Recalling the definition of the matrix $G$, we can write
\begin{equation*}
\begin{aligned}
 \tau \scalarproductO{\nabla \mm_h^{j+1}}{\nabla \vv_h^j} &=\big[g_{11}\normO{\nabla \mm_h^{j}}^2 + 
 2 g_{12}\scalarproductO{\nabla \mm_h^{j+1}}{\nabla \mm_h^{j}} + 
 g_{22}\normO{\nabla \mm_h^{j+1}}^2\big]
 \\&
 \mkern+50mu\big[g_{11}\normO{\nabla \mm_h^{j-1}}^2 + 
 2 g_{12}\scalarproductO{\nabla \mm_h^{j}}{\nabla \mm_h^{j-1}} + 
 g_{22}\normO{\nabla \mm_h^{j}}^2\big]\\
 &\mkern+50mu+ \frac{1}{4}\normO{\tau^2\nabla\d_t^2\mm_h^{j+1}}^2.
\end{aligned}
\end{equation*}			
		Together with~\eqref{eq: energy_inequality_halfproof}, we obtain
	\begin{align}
		& \alpha \tau\normO{\vv_h^{j}}^2+\lamex^2\big[g_{11}\normO{\nabla \mm_h^{j}}^2+2 g_{12}\scalarproductO{\nabla \mm_h^{j+1}}{\nabla \mm_h^{j}}+g_{22}\normO{\nabla \mm_h^{j+1}}^2\big] +\frac{\lamex^2}{4}\tau^4\normO{\nabla\d_t^2\mm_h^{j+1}}^2 \nonumber\\
		& \qquad= \lamex^2\big[g_{11}\normO{\nabla \mm_h^{j-1}}^2+2 g_{12}\scalarproductO{\nabla \mm_h^{j}}{\nabla \mm_h^{j-1}}+g_{22}\normO{\nabla \mm_h^{j}}^2\big]+\tau\scalarproductO{\ff_h^{j+1}}{\vv_h^j}. \label{eq: energy_inequality_halfproof2}
	\end{align}
	
	Using the Young inequality $	\scalarproductO{\ff_h^{j+1}}{\vv_h^j} 
 		\le \dfrac{\alpha}{2} \normO{\vv_h^j}^2 + \dfrac{1}{2\alpha} \normO{\ff_h^{j+1}}^2$
	and summing over $j=1,\dots, n-1$, we obtain
		\[
		\begin{aligned}
		 \frac{\alpha \tau}{2} \sum_{j=1}^{n-1}\normO{\vv_h^j}^2+&\lamex^2\big[g_{11}\normO{\nabla \mm_h^{n-1}}^2+2 g_{12}\scalarproductO{\nabla \mm_h^n}{\nabla \mm_h^{n-1}}+g_{22}\normO{\nabla \mm_h^n}^2\big]+\frac{\lamex^2}{4}\tau^4\sum_{j=1}^{n-1}\normO{\nabla\d_t^2\mm_h^{j+1}}^2\\
		& \le \lamex^2\big[g_{11}\normO{\nabla \mm_h^0}^2+2 g_{12}\scalarproductO{\nabla \mm_h^1}{\nabla \mm_h^0}+g_{22}\normO{\nabla \mm_h^1}^2\big]+\frac{\tau}{2 \alpha} \sum_{j=1}^{n-1}\normO{\ff_h^{j+1}}^2.
		\end{aligned}
		\]
		Since $\gamma^{-}|x|^2 \leq G x \cdot x \leq \gamma^{+}|x|^2$ for every $x \in \mathbb{R}^2$, it follows that
		\[
		\begin{aligned}
		g_{11}\normO{\nabla \mm_h^{n-1}}^2+2 g_{12}\scalarproductO{\nabla \mm_h^n}{\nabla \mm_h^{n-1}}+g_{22}\normO{\nabla \mm_h^n}^2 & \geq \gamma^-\big(\normO{\nabla \mm_h^n}^2+\normO{\nabla \mm_h^{n-1}}^2\big)
		\end{aligned}
		\]
		and
		\[
		\begin{aligned}
			&g_{11}\normO{\nabla \mm_h^{0}}^2+2 g_{12}\scalarproductO{\nabla \mm_h^0}{\nabla \mm_h^{1}}+g_{22}\normO{\nabla \mm_h^1}^2 \leq \gamma^+\big(\normO{\nabla \mm_h^0}^2+\normO{\nabla \mm_h^{1}}^2\big).
		\end{aligned}
		\]
		Combining the last three formulae, we obtain~\eqref{eq: energy_inequality} and conclude the proof.
\end{proof}
\begin{remark}
		Note that the discrete energy identity \eqref{eq: energy inequality 2} contains a physical dissipation term $\tau \sum_{j=0}^{n-1}\normO{\vv_h^j}^2$, which converges in the limit $(h,\tau)\to(0,0)$ to the physical dissipation term $\int_0^{t_n}\normO{\partial_t\mm}^2\d t$ of the continuous energy inequality of Definition~\ref{def: weak solution}\eqref{def: weak solution iv}. Moreover, the terms $\normO{\nabla\mm_h^n}^2$, $\normO{\nabla\mm_h^0}^2$ and $\tau\sum_{j=0}^{n-1}\scalarproductO{\ff_h^{j+1}}{\vv_h^j}$ contribute to the energy terms at times $t'$ and $t_0$ and to the integral including the applied field, i.e., first, third and fourth term of Definition~\ref{def: weak solution}\eqref{def: weak solution iv}. The additional terms $\tau^2\normO{\d_t\nabla\mm_h^n}^2$, $\tau^4\sum_{j=1}^{n-1}\normO{\nabla\d_t^2\mm_h^{j+1}}^2$ and $\eta_0+\eta_n$ in \eqref{eq: energy inequality 2} can be interpreted as numerical dissipation terms, coming from the BDF2 method, which vanish in the limit $(h,\tau)\to(0,0)$. We refer to the final part of the proof of Theorem~\ref{thm: main result} for more details on these convergences to the continuous energy inequality.

		The analogous energy identity in the case of the projection-free first-order tangent plane scheme reads 
		$$
			\frac{\lamex^2}{2}\normO{\nabla\mm_h^n}^2 + \alpha \tau \sum_{j=0}^{n-1}\normO{\vv_h^j}^2 +\frac{\lamex^2}{2} \tau^2\sum_{j=0}^{n-1}\normO{\nabla\vv_h^j}^2 = \frac{\lamex^2}{2}\normO{\nabla\mm_h^0}^2 + \tau \sum_{j=0}^{n-1}\scalarproductO{\ff_h^{j}}{\vv_h^j}.
		$$
		As before, this identity involves terms that converge to their continuous counterpart and a different numerical dissipation term $\tau^2\sum_{j=0}^{n-1}\normO{\nabla\vv_h^j}^2$, which vanishes in the limit $(h,\tau)\to(0,0)$; see e.g., \cite{Spin-polarized} or \cite{hpprss2019} for details.
	\end{remark}

\subsection{Boundedness of discrete functions}

Lemma~\ref{lem: discrete energy identity} and assumption~\eqref{enum: punto A2}--\eqref{enum: punto A1} allow us to prove the following proposition.

\begin{proposition}[Boundedness of discrete gradients]\label{prop: bounds}
For every $0\le n\le N$, Algorithm \textnormal{\ref{alg: full discr}} guarantees
\begin{equation}\label{eq: uniform bound}
 \normO{\nabla \mm_h^n}^2 + \tau \sum_{j=0}^{n-1}\normO{\vv_h^j}^2
 + \tau^4 \sum_{j=1}^{n-1} \normO{\nabla\d_t^2\mm_h^{j+1}}^2
 \leq 
 C_{0,f},
\end{equation}
where $C_{0,f}$ depends only on $\gamma^\pm, \lamex, \alpha, C_0$, and $C_f$, but is independent of $h$, $\tau$, and $n$.
\end{proposition}

\begin{proof}
By assumption, it holds $\normO{\nabla \mm_h^0}^2\leq C_0$ and $\tau \normO{\ff_h^1}^2 \le \normOT{\ff_{h\tau}^+}^2\le C_f$.
For $n=1$, the 
bound~\eqref{eq: bound_3} from Lemma~\ref{lem: discrete energy identity} thus implies 
that $\normO{\nabla \mm_h^0}^2 + \normO{\nabla \mm_h^1}^2 + \tau\normO{\vv_h^0}^2 \lesssim 1$. Therefore, let $n\in\{2,\dots,N\}$. With~\eqref{eq: energy_inequality} from Lemma~\ref{lem: discrete energy identity},
this leads to
\begin{align*}
&
\normO{\nabla \mm_h^n}^2 + \tau \sum_{j=1}^{n-1}\normO{\vv_h^j}^2 + \tau^4 \sum_{j=1}^{n-1} \normO{\nabla\d_t^2\mm_h^{j+1}}^2
\\& \quad
\stackrel{\eqref{eq: energy_inequality}}{\lesssim}
\normO{\nabla \mm_h^0}^2 + \normO{\nabla \mm_h^1}^2 + \tau \sum_{j=1}^{n-1}\normO{\ff_h^{j+1}}^2
\lesssim
1 + \tau \sum_{j=1}^n\normO{\ff_h^j}^2
\le 
1 + \normOT{\ff_{h\tau}^+}^2
\lesssim
1.
\end{align*}
This concludes the proof.
\end{proof}

The next result provides constraint violation estimates and proves, in particular, that $\normO{\mm_h^n}$ is uniformly bounded. We note that its proof adapts \cite[Proposition~3.4]{abp2024}, which guarantees constraint violation estimates of the numerical approximation of harmonic maps into spheres.

\begin{proposition}[Discrete modulus constraint]\label{prop: constraint violation}
	For every $n\ge 2$ and every $\vec{z}\in \NN_h$, we have that
	\begin{equation}\label{eq: first constraint relation}
		|\mm_h^n(\vec{z})|^2=1+\frac{3}{2}\Big(1-\frac{1}{3^n}\Big)\tau^2 |\vv_h^0(\vec{z})|^2 + \frac{3}{2}\tau^4\sum_{i=2}^n\Big(1-\frac{1}{3^{n+1-i}}\Big)|\d_t^2\mm_h^i(\vec{z})|^2.
	\end{equation}
	Moreover, it holds that 
	\begin{equation}\label{eq: constraint violation}
	 \norm[\big]{|\mm_h^n|^2-1}_{\mathbf{L}^1(\Omega)}\le C\tau,
	\end{equation}
	where the constant $C>0$ depends on $\alpha$, $\gamma^+$, $\lamex$, $C_0$, $C_f$, and the $\gamma$-quasi-uniformity of $\TT_h$.
	For every $n\ge 1$, it thus follows that 
	\begin{equation}\label{eq2: constraint violation}
	 \normO{\mm_h^n}^2\le C\tau+|\Omega|.
	\end{equation} 
	\end{proposition}

\begin{proof}
Arguing as in the proof of~\eqref{eq: dot_relation_grad} in Lemma~\ref{lem: discrete energy identity}, but with $\mm_h^{j+1}(\vec{z})\cdot\vv_h^{j}(\vec{z})$ instead of $\scalarproductO{\nabla\mm_h^{j+1}}{\nabla\vv_h^{j}}$, we show that, for all $\vec{z}\in\NN_h$ and all $j = 1,\dots,n-1$,
\begin{equation}\label{eq: nodes useful relation}
\begin{split}
 \mm_h^{j+1}(\vec{z})\cdot\vv_h^{j}(\vec{z}) &= \frac{1}{\tau} \Big(\frac{1}{4}|\mm_h^{j}(\vec{z})|^2 - \mm_h^{j+1}(\vec{z})\cdot\mm_h^{j}(\vec{z}) + \frac{5}{4}|\mm_h^{j+1}(\vec{z})|^2\Big) \\
 &\mkern-80mu- \frac{1}{\tau} \Big(\frac{1}{4}|\mm_h^{j-1}(\vec{z})|^2 - \mm_h^{j}(\vec{z})\cdot\mm_h^{j-1}(\vec{z}) + \frac{5}{4}|\mm_h^{j}(\vec{z})|^2\Big) + \frac{\tau^3}{4}|\d _t^2 \mm_h^{j+1}(\vec{z})|^2
\end{split}
\end{equation}
With $\vv_h^{j} \in \boldsymbol{T}_h(\widehat{\mm}_h^{j+1})$ and hence $\widehat{\mm}_h^{j+1}(\vec{z}) \cdot \vv_h^{j}(\vec{z})=0$, this leads to
	\begin{align}
		\frac{1}{2}|\mm_h^{j}(\vec{z})|^2-2& \mm_h^{j+1}(\vec{z}) \cdot \mm_h^{j}(\vec{z})+\frac{5}{2}|\mm_h^{j+1}(\vec{z})|^2\nonumber\\&-\frac{1}{2}|\mm_h^{j-1}(\vec{z})|^2+2 \mm_h^{j}(\vec{z}) \cdot \mm_h^{j-1}(\vec{z})-\frac{5}{2}|\mm_h^{j}(\vec{z})|^2 \nonumber\\ 
		&\mkern-50mu \overset{\eqref{eq: nodes useful relation}}{=}2 \tau \mm_h^{j+1}(\vec{z}) \cdot \vv_h^{j}(\vec{z})-\frac{\tau^4}{2}|\d_t^2 \mm_h^{j+1}(\vec{z})|^2 \nonumber\\ 
		& \mkern-50mu\overset{\eqref{eq: identity5}}{=}2 \tau\big(\widehat{\mm}_h^{j+1}(\vec{z})+\tau^2 \d _t^2 \mm _h ^{j+1}(\vec{z})\big) \cdot \vv_h^{j}(\vec{z})-\frac{\tau^4}{2}|\d_t^2 \mm_h^{j+1}(\vec{z})|^2\nonumber\\
		& \mkern-50mu \,\,= 2 \tau^3 \d _t ^2 \mm_h^{j+1}(\vec{z}) \cdot \vv_h^{j}(\vec{z})-\frac{\tau^4}{2}|\d _t ^2 \mm_h^{j+1}(\vec{z})|^2\nonumber\\ 
		& \mkern-50mu\overset{\eqref{eq: identity3}}{\,=\,}2 \tau^4 |\d _t ^2 \mm_h^{j+1}(\vec{z})|^2+\tau^3\d_t^2\mm_h^{j+1}(\vec{z})\cdot[\d _t\mm_h^{j+1}(\vec{z})+\d _t\mm_h^j(\vec{z})]-\frac{\tau^4}{2}|\d_t^2\mm_h^{j+1}(\vec{z})|^2\nonumber\\
		& \mkern-50mu\overset{\eqref{eq: second order discrete time derivative}}{\,=\,}\frac{3}{2} \tau^4|\d _t ^2 \mm_h^{j+1}(\vec{z})|^2+\tau^2\big(|\d _t \mm_h^{j+1}(\vec{z})|^2-|\d _t\mm_h^{j}(\vec{z})|^2\big)\label{eq: last line}.
	\end{align}
Moreover, the left-hand side of~\eqref{eq: last line} satisfies
	\begin{align}
		&\frac{1}{2}|\mm_h^{j}(\vec{z})|^2-2\mm_h^{j+1}(\vec{z}) \cdot \mm_h^{j}(\vec{z})+\frac{5}{2}|\mm_h^{j+1}(\vec{z})|^2\nonumber\\
		&\qquad-\frac{1}{2}|\mm_h^{j-1}(\vec{z})|^2+2 \mm_h^{j}(\vec{z}) \cdot \mm_h^{j-1}(\vec{z})-\frac{5}{2}|\mm_h^{j}(\vec{z})|^2 \nonumber\\
		&=|\mm_h^{j+1}(\vec{z})-\mm_h^{j}(\vec{z})|^2-\frac{1}{2}|\mm_h^{j}(\vec{z})|^2+\frac{3}{2}|\mm_h^{j+1}(\vec{z})|^2\nonumber\\
		&\qquad-\Big[|\mm_h^{j}(\vec{z})-\mm_h^{j-1}(\vec{z})|^2-\frac{1}{2}|\mm_h^{j-1}(\vec{z})|^2+\frac{3}{2}|\mm_h^{j}(\vec{z})|^2\Big]\nonumber\\
		&=\tau^2\big[|\d _t \mm_h^{j+1}(\vec{z})|^2-|\d _t\mm_h^{j}(\vec{z})|^2\big]
		-2 |\mm_h^{j}(\vec{z})|^2+\frac{3}{2}|\mm_h^{j+1}(\vec{z})|^2
		+\frac{1}{2}|\mm_h^{j-1}(\vec{z})|^2
		.\label{eq: last line 2}
	\end{align}	
	The equality of the right-hand sides of~\eqref{eq: last line}--\eqref{eq: last line 2} shows that 
	\begin{equation}\label{eq: relazione intermedia}
		\frac{3}{2}|\mm_h^{j+1}(\vec{z})|^2-2|\mm_h^{j}(\vec{z})|^2+\frac{1}{2}|\mm_h^{j-1}(\vec{z})|^2=\frac{3}{2}\tau^4|\d_t^2 \mm_h^{j+1}(\vec{z})|^2.
	\end{equation}
	We now interpret the final equation as an inhomogeneous linear difference equation in the variable \( |\mm_h^{j+1}(\vec{z})|^2 \) with constant coefficients. 
	
	In order to have a simpler notation, let $a_{j+1}\coloneqq|\mm_h^{j+1}(\vec{z})|^2$ and $b_{j+1}\coloneqq \frac{3}{2}\tau^4|\d_t^2\mm_h^{j+1}(\vec{z})|^2$. Then, \eqref{eq: relazione intermedia} can be rewritten as
	\begin{equation*}
		\frac{3}{2}a_{j+1}-2a_j+\frac{1}{2}a_{j-1}=b_{j+1}\quad \text{for every} \quad j\ge 1.
	\end{equation*}
	Let $\lambda$ with $|\lambda|<1$. It is possible to rewrite this relation in an algebraic form multiplying it by $\lambda^{j+1}$ and summing over $j\ge 1$ as
	\begin{equation*}
		\frac{3}{2} \sum_{j\ge 1} a_{j+1}\lambda^{j+1} - 2 \sum_{j\ge 1} a_j \lambda^{j+1} + \frac{1}{2} \sum_{j\ge 1} a_{j-1}\lambda^{j+1} = \sum_{j\ge 1} b_{j+1}\lambda^{j+1}.
	\end{equation*}
	Defining 
	\begin{equation}\label{eq: recursion terms definition}
	P(\lambda) \coloneqq \frac{3}{2} - 2 \lambda + \frac{1}{2} \lambda^2 = \frac{(1-\lambda)(3-\lambda)}{2},\qquad A(\lambda)\coloneqq \sum_{j\ge 0} a_j\lambda^j, \quad\text{and}\quad B(\lambda)\coloneqq \sum_{j\ge 2} b_j \lambda^j,
	\end{equation}
	we can rearrange the previous relation as
	\begin{equation}\label{eq: algebraic relation}
		P(\lambda) A(\lambda) =  \Big(\frac{3}{2} a_1 - 2 a_0\Big)\lambda + \frac{3}{2} a_0 + B(\lambda).
	\end{equation}
	In order to find the coefficients \( a_j \), we need to invert the polynomial \( P(\lambda) \). Since $P(\lambda)$ has no roots within the open unit disk, the rational function \( {1}/{P(\lambda)} \) is holomorphic there. Thus, we can write it as 
	\begin{align*}
		\frac{1}{P(\lambda)}&=\sum_{k\ge 0} \gamma_k \lambda^k.
	\end{align*}
	Alternatively, we can use partial fraction decomposition to obtain
	\begin{equation}
		\frac{1}{P(\lambda)} = \frac{2}{(1-\lambda)(3-\lambda)} = \frac{1}{1-\lambda} - \frac{1}{3-\lambda} = \sum_{k\ge0} \big(1 - 3^{-(k+1)}\big) \lambda^k.
	\end{equation}
	This allows us to identify the coefficients as $\gamma_k = 1 - 3^{-(k+1)}$ for every $k\ge 0$.
	Therefore,~\eqref{eq: algebraic relation} reads as
	\begin{equation}
		A(\lambda) \stackrel{\eqref{eq: recursion terms definition}}{=} \sum_{j\ge 0} a_j \lambda^j \stackrel{\eqref{eq: algebraic relation}}{=} \Big(\sum_{j\ge 0} \gamma_j \lambda^j\Big) \Big(\Big(\frac{3}{2} a_1 - 2 a_0\Big)\lambda + \frac{3}{2} a_0 + \sum_{j\ge 2} b_j \lambda^j\Big).
	\end{equation}
	Hence, the generic coefficient \( a_n \) satisfies, for every \( n\ge 2 \),
	\begin{equation*}
		a_n = \gamma_{n-1}\Big(\frac{3}{2} a_1 - 2 a_0\Big) + \frac{3}{2} \gamma_n a_0 + \sum_{k=2}^n \gamma_{n-k} b_k.
	\end{equation*}
	Recalling the definition of \( a_j = |\mm_h^j(\vec{z})|^2\) and \( b_j = \frac{3}{2}\tau^4 |\d _t^2 \mm_h^{j}(\vec{z})|^2 \) and noticing that $|\mm_h^1(\vec{z})|^2 = 1+\tau^2 |\vv_h^0(\vec{z})|^2\ge 1 = |\mm_h^0(\vec{z})|^2 = a_0$, this can be rewritten as
	\begin{equation}
		|\mm_h^n(\vec{z})|^2 = \gamma_{n-1}\Big(\frac{3}{2}\tau^2 |\vv_h^0(\vec{z})|^2 - \frac{1}{2}\Big)+ \frac{3}{2} \gamma_n +\frac{3}{2} \tau^4  \sum_{k=2}^n \gamma_{n-k} |\d_t^2 \mm_h^k(\vec{z})|^2.
	\end{equation}
	The definition of the coefficients \( \gamma_k \) allows us to write this expression as \eqref{eq: first constraint relation}.
	To prove \eqref{eq: constraint violation}, we note that the right-hand side of \eqref{eq: first constraint relation} is a sum of non-negative terms. Therefore, $|\mm_h^n(\vec{z})|$ is increasing with $n$ (with $|\mm_h^1(\vec{z})|\ge 1$) and it holds that
		\begin{align}
			\norm{|\mm_h^n|^2-1}_{\mathbf{L}^1(\Omega)}&\stackrel{\eqref{eq: norm equivalence 1}}{\lesssim} h^d \sum_{\vec{z}\in\NN_h} \big||\mm_h^n(\vec{z})|^2-1\big| = h^d \sum_{\vec{z}\in\NN_h} \big(|\mm_h^n(\vec{z})|^2-1\big)\nonumber\\
			&\mkern-30mu= h^d \frac{3}{2}\big(1-3^{-n}\big)\tau^2\sum_{\vec{z}\in\NN_h} |\vv_h^0(\vec{z})|^2 + \frac{3}{2}\tau^4\sum_{i=2}^n\big(1-3^{-(n+1-i)}\big)\sum_{\vec{z}\in\NN_h}|\d_t^2\mm_h^i(\vec{z})|^2\nonumber\\
			&\mkern-30mu\stackrel{\eqref{eq: norm equivalence 1}}{\lesssim} \frac{3}{2}\big(1-3^{-n}\big)\tau^2\normO{\vv_h^0}^2 + \frac{3}{2}\tau^4\sum_{i=2}^n\big(1-3^{-(n+1-i)}\big)\normO{\d_t^2\mm_h^i}^2 \label{eq: discrete_condition_intermediate}\\
			&\mkern-30mu\stackrel{\eqref{eq: inverse estimate}}{\lesssim} \tau^2 \normO{\vv_h^0}^2+\tau^2 \sum_{i=0}^{n-1}\normO{\vv_h^i}^2\stackrel{\eqref{eq: uniform bound}}{\lesssim}\tau.\nonumber
		\end{align}
	This concludes the proof of \eqref{eq: constraint violation}. To prove \eqref{eq2: constraint violation}, note that
	$$
		\normO{\mm_h^n}^2 = \norm{|\mm_h^n|^2}_{L^1(\Omega)} \leq \norm{|\mm_h^n|^2-1}_{L^1(\Omega)} + |\Omega| \leq C\tau + |\Omega|.
	$$
	This concludes the proof.
\end{proof}
	
\begin{remark}\label{remark: second-order convergence rate}
		As in \cite{abp2024}, it follows from \eqref{eq: discrete_condition_intermediate} under the mild discrete regularity condition
		$$
		\normO{\vv_h^0}^2 + \tau^2 \sum_{j=2}^{n} \normO{\d_t^2\mm_h^j}^2 \lesssim 1,
		$$
		that \eqref{eq: constraint violation} holds in the improved form 
		$$
		\norm[\big]{|\mm_h^n|^2-1}_{\mathbf{L}^1(\Omega)}\le C\tau^2.
		$$
\end{remark}%
With Proposition~\ref{prop: bounds}--\ref{prop: constraint violation} at hand, 
we can now prove that also the time interpolations defined in~\eqref{eq:mhtau}--\eqref{eq:vhtau-minus} are uniformly bounded. 

\begin{proposition}[Uniform boundedness of interpolands]
	Under the assumptions of Theorem~\ref{thm: main result}, it holds that:
	\begin{equation}\label{eq: uniform bound S-T norms}
		\norm{\mm_{h\tau}}_{\bH^1(\OT)}+\norm{\mm_{h\tau}^{*}}_{L^\infty(0,T;\bH^1(\Omega))}+\norm{\vv_{h\tau}^-}_{\mathbf{L}^2(\OT)}\le C,
	\end{equation}
where $\mm_{h\tau}^*\in\{\mm_{h\tau},\mm_{h\tau}^{\pm},\mhatp\}$ and $C$ depends only on $\gamma^\pm,\lamex,\alpha,C_0,C_f,T$ and the $\gamma$-quasi-uniformity of $\TT_h$.
\end{proposition}
\begin{proof}
Proposition~\ref{prop: constraint violation} guarantees that 
\begin{align*}
	\normOT{\mm_{h\tau}}^2&= \int_{0}^{T}\normO{\mm_{h\tau}(t)}^2\d t= \sum_{j=0}^{N-1}\int_{t_j}^{t_{j+1}}\normO[\Big]{\frac{t-t_j}{\tau}\mm_h^{j+1}+\frac{t_{j+1}-t}{\tau}\mm_h^{j}}^2\d t \\
	&\le 2\sum_{j=0}^{N-1}\int_{t_j}^{t_{j+1}}\Big(\frac{t-t_j}{\tau}\Big)^2\normO{\mm_h^{j+1}}^2+\Big(\frac{t_{j+1}-t}{\tau}\Big)^2\normO{\mm_h^{j}}^2\d t \\
	&= \frac{2\tau}{3}\sum_{j=0}^{N-1}\big[\normO{\mm_h^{j+1}}^2+\normO{\mm_h^{j}}^2\big]
	\eqreff{eq2: constraint violation}\lesssim \tau N(\tau+|\Omega|)\simeq \tau+1.
\end{align*}
Similarly, Proposition~\ref{prop: bounds} guarantees that
\begin{align*}
    \normOT{\nabla\mm_{h\tau}}^2 
    &\le 2\sum_{j=0}^{N-1}\int_{t_j}^{t_{j+1}}\Big(\frac{t-t_j}{\tau}\Big)^2\normO{\nabla\mm_h^{j+1}}^2+\Big(\frac{t_{j+1}-t}{\tau}\Big)^2\normO{\nabla\mm_h^{j}}^2\d t \\
    &= \frac{2\tau}{3}\sum_{j=0}^{N-1}\big[\normO{\nabla\mm_h^{j+1}}^2+\normO{\nabla\mm_h^{j}}^2\big]\stackrel{\eqref{eq: uniform bound}}{\lesssim} \tau N= T.
\end{align*}
With $\partial_t\mm_{h\tau}(t)=\d _t\mm_h^{j+1}$ for every $t\in[t_j,t_{j+1})$, the norm equivalence~\eqref{eq: norm equivalence 2} from Lemma~\ref{lem: norm equivalence 2} and~\eqref{eq: uniform bound}, we see that 
	\begin{equation}
		\begin{split}
		\normOT{\partial_t\mm_{h\tau}}^2 
		=\sum_{j=0}^{N-1}\int_{t_j}^{t_{j+1}}\!\!\normO{\partial_t\mm_{h\tau}(t)}^2\d t
		=\tau\sum_{j=0}^{N-1}\normO{\d _t\mm_h^{j+1}}^2
		\stackrel{\eqref{eq: norm equivalence 2}}{\simeq} \tau\sum_{j=0}^{N-1}\normO{\vv_h^j}^2\stackrel{\eqref{eq: uniform bound}}{\lesssim}1.
		\end{split}
		\label{eq: uniform bound partial_t}
	\end{equation}
	Together, the last three formulas show $\norm{\mm_{h\tau}}_{\bH^1(\OT)} \lesssim 1$. 
	Analogously, it follows that
\begin{align*}
	\normOT{\vv_{h\tau}^-}^2 &= \int_{0}^{T}\normO{\vv_{h\tau}^-(t)}^2\d t = \sum_{j=0}^{N-1}\int_{t_j}^{t_{j+1}}\normO{\vv_h^j}^2\d t = \tau\sum_{j=0}^{N-1}\normO{\vv_h^j}^2\stackrel{\eqref{eq: uniform bound}}{\lesssim}1.
\end{align*}
The bounds $\norm{\mm_{h\tau}^\pm}_{L^\infty(0,T;\bH^1(\Omega))} \lesssim 1$ follow immediately from Proposition~\ref{prop: bounds} and Proposition~\ref{prop: constraint violation}. Thanks to the triangle inequality, the bounds on $\normO{\mm_h^j}$ and $\normO{\nabla\mm_h^j}$ assure a uniform bound also on the 
$\bL^2$-norm of $\widehat{\mm}_h^{j+1}=2\mm_h^{j}-\mm_h^{j-1}$ and of $\nabla\widehat{\mm}_h^{j+1}$. This yields
$\norm{\widehat{\mm}^+_{h\tau}}_{L^\infty(0,T;\bH^1(\Omega))} \lesssim 1$ and concludes the proof.
\end{proof}
In general, the sequence $\vv_{h\tau}^-$ is not uniformly bounded in $L^2(0, T ; \bH^1(\Omega))$. Nevertheless, the following proposition proves that at least $\tau\normOT{\nabla\vv_{h\tau}^-}$ can be controlled.
	
	\begin{proposition}\label{prop: grad v h tau to 0}
		It holds that
\begin{equation}\label{eq: grad v h tau to 0 secondo} 
 \lim_{(h,\tau) \to (0,0)}\tau \normOT{\nabla\partial_t\mm_{h\tau}}=0
 \qquad\text{and}\qquad
 \lim_{(h,\tau) \to (0,0)}\tau \normOT{\nabla\vv_{h\tau}^-}=0.
\end{equation}
	\end{proposition}
	
	\begin{proof}
		Let $n\le N$ be arbitrary and $1\le j \le n$. We test~\eqref{eq: full discr alg B} with $\pphi_h=\vv_h^j$ and multiply it by $\tau$ to obtain
		\begin{equation}\label{eq: discrete weak formulation}
			\begin{split}
			\alpha\tau\normO{{\vv_h^j}}^2&+\lamex^2\big\langle{\nabla\Big(\frac{4}{3}\mm_h^j-\frac{1}{3}\mm_h^{j-1}+\frac{2}{3}\tau\vv_h^j\Big)},{\tau\nabla\vv_h^j}\big\rangle_\Omega=\tau\scalarproductO{\ff_h^{j+1}}{\vv_h^j}.
			\end{split}
		\end{equation}
	We recall~\eqref{eq: identity1} and~\eqref{eq: identity4} from Lemma~\ref{lem: identities} and notice that
	\begin{equation}\label{eq: deriv spezzata}
		\mm_h^{j+1}=\frac{4}{3}\mm_h^j-\frac{1}{3}\mm_h^{j-1}+\frac{2}{3}\tau\vv_h^j=\mm_h^j+\frac{1}{3}\tau\d_t\mm_h^j+\frac{2}{3}\tau \vv_h^j.
	\end{equation}
	The relation $\langle a+b,a \rangle = \dfrac{1}{2}\big[\norm{a}^2+\norm{a+b}^2-\norm{b}^2\big]$ with $a=\tau\nabla\Big(\dfrac{1}{3}\d_t\mm_h^{j}+\dfrac{2}{3}\vv_h^j\Big)\stackrel{\eqref{eq: identity1}}{=}\tau\nabla\d_t\mm_h^{j+1}$ and $b=\nabla\mm_h^j$ gives
		\begin{equation*}
			\begin{split}
		&\big\langle{\nabla\Big(\frac{4}{3}\mm_h^j-\frac{1}{3}\mm_h^{j-1}+\frac{2}{3}\tau\vv_h^j\Big)},{\tau\nabla\vv_h^j}\big\rangle_\Omega\stackrel{\eqref{eq: deriv spezzata}}{=}\big\langle{\nabla\Big(\mm_h^j+\frac{1}{3}\tau\d_t\mm_h^j+\frac{2}{3}\tau \vv_h^j\Big)},{\tau\nabla\vv_h^j}\big\rangle_\Omega\\
		&\mkern+50mu=\big\langle{\nabla\Big(\mm_h^j+\frac{1}{3}\tau\d_t\mm_h^j+\frac{2}{3}\tau \vv_h^j\Big)},{\tau\nabla\Big(\frac{1}{3}\d_t\mm_h^j+\frac{2}{3}\vv_h^j\Big)}\big\rangle_\Omega\\
		&\mkern+120mu-\tau\big\langle{\nabla\Big(\mm_h^j+\frac{1}{3}\tau\d_t\mm_h^j+\frac{2}{3}\tau \vv_h^j\Big)},{\frac{1}{3}\nabla\big(\d_t\mm_h^j-\vv_h^j\big)}\big\rangle_\Omega\\
		&\mkern+50mu\stackrel{\eqref{eq: identity1},\eqref{eq: identity4}}{=}\frac{1}{2}\big[\normO{\tau\nabla\d _t\mm_h^{j+1}}^2+\normO{\nabla\mm_h^{j+1}}^2-\normO{\nabla\mm_h^j}^2\big]+\frac{\tau^2}{2}\big\langle{\nabla\mm_h^{j+1}},{\nabla\d_t^2\mm_h^{j+1}}\big\rangle_\Omega.
			\end{split}
		\end{equation*}
		We insert this relation in~\eqref{eq: discrete weak formulation} and multiply it by ${\tau}^{1/2}$ to obtain
		\begin{equation}\label{eq: discr energy inequality}
			\begin{split}
			\alpha\tau^{3/2} \normO{\vv_h^j}^2&+\frac{\lamex^2}{2}\tau^{5/2}\normO{\nabla\d_t\mm_h^{j+1}}^2+\frac{\lamex^2}{2}\tau^{1/2}\big(\normO{\nabla\mm_h^{j+1}}^2-\normO{\nabla\mm_h^{j}}^2\big)\\
			&=\tau^{3/2}\scalarproductO{\ff_h^{j+1}}{\vv_h^j}-\frac{\lamex^2}{2}\tau^{5/2}\scalarproductO{\nabla\mm_h^{j+1}}{\nabla\d_t^2\mm_h^{j+1}}.
			\end{split}
		\end{equation}
	The Young inequality yields
	$$
		\tau^{3/2}\scalarproductO{\ff_h^{j+1}}{\vv_h^j}\le\frac{\tau^{3/2}}{4\alpha}\normO{\ff_h^{j+1}}^2+\alpha\tau^{3/2} \normO{\vv_h^j}^2
	$$
	and
	$$
		\scalarproductO{\tau^{1/2}\nabla\mm_h^{j+1}}{\tau^{2}\nabla\d_t^2\mm_h^{j+1}}\le\frac{\tau}{2}\normO{\nabla\mm_h^{j+1}}^2+\frac{\tau^4}{2}\normO{\nabla\d_t^2\mm_h^{j+1}}^2.
	$$
	Putting all together and summing over $j=1,\dots, n-1$, we obtain
	\begin{equation*}
		\begin{split}
		\frac{\lamex^2}{2}\tau^{1/2}\normO{\nabla\mm_h^{n}}^2+\frac{\lamex^2}{2}\tau^{5/2}&\sum_{j=1}^{n-1}\normO{\nabla\d_t\mm_h^{j+1}}^2\le\frac{\tau^{3/2}}{4\alpha}\sum_{j=1}^{n-1}\normO{\ff_h^{j+1}}^2+\frac{\lamex^2}{2}\tau^{1/2}\normO{\nabla\mm_h^{1}}^2\\
		&+\frac{\lamex^2}{4}\tau\sum_{j=1}^{n-1}\normO{\nabla\mm_h^{j+1}}^2+\frac{\lamex^2}{4}\tau^4\sum_{j=1}^{n-1}\normO{\nabla\d_t^2\mm_h^{j+1}}^2.
		\end{split}
	\end{equation*}
Together with~\eqref{enum: punto A1},~\eqref{eq: uniform bound}, and~\eqref{eq: uniform bound S-T norms}, this implies 
	\begin{equation*}
		\begin{split}
		\tau^{5/2}\sum_{j=1}^{n-1}\normO{\nabla\d_t\mm_h^{j+1}}^2
		&\lesssim
		\tau^{1/2} \normOT{\ff_{h\tau}^{+}}^2 + \tau^{1/2}\normO{\nabla\mm_h^{1}}^2
		+ \normOT{\nabla\mm_{h\tau}^{+}}^2 + \tau^4\sum_{j=1}^{n-1}\normO{\nabla\d_t^2\mm_h^{j+1}}^2
		\\& 
		\lesssim 1 + \tau^{1/2} \lesssim 1.
		\end{split}
	\end{equation*}%
	The first convergence in~\eqref{eq: grad v h tau to 0 secondo} now follows immediately as
	\begin{equation*} \tau^2\normOT{\nabla\partial_t\mm_{h\tau}}^2=\tau^3 \sum_{j=0}^{N-1}\normO[\big]{\nabla\d_t\mm_h^{j+1}}^2
	\lesssim\tau^{1/2} \xrightarrow{\tau\to 0} 0.
	\end{equation*}
	Using the norm equivalence~\eqref{eq: norm equivalence grad 2}, this also yields
	\begin{equation*} \tau^2\normOT{\nabla\vv_{h\tau}^-}^2=\tau^3 \sum_{j=0}^{N-1}\normO{\nabla\vv_h^j}^2\stackrel{\eqref{eq: norm equivalence grad 2}}{\simeq}\tau^3\sum_{j=0}^{N-1}\normO[\big]{\nabla\d_t\mm_h^{j+1}}^2\lesssim \tau^{1/2} \xrightarrow{\tau\to 0} 0.\tag*{\qedhere}
	\end{equation*}
\end{proof}
\subsection{Extraction of weakly convergent subsequences}

Uniform boundedness~\eqref{eq: uniform bound S-T norms} allows us to state the following proposition.

\begin{proposition}[Existence of weakly convergent subsequences]\label{prop: weak convergence extraction}
Under the assumptions of Theorem~\ref{thm: main result}, there exists $\mm\in \bH^1(\OT)\cap L^{\infty}(0,T;\bH^1(\Omega))$ with $|\mm|=1$ a.e. in $\OT$ such that, up to the extraction of the same subsequence, all convergences
\begin{subequations}\label{eq: conv}
\begin{align}
 \mm_{h \tau} \rightharpoonup \mm& \qquad in\,\, \bH^1(\OT),
 \label{eq: conv H1}\\
 \mm_{h \tau} \to \mm & \qquad in\,\, \bH^s(\OT) \,\, for\,all \,\, 0<s<1,
 \label{eq: conv Hs}\\
 \mm_{h \tau}, \mm_{h \tau}^{ \pm}, \widehat{\mm}_{h\tau}^+ \to \mm & \qquad in\,\, L^2(0,T;\bH^s(\Omega)) \,\, for\,all \,\, 0<s<1,
 \label{eq: conv L2 0,T,Hs}\\
 \mm_{h \tau}, \mm_{h \tau}^{ \pm}, \widehat{\mm}_{h\tau}^+ \rightharpoonup \mm& \qquad in\,\,{L}^2(0, T ; \bH^1(\Omega))
 \label{eq: conv L2 0,T,H1}\\
 \mm_{h \tau}, \mm_{h \tau}^{ \pm}, \widehat{\mm}_{h\tau}^+ \rightarrow \mm& \qquad in\,\,\mathbf{L}^2(\OT),
 \label{eq: conv L2}\\
 \mm_{h \tau}, \mm_{h \tau}^{ \pm}, \widehat{\mm}_{h\tau}^+\rightarrow \mm& \qquad pointwise\,\,a.e.\,\,in\,\,\OT,
 \label{eq: conv pointwise}\\
 \mm_{h \tau}, \mm_{h \tau}^{\pm}, \widehat{\mm}_{h\tau}^+ \overset{*}{\rightharpoonup} \mm& \qquad in\,\,L^{\infty}(0, T ; \bH^1(\Omega))
 \label{eq: conv L infty},\\
 \vv^-_{h \tau} \rightharpoonup \partial_t \mm& \qquad in\,\,\mathbf{L}^2(\OT),
 \label{eq: conv L2 derivatives}
\end{align}
\end{subequations}
hold simultaneously as $h, \tau \rightarrow 0$.
\end{proposition}

\begin{proof}
The proof is split into five steps. We note that the construction of the subsequence satisfying all the limits~\eqref{eq: conv} follows by successive extraction of subsequences, i.e., Step~2 extracts a further subsequence of the subsequence constructed in Step~1, etc.

{\bf Step~1 (Proof of limits~(\ref{eq: conv H1})--(\ref{eq: conv L infty}) for $\mm_{\boldsymbol{h\tau}}$).}
By uniform boundedness~\eqref{eq: uniform bound S-T norms} and the Eberlein-Šmulian theorem (see, e.g., \cite[page 163]{Conway_functional}), we may extract a weakly convergent subsequence such that $\mm_{h \tau} \rightharpoonup \mm$ in $\bH^1(\OT)$. Due to the compact inclusions $\bH^1(\OT) \Subset \bH^s(\OT) \Subset \bL^2(\OT)$ for all $0 < s < 1$, this yields norm convergence $\mm_{h \tau} \to \mm$ in $\bH^s(\OT)$ and $\bL^2(\OT)$. Moreover, the continuous inclusions $\bH^1(\OT) \subset L^2(0,T;\bH^1(\Omega))$ and $\bH^s(\OT) \subset L^2(0,T;\bH^s(\Omega))$ prove that $\mm_{h \tau} \rightharpoonup \mm$ in $L^2(0,T;\bH^1(\Omega))$ and $\mm_{h \tau} \to \mm$ in $L^2(0,T;\bH^s(\Omega))$. Next, convergence in $\bL^2(\OT)$ allows to extract a further subsequence such that $\mm_{h \tau} \to \mm$ pointwise a.e.\ in $\OT$. Finally, uniform boundedness~\eqref{eq: uniform bound S-T norms} and the Banach--Alaoglu theorem show existence of a further convergent subsequence such that $\mm_{h \tau} \stackrel{*}{\rightharpoonup} \widetilde\mm$ in $L^\infty(0,T;\bH^1(\Omega))$. With the continuous inclusion $L^\infty(0,T; \bH^1(\Omega)) \subset L^2(0,T; \bH^s(\Omega))$ and weak convergence $\mm_{h \tau} \rightharpoonup \mm$ in $L^2(0,T;\bH^1(\Omega))$, uniqueness of the limit allows to identify $\widetilde\mm = \mm$. This concludes the proof of~\eqref{eq: conv H1}--\eqref{eq: conv L infty} for $\mm_{h\tau}$.

{\bf Step~2 (Proof of limits~(\ref{eq: conv L2 0,T,Hs})--(\ref{eq: conv L infty}) for $\mm_{\boldsymbol{h\tau}}^{\boldsymbol\pm}$).}
As in Step~1, we extract a weakly convergent subsequence such that additionally $\mm_{h \tau}^\pm \rightharpoonup \mm^\pm$ in $L^2(0,T;\bH^1(\Omega))$, in $L^2(0,T;\bH^s(\Omega))$ for $0 < s < 1$, and in $\bL^2(\OT)$, and $\mm_{h \tau}^\pm \stackrel{*}{\rightharpoonup} \mm^\pm$ in $L^\infty(0,T;\bH^1(\Omega))$. By continuous inclusion, the limits $\mm^\pm$ in these four spaces coincide (but may differ from $\mm$). To identify $\mm = \mm^\pm$, let us consider $\mm_{h\tau}^-$ (and the analogous argument applies to $\mm_{h\tau}^+$). It holds that
\begin{align*}
 &\norm{\mm_{h \tau}-\mm_{h \tau}^{-}}_{L^2(0,T; \bH^1(\Omega))}^2
 =
 \sum_{j=0}^{N-1} \int_{t_j}^{t_{j+1}} \norm[\Big]{\mm_h^j+\frac{t-t_j}{\tau}(\mm_h^{j+1}-\mm_h^j)-\mm_h^j}_{\bH^1(\Omega)}^2 \d t
 \\& \quad
 = \frac{\tau}{3} \sum_{j=0}^{N-1} \norm{\mm_h^{j+1}-\mm_h^j}_{\bH^1(\Omega)}^2
 \simeq 
 \tau^3 \sum_{j=0}^{N-1} \normO{\d_t\mm_h^{j+1}}^2
 + \tau^3 \sum_{j=0}^{N-1} \normO{\nabla \d_t\mm_h^{j+1}}^2
 \\& \quad
 \le 
 \tau^2 \norm{\partial_t \mm_{h\tau}}_{\OT}^2
 + \tau^2 \norm{\nabla \partial_t \mm_{h\tau}}_{\OT}^2
\xrightarrow{\eqref{eq: uniform bound S-T norms}, \eqref{eq: grad v h tau to 0 secondo}} 0
\quad \text{as } (h, \tau) \to (0,0).
\end{align*}
With $\mm_{h \tau} \rightharpoonup \mm$ and $\mm_{h \tau}^- \rightharpoonup \mm^-$ in $L^2(0,T; \bH^1(\Omega))$, the last estimate thus concludes that $\mm = \mm^-$ by uniqueness of limits. Moreover, norm convergence $\mm_{h\tau} \to \mm$ in $L^2(0,T;\bH^s(\Omega))$ and the continuous inclusion $L^2(0,T;\bH^1(\Omega)) \subset L^2(0,T;\bH^s(\Omega)) \subset \bL^2(\OT)$ proves that $\mm_{h\tau}^- \to \mm$ in $L^2(0,T;\bH^s(\Omega))$ and in $\bL^2(\OT)$. In particular, we may extract a further subsequence such that  $\mm_{h \tau}^- \to \mm$ pointwise a.e.\ in $\OT$. This concludes the proof of~\eqref{eq: conv L2 0,T,Hs}--\eqref{eq: conv L infty} for $\mm_{h\tau}^-$, and the proof for $\mm_{h\tau}^+$ follows analogously.

{\bf Step~3 (Proof of limits~(\ref{eq: conv L2 0,T,Hs})--(\ref{eq: conv L infty}) for $\widehat\mm_{\boldsymbol{h\tau}}^{\boldsymbol+}$).} 
We argue as in Step~2 and further extract a weakly convergent subsequence such that $\widehat\mm_{h \tau}^+ \rightharpoonup \widehat\mm^+$ in $L^2(0,T;\bH^1(\Omega))$, in $L^2(0,T;\bH^s(\Omega))$ for $0 < s < 1$, and in $\bL^2(\OT)$, and $\widehat\mm_{h \tau}^+ \stackrel{*}{\rightharpoonup} \widehat\mm^+$ in $L^\infty(0,T;\bH^1(\Omega))$. As before, the limits $\widehat\mm^+$ coincide. To identify $\widehat\mm^+ = \mm$, we note that
\begin{align*}
 &\norm{\widehat\mm_{h \tau}^+ - \mm_{h \tau}^+}_{L^2(0,T; \bH^1(\Omega))}^2
 = 
 \tau \norm{\mm_h^0 - \mm_h^1}_{\bH^1(\Omega)}^2 
 + \tau \sum_{j=1}^{N-1} \norm{\widehat\mm_h^{j+1} - \mm_h^{j+1}}_{\bH^1(\Omega)}^2
 \\& \quad
 \eqreff{eq: identity5}=
 \tau^3 \norm{\vv_h^0}_{\bH^1(\Omega)}^2
 + \tau^5 \sum_{j=1}^{N-1} \norm{\d_t^2 \mm_h^{j+1}}_{\bH^1(\Omega)}^2
 \stackrel{\eqref{eq: inverse estimate}, \eqref{eq: inverse estimate grad}}\lesssim
 \tau^3 \norm{\vv_h^0}_{\bH^1(\Omega)}^2
 + \tau^3 \sum_{j=0}^{N-1} \norm{\vv_h^j}_{\bH^1(\Omega)}^2
 \\& \quad
 \simeq
 \tau^2 \normOT{\vv_{h\tau}^-}^2
 + \tau^2 \normOT{\nabla \vv_{h\tau}^-}^2
\xrightarrow{\eqref{eq: uniform bound S-T norms}, \eqref{eq: grad v h tau to 0 secondo}} 0
\quad \text{as } (h, \tau) \to (0,0).
\end{align*}
Arguing as before, this concludes the proof of~\eqref{eq: conv L2 0,T,Hs}--\eqref{eq: conv L infty} for $\widehat\mm_{h\tau}^+$.

{\bf Step~4 (Proof of convergence~(\ref{eq: conv L2 derivatives}) for $\vv_{\boldsymbol{h\tau}}^{\boldsymbol{-}}$).}
As in Step~2, we extract a further subsequence such that $\vv_{h\tau}^- \rightharpoonup \vv$ in $\bL^2(\OT)$. 
To show that $\vv = \partial_t \mm$, we use a density argument. From $\mm_{h \tau} \xrightharpoonup{ } \mm$ in $\bH^1(\OT)$, we infer $\partial_t\mm_{h \tau} \xrightharpoonup{ } \partial_t\mm$ in $\mathbf{L}^2(\OT)$ and hence
	$$\scalarproductOT{\vv_{h \tau}^{-} - \partial_t\mm_{h\tau}}{\pphi}\xrightarrow{(h,\tau)\to(0,0)}\scalarproductOT{\vv - \partial_t\mm}{\pphi}\quad\text{for all}\quad \pphi\in \boldsymbol{C}^{\infty}_c(\OT).$$
It thus remains to show that 
	\begin{equation}\label{eq: v = dtm}
		\int_{0}^{T}\scalarproductO{\vv_{h \tau}^{-}(t) - \partial_t\mm_{h\tau}(t)}{\pphi(t)}\d t\xrightarrow{(h,\tau)\to(0,0)}0 \quad \text{for all} \quad \pphi \in \boldsymbol{C}^{\infty}_c(\OT).
	\end{equation}
	Let $\pphi\in \boldsymbol{C}^\infty_c(\OT)$. Equation~\eqref{eq: identity2} and $\vv_{h\tau}^-|_{[t_0,t_1]}=\vv_h^0 = \d _t\mm_h^1=\partial_t\mm_{h\tau}|_{[t_0,t_1]}$ yield
	\begin{align}
		\int_{0}^{T} \scalarproductO{\vv_{h \tau}^{-}(t) - \partial_t\mm_{h\tau}(t)}{\pphi(t)}\d t 
		&= \sum_{j=1}^{N-1}\int_{t_j}^{t_{j+1}}\big\langle\vv_h^{j}-\frac{\mm_h^{j+1}-\mm_h^{j}}{\tau},{\pphi(t)}\big\rangle_\Omega\d t \nonumber\\
		&\mkern-200mu= \frac{1}{2}\sum_{j=1}^{N-1}\int_{t_j}^{t_{j+1}}\big\langle{\frac{\mm_h^{j+1}-\mm_h^j}{\tau}-\frac{\mm_h^j-\mm_h^{j-1}}{\tau}},{\pphi(t)}\big\rangle_\Omega \d t\nonumber\\
		&\mkern-200mu= \frac{1}{2}\sum_{j=1}^{N-1}\int_{t_j}^{t_{j+1}}\scalarproductO{\partial_t\mm_{h\tau}(t)-\partial_t\mm_{h\tau}(t-\tau)}{\pphi(t)}\d t \nonumber\\
		&\mkern-200mu= \frac{1}{2}\int_{t_1}^{T}\scalarproductO{\partial_t\mm_{h\tau}(t)}{\pphi(t)}\d t - \frac{1}{2}\int_{t_1}^{T}\scalarproductO{\partial_t\mm_{h\tau}(t-\tau)}{\pphi(t)}\d t\label{eq: v partial t m line 4}. 
	\end{align}	
	Since $\pphi$ is compactly supported in $\OT$ and $\partial_t\mm_{h\tau}\rightharpoonup\partial_t\mm$ in $\mathbf{L}^2(\OT)$, it holds that
	\begin{equation}\label{eq: comp supp conv 1}\begin{split}
		\int_{t_1}^{T}\scalarproductO{\partial_t{\mm_{h\tau}}(t)}{\pphi(t)}\d t&=\int_{0}^{T}\scalarproductO{\partial_t{\mm_{h\tau}}(t)}{\pphi(t)}\d t
		-\int_{0}^{t_1}\scalarproductO{\partial_t{\mm_{h\tau}}(t)}{\pphi(t)}\d t\\
		&\mkern+120mu\xrightarrow{(h,\tau)\to(0,0)}\int_{0}^{T}\scalarproductO{\partial_t{\mm}(t)}{\pphi(t)}\d t.	
	\end{split}
	\end{equation}
Moreover, $\pphi(\cdot+\tau)\xrightarrow{\tau\to 0}\pphi$ strongly in $\bL^2(\OT)$ and $\partial_t\mm_{h\tau}\rightharpoonup\partial_t\mm$ in $\bL^2(\OT)$ yield
	\begin{equation}\label{eq: comp supp conv 2}
		\begin{split}
		\int_{t_1}^{T}\scalarproductO{\partial_t\mm_{h\tau}(t-\tau)}{\pphi(t)}\d t&=\int_{0}^{T-\tau}\scalarproductO{\partial_t\mm_{h\tau}(s)}{\pphi(s+\tau)}\d s\\
		&\mkern-150mu=\int_{0}^{T}\scalarproductO{\partial_t\mm_{h\tau}(t)}{\pphi(t+\tau)}\d t-\int_{T-\tau}^{T}\scalarproductO{\partial_t\mm_{h\tau}(t)}{\pphi(t+\tau)}\d t\\
		&\xrightarrow{(h,\tau)\to(0,0)}\int_0^T\scalarproductO{\partial_t\mm(t)}{\pphi(t)} \d t.
		\end{split}
	\end{equation}
	Plugging~\eqref{eq: comp supp conv 1} and~\eqref{eq: comp supp conv 2} into~\eqref{eq: v partial t m line 4}, we conclude~\eqref{eq: v = dtm}.
	By density, this implies $\vv=\partial_t\mm$ a.e. in $\OT$. This concludes the proof of~\eqref{eq: conv L2 derivatives}.

\textbf{Step~5 (Proof of modulus constraint $|\mm| \boldsymbol{= 1}$ a.e.\ in $\boldsymbol{\OT}$).}
Thanks to the Hölder inequality and~\eqref{eq: conv L2}, we have that
\begin{align*}
	\norm[\big]{|\mm|^2-|\mm^-_{h\tau}|^2}_{\mathbf{L}^1(\OT)}&\le\normOT{\mm+\mm^-_{h\tau}}\normOT{\mm-\mm^-_{h\tau}}
	\eqreff{eq: uniform bound S-T norms}\lesssim
	\normOT{\mm-\mm^-_{h\tau}}\xrightarrow{(h,\tau)\to(0,0)} 0.
\end{align*}
Moreover, Proposition~\ref{prop: constraint violation} proves
	\begin{align*}
		\norm[\big]{|\mm^-_{h\tau}|^2-1}_{\mathbf{L}^1(\OT)}=&\sum_{j=0}^{N-1}\int_{t_j}^{t_{j+1}}\norm[\big]{|\mm_h^j|^2-1}_{\bL^1(\Omega)}\d t\stackrel{\eqref{eq: constraint violation}}{\lesssim} T\tau\xrightarrow{\tau\to0} 0.
	\end{align*}
Altogether, we thus see
\begin{align*}
	\norm[\big]{|\mm|^2-1}_{\mathbf{L}^1(\OT)}&\le\norm[\big]{|\mm|^2-|\mm^-_{h\tau}|^2}_{\mathbf{L}^1(\OT)}+\norm[\big]{|\mm^-_{h\tau}|^2-1}_{\mathbf{L}^1(\OT)} \xrightarrow{(h,\tau)\to(0,0)} 0
\end{align*}
and hence $|\mm|=1$ a.e.\ in $\OT$.
\end{proof}
%
\subsection{Identification of the limit with a weak solution of LLG}

In this section we will conclude the proof of Theorem~\ref{thm: main result} by showing that the limit function $\mm$ from Proposition~\ref{prop: weak convergence extraction} is a weak solution of LLG~\eqref{eq: LLG system} in the sense of Definition~\ref{def: weak solution}.

\begin{proof}[{\bfseries Proof of Theorem~\ref{thm: main result} (Limit $\mm$ satisfies Definition~{\ref{def: weak solution}\bf(\ref{def: weak solution i})--(\ref{def: weak solution ii})})}]
Proposition~\ref{prop: weak convergence extraction} proves $~\mm \in \bH^1(\OT) \cap L^{\infty}(0, T ; \bH^1(\Omega))$ with $|\mm|=1$ a.e.\ in $\OT$.
Thus, Definition~\ref{def: weak solution}\eqref{def: weak solution i} is verified.
Since $\mm_{h \tau} \rightharpoonup \mm$ in $\bH^1(\OT)$, it follows $\mm_{h \tau}(0) \rightharpoonup \mm(0)$ in $\bH^{1 / 2}(\Omega)$.
By assumption~\eqref{enum: punto A2}, 
	it holds that $\mm_{h \tau}(0)=\mm_h^0 \to \mm^0$ strongly and hence weakly in $\bH^1(\Omega)$. It follows that $\mm(0)=\mm^0$ in the sense of traces, which verifies Definition~\ref{def: weak solution}\eqref{def: weak solution ii}.
\end{proof}

\def\zzeta{\boldsymbol{\zeta}}
\begin{proof}[{\bfseries Proof of Theorem~\ref{thm: main result} (Limit $\mm$ satisfies Definition~{\ref{def: weak solution}\bf(\ref{def: weak solution iii})})}]
The proof is split into six steps, where the core argument is summarized in Step~1. 

\textbf{Step~1 (Discrete variational formulation).}
Let $\vvarphi \in \boldsymbol{C}^{\infty}(\OT)$ be an arbitrary test function.
Define $\zzeta_h\colon[0, T] \rightarrow \SS^1(\TT_h)$ by
$$
 \zzeta_h(t) 
 \coloneqq 
 \III_h[\widehat{\mm}_{h \tau}^{+}(t) \times \vvarphi(t)] \quad \text { for all } t \in[0, T] .
$$
According to the definition of $\widehat{\mm}_{h \tau}^+$, we consider the cases $j=0$ and \mbox{$j=1,\dots,N-1$} separately.
For $j=0$ and all $t \in[t_0, t_1)$ and $\boldsymbol{z} \in \NN_h$, it holds that
$$
\begin{aligned}
 \zzeta_h(t)(\boldsymbol{z}) \cdot {\mm}_h^{0}(\boldsymbol{z}) 
 &=
 \III_h[\widehat{\mm}_{h \tau}^{+}(t) \times \vvarphi(t)](\boldsymbol{z})\cdot \mm_h^0(\boldsymbol{z})
  =
 [\mm_h^0(\boldsymbol{z}) \times \vvarphi(\boldsymbol{z}, t)] \cdot \mm_h^0(\boldsymbol{z})
 =
 0
\end{aligned}
$$
and hence $\zzeta_h(t)\in\bT _h({\mm_h^0})$.
In particular,~\eqref{eq: first step full discr} in Algorithm~\ref{alg: full discr} yields, for all $t \in[t_0, t_1)$, 
\begin{align}\label{eq1:dpr}
 \alpha\scalarproductO{\vv_h^0}{\zzeta_h(t)} + \scalarproductO{\mm_h^0\times\vv_h^0}{\zzeta_h(t)} + \lamex^2\tau\scalarproductO{\nabla\mm_h^1}{\nabla\zzeta_h(t)}
 =
 \scalarproductO{\ff_h^1}{\zzeta_h(t)}.
\end{align}%
For $1 \leq j \leq N-1$ and all $t \in[t_j, t_{j+1})$ and $\boldsymbol{z} \in \NN_h$, it holds that
$$
\begin{aligned}
 \zzeta_h(t)(\boldsymbol{z})\mkern-4mu\cdot\mkern-4mu\widehat{\mm}_h^{j+1}(\boldsymbol{z}) & \mkern-4mu=\III_h[\widehat{\mm}_{h \tau}^{+}(t) \times \vvarphi(t)](\boldsymbol{z}) \cdot\widehat{\mm}_h^{j+1}(\boldsymbol{z})=\mkern-5mu[\widehat{\mm}_h^{j+1}(\boldsymbol{z}) \times \vvarphi(\boldsymbol{z}, t)]\mkern-5mu\cdot\mkern-3mu \widehat{\mm}_h^{j+1}(\boldsymbol{z})=0
\end{aligned}
$$
and hence $\zzeta_h(t)\in\bT_h({\widehat{\mm}_h^{j+1}})$. 
In particular,~\eqref{eq: full discr alg B} in Algorithm~\ref{alg: full discr} yields,  for all $t\in[t_j,t_{j+1})$, 
\begin{align}\label{eq2:dpr}
 \alpha\scalarproductO{\vv_h^j}{\zzeta_h(t)} &+ \scalarproductO{\widehat{\mm}_h^{j+1}\times\vv_h^j}{\zzeta_h(t)} +  \lamex^2\tau\scalarproductO{\nabla\mm_h^{j+1}}{\nabla\zzeta_h(t)}
  =
  \scalarproductO{\ff_h^{j+1}}{\zzeta_h(t)}.
\end{align}
We integrate~\eqref{eq1:dpr} over $[t_0,t_1)$ and~\eqref{eq2:dpr} over $[t_j, t_{j+1})$ and sum these integrals over $j = 0, \dots, N-1$. 
Recalling the definitions~\eqref{eq:mhtau-plus}--\eqref{eq:vhtau-minus} of the interpolations $\mm_{h\tau}^+$, $\widehat{\mm}_{h\tau}^+$ and $\vv_{h\tau}^-$, this leads to
	\begin{equation}\label{eq: weak formulation}
		\begin{aligned}
		 \alpha \int_0^T\langle \vv_{h \tau}^{-}(t), \zzeta_h(t)\rangle _\Omega\d t &+ \int_0^T\langle\widehat{\mm}_{h \tau}^{+}(t) \times \vv^-_{h \tau}(t), \zzeta_h(t)\rangle_\Omega\d t\\
		 &\mkern-50mu+\lambda_{\rm e x}^2 \int_0^T\langle\nabla \mm_{h \tau}^{+}(t), \nabla \zzeta_h(t)\rangle _\Omega \d t =\int_0^T\langle \ff_{h \tau}^{+}(t), \zzeta_h(t)\rangle _\Omega\d t.
		\end{aligned}
	\end{equation}
The subsequent Steps~3--6 consider the convergence of each of these four integrals separately and we will prove that the limit equation reads
\begin{equation}\label{eq3:dpr}
	\begin{aligned}
		&
		-\alpha\int_0^T\langle \mm(t) \times \partial_t \mm(t), \vvarphi(t)\rangle _\Omega\d t
		+ \int_0^T\langle\partial_t \mm(t), \vvarphi(t)\rangle _\Omega\d t
		\\ 
		&\mkern+50mu-\lambda_{\rm e x}^2 \int_0^T\langle\mm(t) \times\nabla\mm(t), \nabla\vvarphi(t)\rangle _\Omega\d t = -\int_0^T\langle \mm(t)\times\ff(t), \vvarphi(t)\rangle _\Omega\d t.
	\end{aligned}
\end{equation}
Overall, we have then proved that~\eqref{eq3:dpr} holds for all $\vvarphi \in \boldsymbol{C}^{\infty}(\OT)$ and hence, by density, for all $\vvarphi \in \bH^1(\OT)$. Consequently, this verifies Definition~\ref{def: weak solution}\eqref{def: weak solution iii}.

\textbf{Step~2 (Auxiliary convergence result).}
For all $v \in H^2(\Omega)$, the nodal interpoland $\III_h v \in \SSS^1(\TT_h)$ with respect to the $\gamma$-quasi uniform mesh $\mathcal{T}_h$, satisfies 
\begin{equation}\label{eq: Lemma 1}
 \normO{v-\III_h v} + h\normO{\nabla(v-\III_h v)}
 \lesssim 
 h^2 \normO{D^2 v};
\end{equation}
see, e.g., \cite[Theorem~4.4.4 and Theorem~4.4.20]{brennerscott}.
Moreover, for all $\pphi_h \in \SSS^1(\TT_h)$ and all $\vvarphi \in \boldsymbol{C}^{\infty}(\overline{\Omega})$, it holds that
	\begin{equation}\label{eq: Lemma 2}
		\sum_{K \in \TT_h}\norm{D^2(\pphi_h \times \vvarphi)}_{\mathbf{L}^2(K)}^2 \lesssim\norm{\pphi_h}_{\bH^1(\Omega)}^2\norm{\vvarphi}_{{\bf W}^{2, \infty}(\Omega)}^2.
\end{equation}
	Therefore, we have that
		$$
		\begin{aligned}
		&\norm[\big]{(\III_h-\mathbf{Id})[\widehat{\mm}_{h\tau}^+(t) \times \vvarphi(t)]}_{\bH^1(\Omega)}^2 =\sum_{K \in \TT_h}\norm[\big]{(\III_h-\mathbf{Id})[\widehat{\mm}_{h\tau}^+(t) \times \vvarphi(t)]}_{\bH^1(K)}^2
		\\& \qquad
		\stackrel{\eqref{eq: Lemma 1}}{\lesssim} h^2 \sum_{K \in \TT_h}\norm{D^2(\widehat{\mm}_{h\tau}^+(t) \times \vvarphi(t))}_{\bL^2(K)}^2 \stackrel{\eqref{eq: Lemma 2}}{\lesssim} h^2\norm{\widehat{\mm}_{h\tau}^+(t)}_{\bH^1(\Omega)}^2\norm{\vvarphi(t)}_{{\bf W}^{2, \infty}(\Omega)}^2 .
		\end{aligned}
		$$
		In particular, it follows that
	\begin{subequations}\label{eq: bound (I-Ih)}
		\begin{equation}
		\norm{(\III_h-\mathbf{Id})[\widehat{\mm}_{h\tau}^+ \times \vvarphi]}_{L^{2}(0, T ; \bH^1(\Omega))} 
		\lesssim 
		h\norm{\widehat{\mm}_{h\tau}^+}_{L^{2}(0, T ; \bH^1(\Omega))}\norm{\vvarphi}_{{\bf W}^{2, \infty}(\OT)} 
		\eqreff{eq: uniform bound S-T norms}\lesssim h \to 0.
	\end{equation}
Note that $\widehat\mm_{h\tau}^+ \to \mm$ in $L^2(0,T; \bH^s(\Omega))$ yields $\widehat\mm_{h\tau}^+ \times \vvarphi \to \mm  \times \vvarphi$ in $L^2(0,T; \bH^s(\Omega))$ for $0 < s < 1$ (with weak convergence for $s = 1$). In particular, we thus see
\begin{align}
 \zzeta_h = \III_h[\widehat{\mm}_{h\tau}^+ \times \vvarphi] 
 \xrightarrow{(h,\tau) \to (0,0)} \mm \times \vvarphi
 \quad \text{in } L^2(0,T; \bH^s(\Omega)).
\end{align}
\end{subequations}%

\textbf{Step~3 (Convergence of first integral in~(\ref{eq: weak formulation})--(\ref{eq3:dpr})).}
Recall that $\langle f_n, g_n \rangle_{L^2} \to \langle f, g \rangle_{L^2}$ provided that $f_n \rightharpoonup f$ and $g_n \to g$ in $L^2$. Applied to our setting with weak convergence~\eqref{eq: conv L2 derivatives} of $\vv_{h\tau}^-$ and strong convergence~\eqref{eq: bound (I-Ih)} of $\zzeta_h$, this shows
\begin{align}\label{eq: conv first term}
\begin{split}
 \int_0^T\langle \vv_{h \tau}^{-}(t), \zzeta_h(t)\rangle _\Omega\d t
 \xrightarrow{(h,\tau)\to(0,0)}
 &\int_0^T\langle \partial_t\mm(t), \mm(t) \times \vvarphi(t) \rangle _\Omega\d t
 \\&
 =
 - \int_0^T\langle \mm(t) \times \partial_t\mm(t), \vvarphi(t) \rangle _\Omega\d t.
\end{split}
\end{align}%

\textbf{Step~4 (Convergence of second integral in~(\ref{eq: weak formulation})--(\ref{eq3:dpr})).}
	The generalized Hölder inequality (with $p_1 = p_2 = 1/4$ and $p_3 = 1/2$) and the Sobolev embedding $\bH^1(\Omega)\subset\bL^4(\Omega)$ prove that
		\begin{equation*}
			\begin{split}
				&\mkern-50mu\bigg|\int_0^T \scalarproductO{\widehat{\mm}_{h \tau}^{+}(t) \times \vv^-_{h \tau}(t)}{(\III_h-\mathbf{Id})[\widehat{\mm}_{h\tau}^+(t)\times\vvarphi(t)]} \d t\bigg|\\
				&\le \int_0^T \norm{\widehat{\mm}_{h \tau}^{+}(t)}_{\bL^4(\Omega)} \norm{\vht(t)}_{\bL^2(\Omega)} \norm[\big]{(\III_h-\mathbf{Id})[\widehat{\mm}_{h\tau}^+(t)\times\vvarphi(t)]}_{\bL^4(\Omega)} \d t\\
				&\lesssim \int_0^T \norm{\widehat{\mm}_{h \tau}^{+}(t)}_{\bH^1(\Omega)} \norm{\vht(t)}_{\bL^2(\Omega)} \norm[\big]{(\III_h-\mathbf{Id})[\widehat{\mm}_{h\tau}^+(t)\times\vvarphi(t)]}_{\bH^1(\Omega)} \d t\\
				&\lesssim \norm{\mhatp}_{L^\infty(0,T;\bH^1(\Omega))} \norm{\vht}_{\bL^2(\OT)} \norm[\big]{(\III_h-\mathbf{Id})[\widehat{\mm}_{h\tau}^+\times\vvarphi]}_{L^2(0,T;\bH^1(\Omega))}\\
				&\stackrel{\eqref{eq: uniform bound S-T norms}}{\lesssim} \norm[\big]{(\III_h-\mathbf{Id})[\widehat{\mm}_{h\tau}^+\times\vvarphi]}_{L^2(0,T;\bH^1(\Omega))}\xrightarrow{\eqref{eq: bound (I-Ih)}} 0\quad\text{as}\,\,(h,\tau)\to(0,0).
			\end{split}
		\end{equation*}
Recall the continuous embedding $\mathbf{W}^{s,p}(\Omega)\subset \bL^q(\Omega)$ for any $q\in[1,\frac{dp}{d-sp}]$ (see, e.g., \cite{HsSpaces}). Hence, if $3/4\le s <1$ (if $d=2$, then $1/2\le s <1$ is enough) we have that $\bH^s(\OT)\subset\bL^4(\OT)$. Thanks to the generalized Hölder inequality and this continuous embedding we have that
\begin{eqnarray*}
\normOT[\big]{|\mhatp|^2-1}^2 
&=& 
\int_{0}^{T} \normO[\big]{|\mhatp(t)|^2-|\mm(t)|^2}^2 \d t
\\
&\le& 
\int_{0}^{T} \norm{\mhatp(t)-\mm(t)}_{\bL^4(\Omega)}^2 \, \norm{\mhatp(t)+\mm(t)}_{\bL^4(\Omega)}^2 \d t
\\
&\lesssim& 
\int_{0}^{T} \norm{\mhatp(t)-\mm(t)}_{\bH^s(\Omega)}^2 \, \norm{\mhatp(t)+\mm(t)}_{\bH^1(\Omega)}^2 \d t
\\
&\le& 
\norm{\mhatp - \mm}^2_{L^2(0,T;\bH^s(\Omega))}\,\norm{\mhatp + \mm}^2_{L^\infty(0,T;\bH^1(\Omega))} 
\\
&\stackrel{\eqref{eq: conv L infty}}{\lesssim}& 
\norm{\mhatp-\mm}^2_{L^2(0,T;\bH^s(\Omega))}\xrightarrow{\eqref{eq: conv L2 0,T,Hs}} 0,
\end{eqnarray*}
i.e., it holds that $|\mhatp|^2\to 1$ in $\bL^2(\OT)$. With the relation
	$$
	(\mathbf{a} \times \mathbf{b}) \cdot(\mathbf{c} \times \mathbf{d})=(\mathbf{a} \cdot \mathbf{c})(\mathbf{b} \cdot \mathbf{d})-(\mathbf{a} \cdot \mathbf{d})(\mathbf{b} \cdot \mathbf{c}) \text { for all } \mathbf{a}, \mathbf{b}, \mathbf{c}, \mathbf{d} \in \mathbb{R}^3 \text {, }
	$$
weak convergence $\vv_{h\tau}^-\cdot\vvarphi\rightharpoonup\partial_t\mm\cdot\vvarphi$ and strong convergence $|\mhatp|^2\to 1$ in $\bL^2(\OT)$, and  orthogonality $\langle \vv_{h\tau}^-(t), \widehat\mm_{h\tau}^+(t) \rangle_\Omega = 0$, it follows that
\begin{align*}
\int_0^T\scalarproductO{\widehat{\mm}_{h \tau}^{+}(t) \times \vv^-_{h \tau}(t)}{\widehat{\mm}_{h \tau}^{+}(t) \times \vvarphi(t)} \d t 
&= 
\int_0^T\scalarproductO{|\mhatp(t)|^2}{\vv_{h\tau}^-(t)\cdot\vvarphi(t)} \d t 
\\& \quad
\xrightarrow{(h, \tau) \to (0,0)}
\int_0^T\langle\partial_t\mm(t),\vvarphi(t)\rangle _\Omega \d t.
\end{align*}
This implies that
		\begin{align}\label{eq: conv second term}
		\int_0^T \scalarproductO{\widehat{\mm}_{h \tau}^{+}(t) \times \vv^-_{h \tau}(t)}{\zzeta_h(t)} \d t&=\int_0^T \scalarproductO{\widehat{\mm}_{h \tau}^{+}(t) \times \vv^-_{h \tau}(t)}{(\III_h-\mathbf{Id})[\widehat{\mm}_{h\tau}^+(t)\times\vvarphi(t)]} \d t\nonumber\\
		&\mkern-210mu+ \int_0^T \scalarproductO{\widehat{\mm}_{h \tau}^{+}(t) \times \vv^-_{h \tau}(t)}{\widehat{\mm}_{h\tau}^+(t)\times\vvarphi(t)} \d t\xrightarrow{(h, \tau) \rightarrow(0,0)}\int_0^T\scalarproductO{\partial_t \mm(t)}{\vvarphi(t)} \d t.
		\end{align}

\textbf{Step~5 (Convergence of third integral in~(\ref{eq: weak formulation})--(\ref{eq3:dpr})).}
As in Step~4, we see
\begin{align*}
 \Big|\int_0^T\langle\nabla \mm_{h \tau}^{+}(t), \nabla(\III_h-\mathbf{Id})&[\widehat{\mm}_{h \tau}^{+}(t) \times \vvarphi(t)]\rangle _\Omega \d t\Big|
 \le\normOT{\nabla\mm_{h\tau}^+}\normOT[\big]{\nabla(\III_h-\mathbf{Id})[\widehat{\mm}_{h\tau}^+\times\vvarphi]}
 \\& \quad
 \stackrel{\eqref{eq: uniform bound S-T norms}}\lesssim\normOT[\big]{\nabla(\III_h-\mathbf{Id})[\widehat{\mm}_{h\tau}^+\times\vvarphi]}\xrightarrow{\eqref{eq: bound (I-Ih)}} 0\quad\text{as}\,\,(h,\tau)\to(0,0).
\end{align*}
Moreover, equation~\eqref{eq: grad v h tau to 0 secondo} from Proposition~\ref{prop: grad v h tau to 0} proves
		\begin{equation*}
			\begin{split}
				\tau\Big|\int_0^T\scalarproductO{\nabla\partial_t\mm_{h\tau}(t)}{\nabla\mhatp(t)\times\vvarphi(t)}\d t\Big|
				&\eqreff{eq: uniform bound S-T norms}\lesssim \tau\normOT{\nabla \partial_t \mm_{h\tau}}\xrightarrow{(h,\tau)\to (0,0)} 0
			\end{split}
		\end{equation*}
		as well as
		\begin{equation*}
			\begin{split}
				&\tau\Big|\int_{t_1}^T\scalarproductO{\nabla\partial_t\mm_{h\tau}(t-\tau)}{\nabla\mhatp(t)\times\vvarphi(t)}\d t\Big|\\
				&=\tau\Big|\int_0^{T-\tau}\scalarproductO{\nabla\partial_t\mm_{h\tau}(t)}{\nabla\mhatp(t)\times\vvarphi(t)}\d t\Big|
				\eqreff{eq: uniform bound S-T norms}\lesssim \tau\normOT{\nabla \partial_t \mm_{h\tau}}\xrightarrow{(h,\tau)\to (0,0)} 0.
			\end{split}
		\end{equation*}
Using the identity~\eqref{eq: identity5} from Lemma~\ref{lem: identities}, it follows that
\begin{equation*}
\begin{split}
&\int_{0}^{T}\scalarproductO{\nabla\mm_{h\tau}^+(t)}{\nabla{\mhatp(t)\times\vvarphi(t)}}\d t
=
\int_{0}^{T}\scalarproductO{\nabla\mm_{h\tau}^+(t)-\nabla\mhatp(t)}{\nabla\mhatp(t)\times\vvarphi(t)}\d t
\\&
\,\,=\,\, 
\int_{t_0}^{t_1}\scalarproductO{\nabla[\mm_h^1-\mm_h^0]}{\nabla\mm_h^0\times\vvarphi(t)}\d t
+ \sum_{j=1}^{N-1} \int_{t_j}^{t_{j+1}}\scalarproductO{\nabla[\mm_h^{j+1}-\widehat{\mm}_h^{j+1}]}{\nabla\widehat{\mm}_h^{j+1}\times\vvarphi(t)}\d t
\\&
\stackrel{\eqref{eq: identity5}}{=}
\int_{t_0}^{t_1}\tau\scalarproductO{\nabla\d_t\mm_h^{1}}{\nabla\mm_h^0\times\vvarphi(t)}\d t
+ \sum_{j=1}^{N-1} \int_{t_j}^{t_{j+1}}\tau\scalarproductO{\nabla\d_t\mm_h^{j+1}}{\nabla\widehat{\mm}_h^{j+1}\times\vvarphi(t)}\d t
\\&\mkern+275mu
- \sum_{j=1}^{N-1} \int_{t_j}^{t_{j+1}}\tau\scalarproductO{\nabla\d_t\mm_h^{j}}{\nabla\widehat{\mm}_h^{j+1}\times\vvarphi(t)}\d t
\\& 
= 
\int_0^T\tau\scalarproductO{\nabla\partial_t\mm_{h\tau}(t)}{\nabla\mhatp(t)\times\vvarphi(t)}\d t-\int_{t_1}^T\tau\scalarproductO{\nabla\partial_t\mm_{h\tau}(t-\tau)}{\nabla\mhatp(t)\times\vvarphi(t)}\d t
\\&
\xrightarrow{(h,\tau)\to(0,0)} 0.
\end{split}
\end{equation*}
From this as well as weak convergence $\nabla\mm_{h\tau}^+\rightharpoonup\nabla\mm$ and strong convergence $\widehat{\mm}_{h \tau}^{+} \times \nabla\vvarphi \to \mm \times \nabla\vvarphi$ in $\bL^2(\OT)$, we derive
\begin{align*}
 &\int_0^T\scalarproductO{\nabla \mm_{h \tau}^{+}(t)}{\nabla [\widehat{\mm}_{h \tau}^{+}(t) \times \vvarphi(t)]} \d t 
 \\& \quad
 = 
 \int_0^T\scalarproductO{\nabla \mm_{h \tau}^{+}(t)}{\widehat{\mm}_{h \tau}^{+}(t) \times \nabla\vvarphi(t)} \d t
 + \int_0^T\scalarproductO{\nabla \mm_{h \tau}^{+}(t)}{\nabla\widehat{\mm}_{h \tau}^{+}(t)\times\vvarphi(t)} \d t
 \nonumber
 \\& \quad
 \xrightarrow{(h, \tau) \rightarrow(0,0)}
 \int_0^T\langle\nabla \mm(t), \mm(t)\times\nabla\vvarphi(t)\rangle _\Omega\d t
 =-\int_0^T\langle\mm(t) \times\nabla\mm(t), \nabla\vvarphi(t)\rangle _\Omega\d t.
		\end{align*}
Altogether, this implies that
		\begin{align}
			&\int_0^T\scalarproductO{\nabla \mm^+_{h\tau}(t)}{\nabla \zzeta_h(t)} \d t = \int_0^T\scalarproductO{\nabla \mm_{h \tau}^{+}(t)}{\nabla\III_h[\widehat{\mm}_{h \tau}^{+}(t) \times \vvarphi(t)]} \d t\nonumber\\
			&= \int_0^T\mkern-7mu\scalarproductO{\nabla \mm_{h \tau}^{+}(t)}{\nabla(\III_h-\mathbf{Id})[\widehat{\mm}_{h \tau}^{+}(t) \times \vvarphi(t)]} \d t+ \int_0^T\mkern-7mu\scalarproductO{\nabla \mm_{h \tau}^{+}(t)}{\nabla[\widehat{\mm}_{h \tau}^{+}(t) \times \vvarphi(t)]} \d t\nonumber\\
			&\xrightarrow{(h, \tau) \rightarrow(0,0)}-\int_0^T\scalarproductO{\mm(t) \times\nabla\mm(t)}{\nabla\vvarphi(t)} \d t.\label{eq: conv third term}
		\end{align}

\textbf{Step~6 (Convergence of fourth integral in~(\ref{eq: weak formulation})--(\ref{eq3:dpr})).}
Thanks to~\eqref{eq: bound (I-Ih)} and $\ff_{h \tau}^{+} \rightharpoonup \boldsymbol{f}$ in $\bL^2(\OT)$ by assumption~\eqref{enum: punto A1}, we have
\begin{align}
\begin{split}\label{eq: conv fourth term}
 \int_0^T\scalarproductO{\ff_{h \tau}^{+}(t)}{\zzeta_h} \d t
 \xrightarrow{(h, \tau) \rightarrow(0,0)} 
 &\int_0^T \scalarproductO{\boldsymbol{f}(t)}{\mm(t)\times \vvarphi(t)} \d t
 \\&
 =
 - \int_0^T\scalarproductO{\mm(t)\times\ff(t)}{\vvarphi(t)} \d t.
\end{split}
\end{align}
As already noted in Step~1, the convergences~\eqref{eq: conv first term}--\eqref{eq: conv fourth term} prove that the discrete variational formulation~\eqref{eq: weak formulation} tends to its continuous counterpart~\eqref{eq3:dpr}. This concludes the proof.
\end{proof}


\begin{proof}[{\bfseries Proof of Theorem~\ref{thm: main result} (Limit $\mm$ satisfies Definition~{\ref{def: weak solution}\bf(\ref{def: weak solution iv})})}]
The proof is split into five steps.

\textbf{Step~1 (Discrete energy inequality).}
For all $0\le j \le N$, we use the notation $\ff^j\coloneqq \ff(t_j)$. In analogy to $\mm_{h\tau}$ and $\mm_{h\tau}^\pm$, we define approximations $\ff_\tau \in \bH^1(\OT)$ and $\ff^\pm_\tau \in \bL^2(\OT)$ of the external field $\ff$ by
\begin{equation*}
 \ff_\tau(t) \coloneqq \frac{t-t_j}{\tau} \ff^{j+1}+\frac{t_{j+1}-t}{\tau}\ff^j, 
 \,\,\,\, 
 \ff^-_\tau(t) \coloneqq \ff^j, 
 \,\, \text{and} \,\,
 \ff^+_\tau(t) \coloneqq \ff^{j+1}
 \,\,\,\, \text{for all } t \in [t_j,t_{j+1})
\end{equation*}
and all $0\le j \le N-1$.
Since $\ff\in C^1([0,T], \bL^2(\Omega))$, it holds that $\ff_\tau, \ff_\tau^{\pm}\to \ff$ in $\bL^2(\OT)$ and $\partial_t\ff_\tau\to\partial_t\ff$ in $\bL^2(\OT)$ as $\tau\to 0$.

For all $n = 1, \dots, N$, Lemma~\ref{lem: discrete energy identity} and the definition of the energy~\eqref{eq: energy} give
\begin{equation*}
	\begin{split}
		\EE(\mm_h^{n},\ff^{n})-\EE(\mm_h^0,\ff^0)&=\frac{\lamex^2}{2}\big[\normO{\nabla\mm_h^{n}}^2-\normO{\nabla\mm_h^{0}}^2\big]-\scalarproductO{\ff^{n}}{\mm_h^{n}}+\scalarproductO{\ff^{0}}{\mm_h^{0}}\\
		&\mkern-180mu\stackrel{\eqref{eq: energy inequality 2}}{\le} -\alpha \tau \sum_{j=0}^{n-1}\normO{\vv_h^j}^2+\sum_{j=0}^{n-1}\big[\tau\scalarproductO{\ff_h^{j+1}}{\vv_h^j}-\scalarproductO{\ff^{j+1}}{\mm_h^{j+1}}+\scalarproductO{\ff^{j}}{\mm_h^j}\big]+\frac{\lamex^2}{4}(\eta_0+\eta_n).
	\end{split}
\end{equation*}
Thanks to~\eqref{eq: identity2}, the second sum of the right-hand side can be written as 
\begin{equation*}
	\begin{split}
		&\sum_{j=0}^{n-1}\big[\tau\scalarproductO{\ff_h^{j+1}}{\vv_h^j}-\scalarproductO{\ff^{j+1}}{\mm_h^{j+1}}+\scalarproductO{\ff^{j}}{\mm_h^j}\big]\\
		&=\sum_{j=0}^{n-1}\big[
			-\scalarproductO{\ff^{j+1}-\ff^j}{\mm_h^j}-\tau\scalarproductO{\ff^{j+1}-\ff_h^{j+1}}{\vv_h^j}-\scalarproductO{\ff^{j+1}}{\mm_h^{j+1}-\mm_h^j-\tau\vv_h^j}\big]\\
		&\stackrel{\eqref{eq: identity2}}{=}-\sum_{j=0}^{n-1}\big[
			\tau\scalarproductO{\d_t\ff^{j+1}}{\mm_h^j}+\tau\scalarproductO{\ff^{j+1}-\ff_h^{j+1}}{\vv_h^j}\big]+\frac{\tau}{2}\sum_{j=1}^{n-1}\scalarproductO{\ff^{j+1}}{\d_t\mm_h^{j+1}-\d_t\mm_h^j}.\\
	\end{split}
\end{equation*}
Therefore, we obtain that 
\begin{equation*}
	\begin{split}
		\EE(\mm_h^{n},\ff^{n})-\EE(\mm_h^0,\ff^0)+\alpha \tau \sum_{j=0}^{n-1}\normO{\vv_h^j}^2&+\sum_{j=0}^{n-1}\big[
			\tau\scalarproductO{\d_t\ff^{j+1}}{\mm_h^j}+\tau\scalarproductO{\ff^{j+1}-\ff_h^{j+1}}{\vv_h^j}\big]\\
		&\mkern-50mu\le\frac{\tau}{2}\sum_{j=1}^{n-1}\scalarproductO{\ff^{j+1}}{\d_t\mm_h^{j+1}-\d_t\mm_h^j}+\frac{\lamex^2}{4}(\eta_0+\eta_n).
	\end{split}
\end{equation*}
For any $t'\in[0,T)$ and $1\le n\le N$ such that $t'\in[t_{n-1},t_n)$, this can be rewritten as
\begin{align}\label{eq: energy in progress}
 &\EE(\mm_{h\tau}^+(t'),\ff_\tau^+(t')) 
 - \EE(\mm_{h\tau}^-(0),\ff(0))
 + \alpha\int_{0}^{t'} \normO{\vv_{h\tau}^-(t)}^2\d t
 + \int_0^{t'}\scalarproductO{\partial_t\ff_\tau(t)}{\mm_{h\tau}^-(t)}\d t
 \nonumber
 \\ &\qquad
 + \int_{t'}^{t_n}\scalarproductO{\partial_t\ff_\tau(t)}{\mm_{h\tau}^-(t)}\d t
 + \int_{0}^{t_n}\scalarproductO{\ff_\tau^+(t)-\ff_{h\tau}^+(t)}{\vv_{h\tau}^-(t)}\d t
 \\& \quad
 \le
 \frac{1}{2}\int_{t_1}^{t_n}\scalarproductO{\ff_\tau^+(t)}{\partial_t\mm_{h\tau}(t)-\partial_t\mm_{h\tau}(t-\tau)}\d t
 + \frac{\lamex^2}{4}(\eta_0+\eta_n)\nonumber.
\end{align}

\textbf{Step~2 (Convergence of left-hand side of~(\ref{eq: energy in progress})).}
Since $\mm_{h\tau}^+\rightharpoonup\mm$ in $L^2(0,T;\bH^1(\Omega))$ by~\eqref{eq: conv L2 0,T,H1} and $\ff_\tau^+\to\ff$ in $\bL^2(\OT)$ by assumption~\eqref{enum: punto A1}, weakly lower semicontinuity proves
\begin{align*}
 \int_I\EE(\mm(t'),\ff(t'))\d t'
 \le
 \liminf_{(h,\tau)\to (0,0)}\int_I\EE(\mm_{h\tau}^+(t'),\ff_\tau^+(t'))\d t'
 \,\, \text{for any measurable } I \subseteq [0,T].
\end{align*}
The strong convergence $\mm_h^0\to \mm^0$ in $\bH^1(\Omega)$ from assumption~\eqref{enum: punto A2} implies
$$
\EE(\mm^0,\ff(0))
= 
\lim_{h\to 0}\EE(\mm_h^0,\ff(0))
=
\liminf_{(h,\tau)\to(0,0)}\EE(\mm_{h\tau}^-(0),\ff(0)).
$$
By convergence $\vv_{h\tau}^-\rightharpoonup\partial_t\mm$ in $\bL^2(\OT)$ from~\eqref{eq: conv L2 derivatives}, weakly lower semicontinuity yields 
$$
	\int_{0}^{t'}\normO{\partial_t\mm(t)}^2\d t\le\liminf_{(h,\tau)\to(0,0)}\int_{0}^{t'}\normO{\vv_{h\tau}^-(t)}^2\d t.
$$
The convergence $\partial_t\ff_\tau\to\partial_t\ff$ in $\bL^2(\OT)$  by assumption~\eqref{enum: punto A1} and $\mm_{h\tau}^-\to\mm$ in $\bL^2(\OT)$ \mbox{by~\eqref{eq: conv L2}} imply that
$$
\int_{0}^{t'} \scalarproductO{\partial_t\ff(t)}{\mm(t)}\d t
= \lim_{(h,\tau) \to (0,0)}
\int_{0}^{t'}\scalarproductO{\partial_t\ff_\tau(t)}{\mm_{h\tau}^-(t)}\d t.
$$
By no concentration of Lebesgue functions and $t'\le t_n\le t'+\tau$, this also shows that 
$$
0 = 
\lim_{(h,\tau) \to (0,0)}
\int_{t'}^{t_n}\scalarproductO{\partial_t\ff_\tau(t)}{\mm_{h\tau}^-(t)}\d t.
$$
Strong convergence $\ff_\tau^+\to \ff$ and $\ff_{h\tau}^+\to\ff$ in $\bL^2(\OT)$ by Assumption~\eqref{enum: punto A1} prove that
$$
0 = 
\lim_{(h,\tau) \to (0,0)}
\int_0^{t_n}\scalarproductO{\ff_\tau^+(t)-\ff_{h\tau}^+(t)}{\vv_{h\tau}^-(t)}\d t.
$$
This concludes the consideration of the left-hand side of~\eqref{eq: energy in progress}.

\textbf{Step~3 (Convergence of right-hand side of~(\ref{eq: energy in progress})).}
Note that $\ff_\tau^+\to\ff$ in $\bL^2(\OT)$ and $\ff\in C^1(0,T,\bL^2(\Omega))$.
By this and no-concentration of Lebesgue functions, the integral on the right-hand side of~\eqref{eq: energy in progress} satisfies
\begin{align*}
&\int_{t_1}^{t_n}\scalarproductO{\ff_\tau^+(t)}{\partial_t\mm_{h\tau}(t)-\partial_t\mm_{h\tau}(t-\tau)}\d t
\\&
=
\int_{t_1}^{t_n}\scalarproductO{\ff_\tau^+(t)}{\partial_t\mm_{h\tau}(t)}\d t
- \int_{t_1}^{t_n}\scalarproductO{\ff_\tau^+(t)}{\partial_t\mm_{h\tau}(t-\tau)}\d t
\\&
=
\int_{t_1}^{t_n}\scalarproductO{\ff_\tau^+(t)}{\partial_t\mm_{h\tau}(t)}\d t
- \int_{t_0}^{t_{n-1}}\scalarproductO{\ff_\tau^+(t+\tau)}{\partial_t\mm_{h\tau}(t)}\d t
\\&
=
\int_{0}^{t_n}\scalarproductO{\ff_\tau^+(t) - \ff_\tau^+(t+\tau)}{\partial_t\mm_{h\tau}(t)}\d t
\\& \quad
- \int_{0}^\tau \scalarproductO{\ff_\tau^+(t)}{\partial_t\mm_{h\tau}(t)}\d t
+ \int_{t_n-\tau}^{t_n}\scalarproductO{\ff_\tau^+(t+\tau)}{\partial_t\mm_{h\tau}(t)}\d t
\xrightarrow{(h,\tau)\to(0,0)} 0.
\end{align*}
Convergence of the last term $\eta_0+\eta_n\to 0$ follows by assumption \eqref{eq: condition eta0n}. This concludes the consideration of the right-hand side of~\eqref{eq: energy in progress}.

\textbf{Step~4 (Continuous energy inequality).}
For every measurable set $I\subseteq[0,T]$, the convergence statements from Step~2--3 prove
\begin{equation*}
	\begin{aligned}
		& \int_I \Big[\EE(\mm(t'), \boldsymbol{f}(t'))+\alpha\int_0^{t'}\normO{\partial_t \mm(t)}^2 \, \mathrm{d} t \Big] \, \mathrm{d} t' \\
		& \leq \liminf _{(h,\tau)\to(0,0)} \int_I \Big[\EE\big(\mm_{h \tau}^{+}(t'), \boldsymbol{f}_\tau^{+}(t')\big)+\alpha\int_0^{t'}\normO{\vv_{h \tau}^{-}(t)}^2 \, \mathrm{d} t \Big] \, \mathrm{d} t' \\
		& \stackrel{\eqref{eq: energy in progress}}{\leq} \liminf _{(h,\tau)\to(0,0)} \int_I \Big[\EE\big(\mm_{h \tau}^{-}(0), \boldsymbol{f}(0)\big) - \int_0^{t'} \scalarproductO{\partial_t \boldsymbol{f}_\tau(t)}{\mm_{h \tau}^{-}(t)} \, \mathrm{d} t - \int_{t'}^{t_n} \scalarproductO{\partial_t \boldsymbol{f}_\tau(t)}{\mm_{h \tau}^{-}(t)} \, \mathrm{d} t \\
		& \qquad - \int_0^{t_n} \scalarproductO{\boldsymbol{f}_\tau^+(t) - \boldsymbol{f}_{h \tau}^+(t)}{\vv_{h \tau}^-(t)} \, \mathrm{d} t + \frac{1}{2} \int_{t_1}^{t_n} \scalarproductO{\boldsymbol{f}_\tau^+(t)}{\partial_t \mm_{h \tau}(t) - \partial_t \mm_{h \tau}(t-\tau)} \, \mathrm{d} t \\
		& \qquad + \frac{\lamex^2}{4} (\eta_0 + \eta_n) \Big] \, \mathrm{d} t' \\
		& = \int_I \Big[\EE\big(\mm^0, \boldsymbol{f}(0)\big) - \int_0^{t'} \scalarproductO{\partial_t \boldsymbol{f}(t)}{\mm(t)} \, \mathrm{d} t \Big] \, \mathrm{d} t'.
	\end{aligned}
\end{equation*}
From this, we deduce that the integrand satisfies the inequality a.e.\ in $(0, T)$, i.e., 
$$
\EE(\mm(t'), \boldsymbol{f}(t'))
+ \alpha \int_0^{t'} \!\!\!\!  \normO{\partial_t \mm(t)}^2 \! \d t
+ \int_0^{t'} \!\!\!\! \langle\partial_t \boldsymbol{f}(t), \mm(t)\rangle_\Omega \! \d t 
\leq
\EE(\mm^0, \boldsymbol{f}(0))
\, \text{for a.e.\ } t' \in [0,T].
$$
This finally is the energy inequality of Definition~\ref{def: weak solution}\eqref{def: weak solution iv}.

{\bf Step~5 (CFL condition yields convergence (\ref{eq: condition eta0n})).} To prove $\eta_0+\eta_n\to 0$, first note that
\begin{equation}\label{eq63}
 \tau\normO{\d_t\mm_h^j}^2
 \le
 \normOT{\partial_t\mm_{h\tau}}^2
 \le
 \norm{\mm_{h\tau}}^2_{\bH^1(\Omega)}
 \eqreff{eq: uniform bound S-T norms}\lesssim
 1
 \quad \text{for all } j = 1,\dots,N.
\end{equation}
Using boundedness~\eqref{eq: uniform bound} and an inverse estimate, the CFL condition $\tau=o(h^2)$ proves
\begin{equation*}
	\begin{split}
\eta_0+\eta_n&=\frac{1}{4}\big[\normO{\nabla\mm_h^1}^2-\normO{\nabla\mm_h^0}^2+\normO{\nabla\mm_h^{n-1}}^2-\normO{\nabla\mm_h^n}^2\big]\\
&=\frac{\tau}{4}\big[ \scalarproductO{\nabla\mm_h^1+\nabla\mm_h^0}{\nabla\d_t\mm_h^1} -\scalarproductO{\nabla\mm_h^{n-1}+\nabla\mm_h^n}{\nabla\d_t\mm_h^{n}}\big]\\
&\le \frac{\tau}{4}\big[\normO{\nabla\mm_h^1+\nabla\mm_h^0}\normO{\nabla\d_t\mm_h^1} + \normO{\nabla\mm_h^{n-1}+\nabla\mm_h^n}\normO{\nabla\d_t\mm_h^n} \big]\\
&\eqreff{eq: uniform bound}\lesssim 
\tau h^{-1}\big[\normO{\d_t\mm_h^1}+\normO{\d_t\mm_h^n}\big]\eqreff{eq63}\lesssim \sqrt{\tau}h^{-1}\xrightarrow{(h,\tau)\to(0,0)}0.
	\end{split}
\end{equation*}
This concludes the proof.
\end{proof}

\begin{proof}[{\bfseries Proof of Corollary~\ref{cor: main result} (convergence of full sequence)}]
Recall that the proof of Theorem~\ref{thm: main result} exploits the well-known three steps of the energy method in the context of PDEs: 
show that the approximations are uniformly bounded; 
extract weakly convergent subsequences; 
prove that the weak limits are weak solutions of the PDE.
If the weak solution $\mm$ to LLG is unique, then any limit constructed above coincides with $\mm$. In particular, any subsequence of $\mm_{h\tau}$ admits a further subsequence that converges to $\mm$. According to basic calculus, this is equivalent to the fact that the full sequence converges to $\mm$.
\end{proof}

\section{Numerical experiments}\label{section:numerics}

To support our theoretical findings, this section presents some numerical experiments.
Throughout, we restrict to $\Omega\subset \R^2$ and report on the empirical results obtained by our Matlab implementation of Algorithm \ref{alg: full discr}. In Section \ref{subsec: first}--\ref{subsec: second}, we will show numerical evidence of the formal second-order convergence in time and first-order convergence in space of the proposed BDF2-type integrator (Algorithm \ref{alg: full discr}, subsequently abbreviated BDF2). In Section \ref{subsec: third}, we will investigate the behavior for a non-smooth singular solution. In Section \ref{subsec: fourth}, we will investigate the numerical deviation from the unit-length constraint, while the final Section \ref{subsec: CFL} discusses the role of the CFL condition $\tau = o(h^2)$. To this end, we also compare BDF2 with the projection-free first-order tangent-plane scheme (subsequently abbreviated TPS) from \cite{Spin-polarized} and the second-order in time implicit midpoint scheme (subsequently abbreviated MID) from \cite{bp2006}.

\subsection{Experimental convergence rates for time-independent applied field}
\label{subsec: first}

\begin{figure}[t]
\centering
\begin{subfigure}[b]{0.49\linewidth} 
   \centering
   \includegraphics[width=1\linewidth]{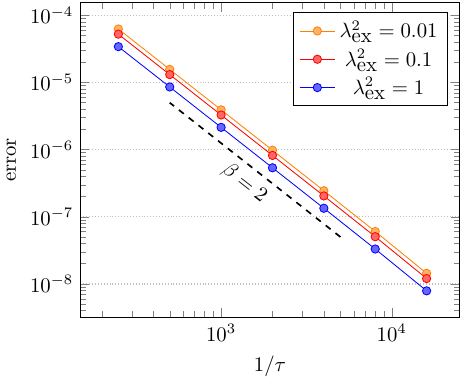}
   \caption{Second-order convergence in time.}
\end{subfigure} 
\begin{subfigure}[b]{0.49\linewidth}
  \centering
  \includegraphics[width=1\linewidth]{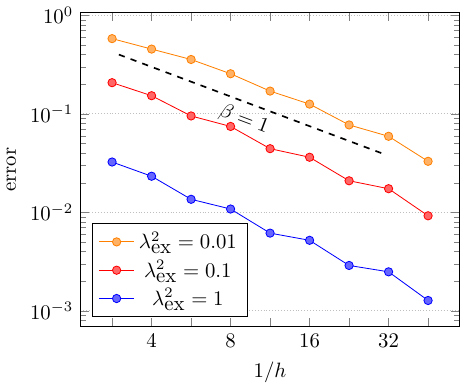}
  \caption{First-order convergence in space.} 
\end{subfigure}
\caption{Empirical convergence rates in of the $\bH^1$-error $\norm{\mm_{\text{ref}}(T)-\mm_{h\tau}(T)}_{\bH^1(\Omega)}$ at final time in the experiment from Section~\ref{subsec: first}.}
\label{fig:example_1_conv_rates}
\end{figure}

We solve~\eqref{eq: LLG system} for $\Omega = (0,1)^2$ and $T = 1$. We choose a constant initial condition $\mm^0 \equiv (0,1,0)$ and a radial external field $\ff(x,y) = (x,y,0)/|(x,y,0)|$ that is constant in time. We fix the damping parameter $\alpha = 0.25$ and choose the exchange constant $\lamex^2 \in \{0.01; 0.1; 1\}$ so that either the external field or the exchange contribution of the effective field dominate. 

To test second-order convergence in time, we fix a uniform triangulation of $\Omega$ consisting of 64 triangles, corresponding to a mesh size of $h = 0.25$. For the time discretization, we employ a uniform partition of the time interval with time-step size $\tau_0 = 4 \cdot 10^{-3}$ that is successively bisected, i.e., $\tau_k = 2^{-k} \tau_0$ for $k = 0, \ldots, 6$. Since the exact solution $\mm$ is unknown, we use a reference solution $\mm_\text{ref}$ obtained by $\tau_{\text{ref}} = 2^{-8} \tau_0$.

To test first-order convergence in space, we fix a time-step size $\tau = 4\cdot 10^{-3}$ and consider a sequence of uniform triangulations of $\Omega$. We start from a coarse initial mesh $\mathcal{T}_h$ consisting of 32 triangles, corresponding to a mesh-size of $h=\sqrt{2}/4\approx 0.3536$, that is successively refined by uniform longest-edge bisection.
As reference solution, we employ the discrete solution obtained from a mesh that is additionally refined twice, consisting of $32\,768$ uniform triangles with mesh-size $h=\sqrt{2}/128\approx 0.01105$.

Figure \ref{fig:example_1_conv_rates} shows the convergence rates in $\tau$ and $h$ of the computed $\bH^1$-error~$\norm{\mm_{\text{ref}}(T)-\mm_{h\tau}(T)}_{\bH^1(\Omega)}$ at final time $T = 1$. As expected, we observe first-order convergence in space and second-order convergence in time.

\subsection{Experimental convergence rates for time-dependent applied field}
\label{subsec: second}

\begin{figure}[ht!]
\centering
\begin{subfigure}[b]{0.49\linewidth} 
  \centering
  \includegraphics[width=1\linewidth]{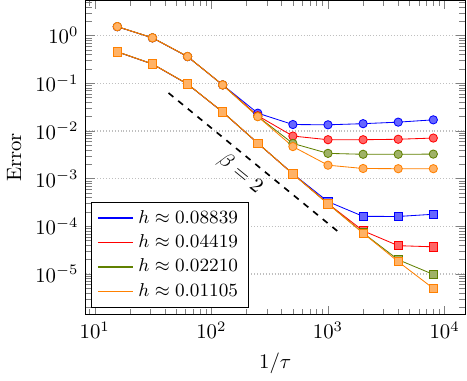}
  \caption{Second-order convergence in time.}
\end{subfigure} 
\begin{subfigure}[b]{0.49\linewidth}  
  \centering
  \includegraphics[width=1\linewidth]{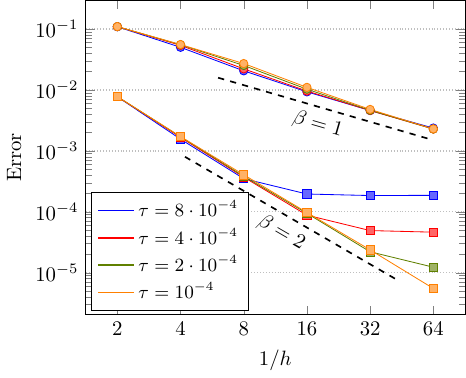}
  \caption{Convergence orders in space.}
\end{subfigure}
\caption{Empirical convergence rates of $\max\limits_{j=1,\dots,N} \norm{\mm(t_j) - \mm_h^j}_{\bL^2(\Omega)}$
 (\protect\tikz[baseline=-0.5ex]{
  \protect\draw[line width=1pt] (-0.4,0) -- (-0.15,0); 
  \protect\draw[line width=1pt] (-0.15,-0.15) rectangle (0.15,0.15); 
  \protect\draw[line width=1pt] (0.15,0) -- (0.4,0); 
 })
 and $\max\limits_{j=1,\dots,N} \norm{\mm(t_j) - \mm_h^j}_{\bH^1(\Omega)}$   
 (\protect\tikz[baseline=-0.5ex]{
  \protect\draw[line width=1pt] (-0.4,0) -- (-0.15,0); 
  \protect\draw[line width=1pt] (0,0) circle [radius=0.15]; 
  \protect\draw[line width=1pt] (0.15,0) -- (0.4,0); 
 })
 in the experiment of Section~\ref{subsec: second}.}
\label{fig:example_2_conv_rates}
\end{figure}

We consider an experiment proposed in~\cite{feischl}, where the exact solution is prescribed by
$$
\mm(\xx,t) \coloneq \begin{pmatrix}
    -(x_1^3-3x_1^2/2+1/4)\sin(3\pi t/T) \\
    \sqrt{1-(x_1^3-3x_1^2/2+1/4)^2} \\
    -(x_1^3-3x_1^2/2+1/4)\cos(3\pi t/T)
\end{pmatrix}.
$$ 
We solve \eqref{eq: LLG system} for $\Omega = (0,1)^2$ and $T = 0.2$. We fix the damping parameter $\alpha=0.2$ and the exchange constant $\lamex^2 = 1$. We compute the external field $\ff$ analytically by inserting the prescribed solution into equation \eqref{eq: alternative LLG}, i.e., 
$$
\ff(\xx,t) = \alpha\partial_t \mm(\xx,t) +\mm(\xx,t)\times \partial_t\mm(\xx,t)- \lamex^2\Delta\mm(\xx,t).
$$
Note that one may ignore the term $(\mm\cdot\heff)\mm$, since the numerical scheme solves the equation in the tangent space. 

To test second-order convergence in time, we fix a spatial mesh $\TT_h$ and consider step-sizes $\tau_k=2^{-k}\tau_0$ with $\tau_0= 0.064$. 
In Figure \ref{fig:example_2_conv_rates} (left), we plot the convergence of the errors 
$$
\max_{j=1,\dots,N}\norm{{\mm^j_{h}}-\mm(t_j)}_{\bH^1(\Omega)}\quad\text{and}\quad\max_{j=1,\dots,N}\norm{{\mm^j_{h}}-\mm(t_j)}_{\bL^2(\Omega)}
$$
over $1/\tau$.
Initially, the error converges with order $2$ in time. However, the error stagnates after a certain time-step, since the spatial error starts to dominate over the temporal error. To investigate this behavior, we consider four different fixed meshes with $512$, $2048$, $8192$, and $32768$ triangles, respectively.

To test first-order convergence in space, we fix the time-step size and consider spatial meshes $\TT_h$ consisting of uniform triangles with mesh sizes $h_k = 2^{-k}h_0$ and $h_0=1/2$, respectively.
In Figure \ref{fig:example_2_conv_rates} (right), we plot the convergence of the $\bL^2$- and $\bH^1$-errors over $1/h$.
Initially, the error converges with order $1$ in $h$ for the $\bH^1$-norm and order $2$ in $h$ for the $\bL^2$-norm. Again, however, the error stagnates after a certain refinement, since the temporal error starts to dominate over the spatial error.

\subsection{Experimental convergence behavior for a problem with finite-time blow-up}
\label{subsec: third}

As a third experiment, we consider the example proposed in \cite{bp2006,bkp2008}, which is explicitly constructed in order to induce a singularity appearing in the center of the domain. The setting proposed leads to numerical approximations with a large gradient, i.e., numerical occurrence of a finite-time blow-up. We refer to \cite{bkp2008} for more details on the construction of the problem and for a visual representation of the evolution in time of this solution.
We note, however, that \cite{bdfk2023adaptive} proposes an adaptive strategy in time and space, which prevents the singularity from forming; see \cite[Example 3]{bdfk2023adaptive}.

Let $\Omega=(-1/2,1/2)^2$. With $A = (1-2\abs{x})^4$, define the initial data $\mm^0\colon\Omega\to\R^3$ by
\begin{equation}\label{eq: initial data blow up}
\mm^0(\xx) := 
\begin{cases}
(0, 0, -1) & \text{if } |\xx| \ge \dfrac{1}{2}, \\
\dfrac{\Big(2x_1A, 2x_2A, A^2-|\xx|^2\Big)}{A+|\xx|^2} & \text{if } |\xx| \le \dfrac{1}{2}.
\end{cases}
\end{equation}
We choose $\lamex^2=1$, $\alpha = 0.25$ and let the dynamics evolve without the presence of an external field, i.e., $\ff\equiv 0$. We choose the final time $T=0.3$.

\begin{figure}[ht!]
    \centering
    \begin{subfigure}[b]{0.52\linewidth}
        \centering
        \includegraphics[width=1\linewidth]{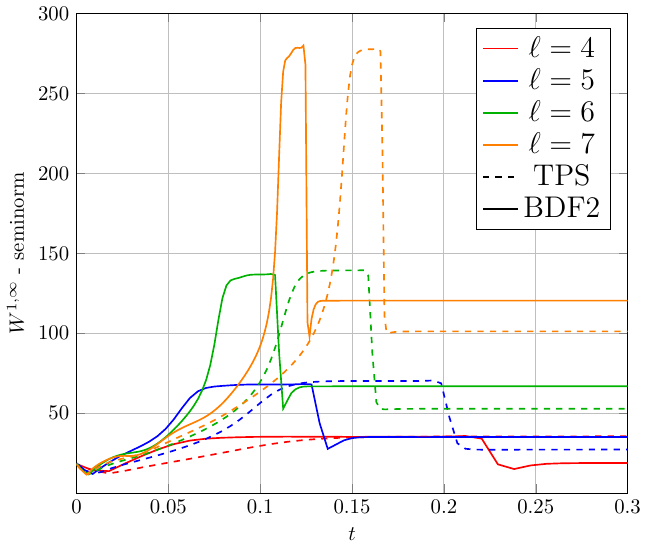}
    \end{subfigure}
    \begin{subfigure}[b]{0.52\linewidth}  
        \centering
        \includegraphics[width=1\linewidth]{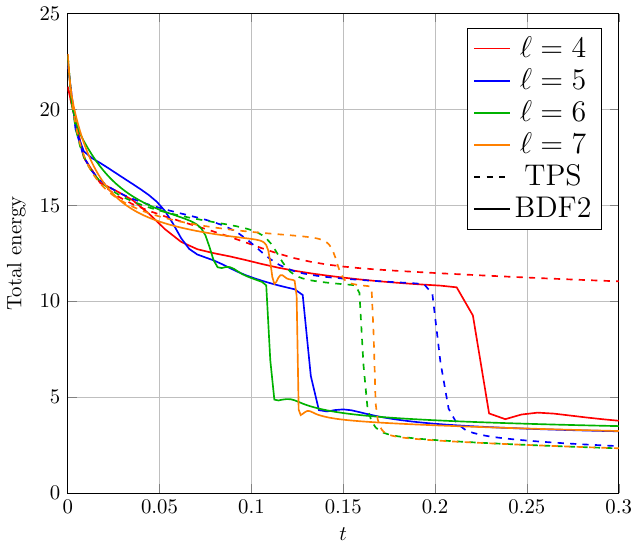}
    \end{subfigure}
    \caption{Evolution of the $\mathbf{W}^{1,\infty}$-seminorm in time (left) and energy decay (right) for the singular behavior experiment, depending on the mesh size $h$.}
    \label{fig:example_blowup_h} 
\end{figure}

\begin{figure}[ht!]
    \centering
    \begin{subfigure}[b]{0.52\linewidth}
        \centering
        \includegraphics[width=1\linewidth]{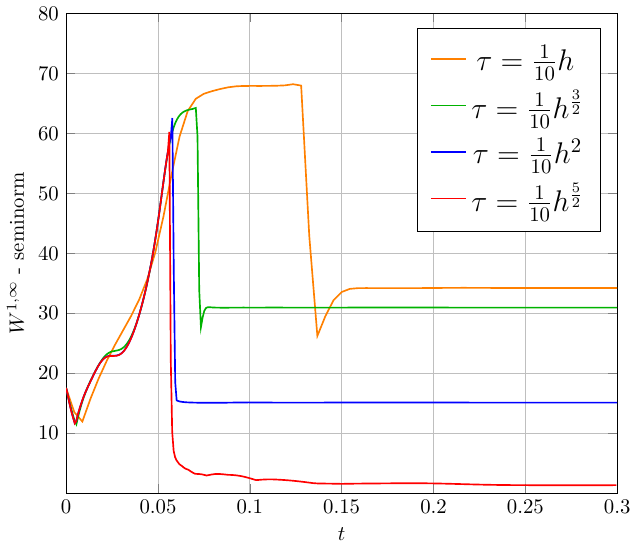}
    \end{subfigure}
    \begin{subfigure}[b]{0.52\linewidth}  
        \centering
        \includegraphics[width=1\linewidth]{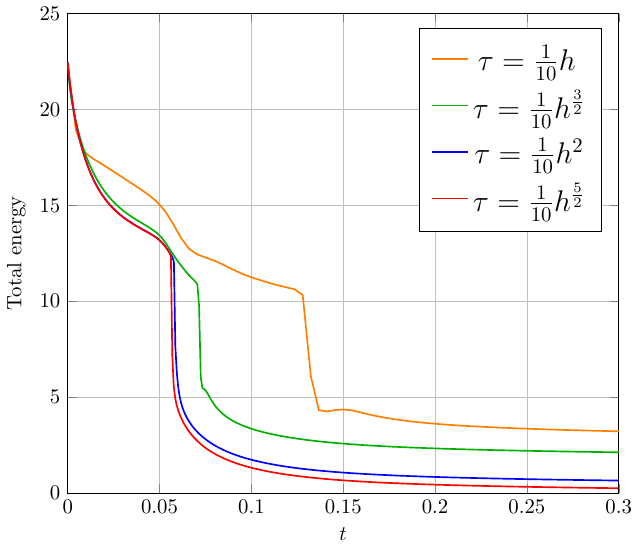}
    \end{subfigure}
    \centering
    \begin{subfigure}[b]{0.52\linewidth}
        \centering
        \includegraphics[width=1\linewidth]{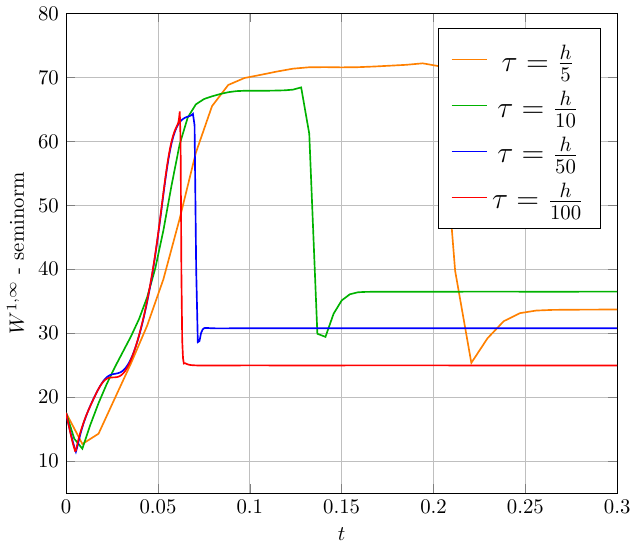}
    \end{subfigure}
    \begin{subfigure}[b]{0.52\linewidth}  
        \centering
        \includegraphics[width=1\linewidth]{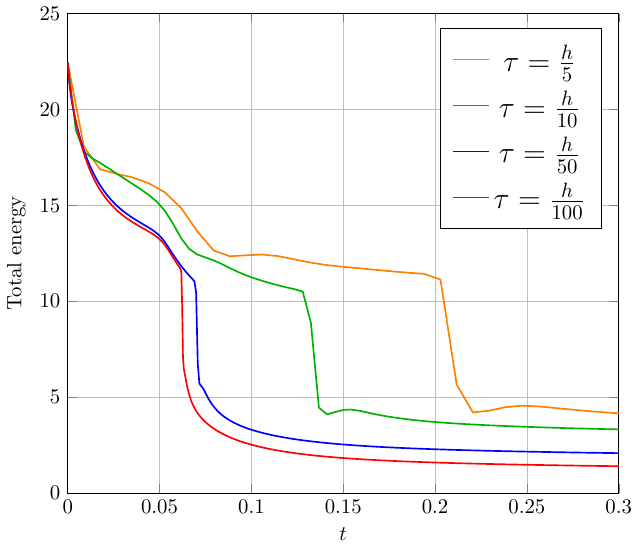}
    \end{subfigure}
    \caption{Evolution of the $\mathbf{W}^{1,\infty}$-seminorm in time (left) and energy decay (right) for the singular behavior experiment, depending on the time-step size $\tau$. Top line: $\tau~\in~\big\{\frac{h}{10}, \frac{h^{3/2}}{10}, \frac{h^{2}}{10}, \frac{h^{5/2}}{10}\big\}$ ; Bottom line: $\tau \in\big\{\frac{h}{5}, \frac{h}{10}, \frac{h}{50}, \frac{h}{100}\big\}$.}
    \label{fig:example_blowup_taus} 
\end{figure}

In Figure \ref{fig:example_blowup_h}, we display the evolution of the $\mathbf{W}^{1,\infty}$-seminorm $\norm{\nabla \mm}_{\bL^\infty(\Omega)}$ in time and the decay of the energy \eqref{eq: energy}, comparing the results coming from BDF2 to those coming from TPS. 
Similarly to the setting of \cite{bp2006}, we consider triangulations consisting of $2^{2\ell+1}$ halved squares with edge length $\lambda\coloneq 2^{-\ell}$, for $\ell\in~\{4,5,6,7\}$, mesh size $h_\ell=~2^{-\ell+1/2}$ and set the time-step size as $\tau=h_\ell/10$. We observe that for decreasing mesh sizes the blow-up time (i.e., the time when $\norm{\nabla\mm(t)}_{\bL^\infty(\Omega)}$ reaches its maximum) increases, as also observed in \cite{bp2006,bkp2008,bppr2014}.
Moreover, we notice that for the same choice of the discretization parameters, BDF2 yields an earlier blow-up time than TPS. To investigate how the blow-up time and the stagnation of the $\mathbf{W}^{1,\infty}$-seminorm depend on the discretization parameters, we plot the same quantities for a fixed triangulation corresponding to $\ell = 5$ (i.e., $\#\mathcal{T}_h = 2^{11} = 2048$ elements; $h_5 = 2^{-5}\sqrt{2}\approx 0.0442$) and different time-step sizes $\tau~\in\big\{\frac{h}{10}, \frac{h^{3/2}}{10}, \frac{h^{2}}{10}, \frac{h^{5/2}}{10}\big\}$ (Figure \ref{fig:example_blowup_taus}, top line) and $\tau~\in\big\{\frac{h}{5}, \frac{h}{10}, \frac{h}{50}, \frac{h}{100}\big\}$ (Figure \ref{fig:example_blowup_taus}, bottom line). In order to avoid an overloading of plots, only the results obtained with BDF2 are reported. We notice that for finer stepsize a limit behavior in the energy decay seems to be reached, and the stagnation of the $\mathbf{W}^{1,\infty}$-seminorm does not seem to appear anymore, independently of the choice of $\tau$ respecting or not the CFL condition $\tau = o(h^2)$.

\begin{figure}[ht!]
    \centering
    \begin{subfigure}[b]{0.52\linewidth}
        \centering
        \includegraphics[width=1\linewidth]{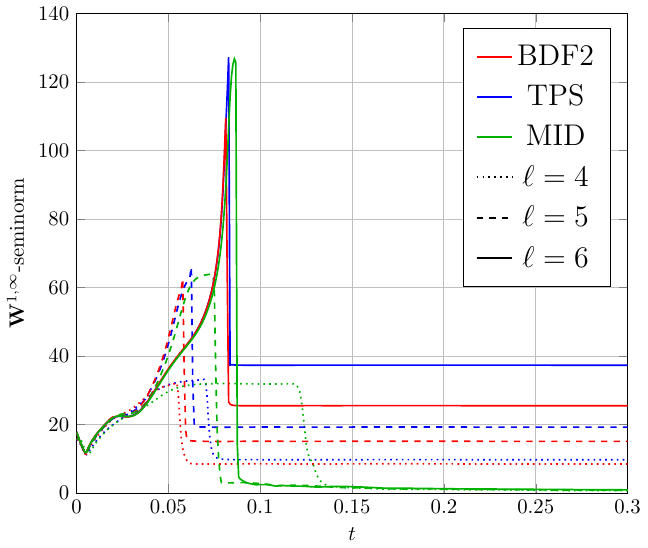}
    \end{subfigure}
    \begin{subfigure}[b]{0.52\linewidth}  
        \centering
        \includegraphics[width=1\linewidth]{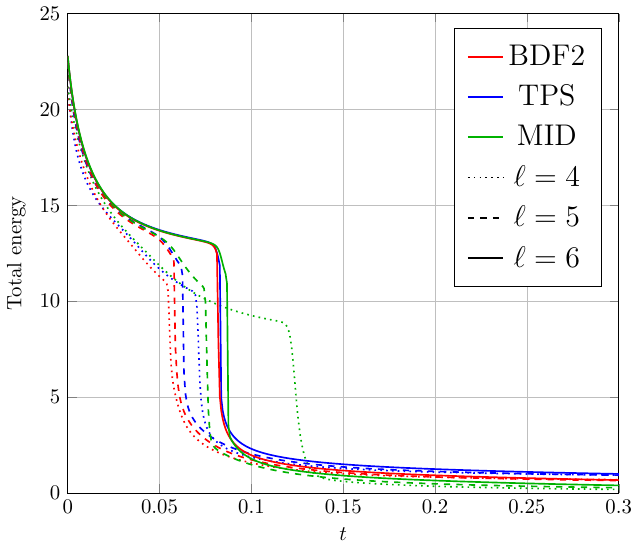}
    \end{subfigure}
    \centering
    \begin{subfigure}[b]{0.52\linewidth}
        \centering
        \includegraphics[width=1\linewidth]{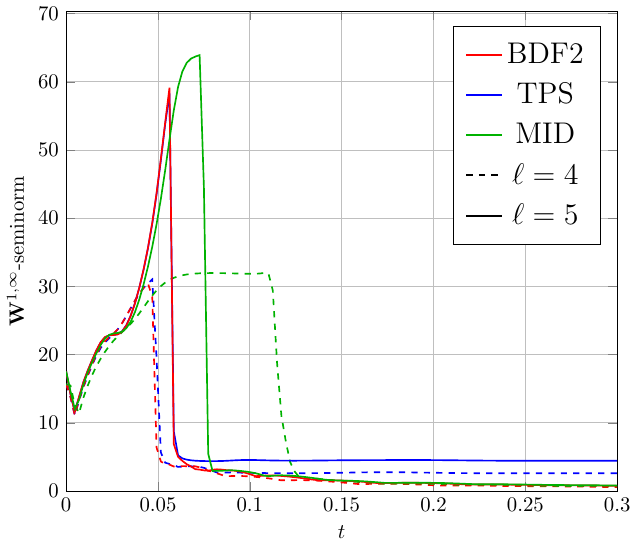}
    \end{subfigure}
    \begin{subfigure}[b]{0.52\linewidth}  
        \centering
        \includegraphics[width=1\linewidth]{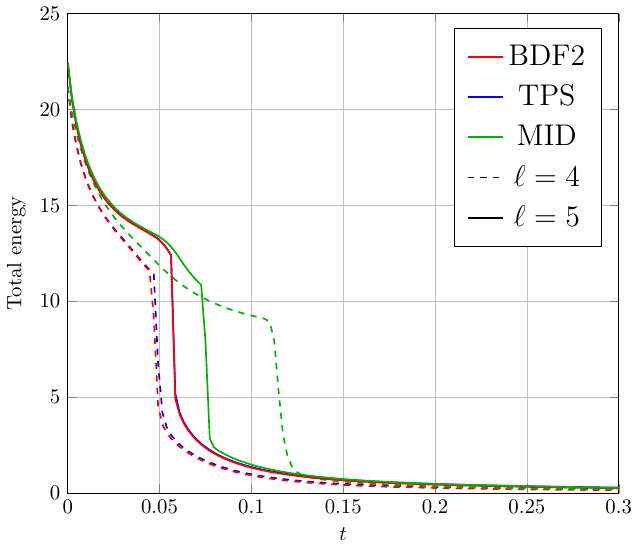}
    \end{subfigure}
    \caption{Evolution of the $\mathbf{W}^{1,\infty}$-seminorm in time (left) and energy decay (right) for $\tau = h^2/10$ (top) and $\tau = h^2/100$ (bottom), depending on the numerical scheme.}
    \label{fig:example_blowup_schemes} 
\end{figure}

Finally, we compare the evolution of the energy and the $\mathbf{W}^{1,\infty}$-seminorm for three different schemes, i.e., BDF2, MID, and TPS. As before, we consider a uniform mesh consisting of $2^{2 \ell+1}$ rectangular triangles, with $\ell\in\{4,5,6\}$, i.e., corresponding to mesh sizes $h_4~=~2^{-4}\sqrt{2}\approx 0.0884$, $h_5 = 2^{-5}\sqrt{2}\approx 0.0442$, and $h_5 = 2^{-6}\sqrt{2}\approx 0.0221$, respectively, and we set the time-step size as $\tau_\ell = h_\ell^2/10$ and $\tau_\ell = h_\ell^2/100$. The results are reported in Figure \ref{fig:example_blowup_schemes}. 

We remark that MID requires the solution of a nonlinear system at each time-step, which we solve by a fixed-point iteration as proposed in \cite{bp2006,mrs2018}, requiring the CFL condition $\tau = o(h^2)$ for the convergence analysis. Performing the numerical experiments, we noticed that the relation $\tau = h^2$ is not enough to guarantee the convergence of the fixed-point iteration, and a smaller time-step size is needed. This is probably due to the choice of the fixed-point iteration, rather than MID itself, and other nonlinear solvers could probably improve this aspect. In that sense, the BDF2 scheme seems to allow more freedom in the choice of the time-step size.

Nevertheless, we observe that TPS and the BDF2 scheme coincide, while MID yields a blow-up time that is slightly delayed.

\subsection{Numerical deviation from the unit-length constraint}
\label{subsec: fourth}
In this section, we investigate the deviation from the unit-length constraint $|\mm(t)|~=~1$ for every $t \in [0,T]$. To this end, we compare the growth of $t \mapsto \norm{\mm_{h\tau}(t)}_{\bL^\infty(\Omega)}$ for TPS and for BDF2.

As a first test case, we consider the setting from Section~\ref{subsec: first}. We fix the exchange constant $\lamex^2 = 0.01$ so that the external field dominates on the exchange contribution of the effective field. 
As a second test case, we consider the experiment from Section~\ref{subsec: second}.

For both examples, we employ a fixed uniform triangulation with mesh size $h = 2^{-4}\sqrt{2}\approx 0.0884$ and a uniform partition of the time interval with time-step size $\tau_0= 4 \cdot 10^{-3}$ that is successively bisected, i.e., $\tau_k = 2^{-k} \tau_0$ for $k = 0, \ldots, 4$.

\begin{figure}[ht!]
    \centering
    \begin{subfigure}[b]{0.52\linewidth}
        \centering
         \includegraphics[width=0.95\linewidth]{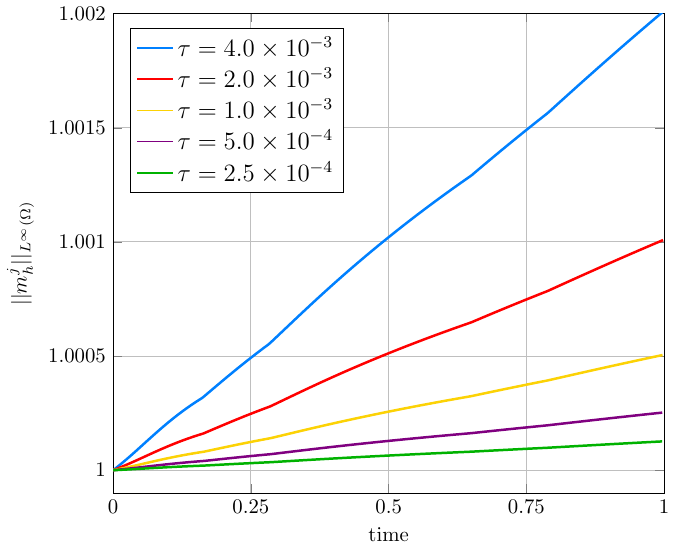}
         \caption{TPS, Experiment from Section~\ref{subsec: first}.}
    \end{subfigure}
    \begin{subfigure}[b]{0.53\linewidth}  
        \centering
        \includegraphics[width=0.95\linewidth]{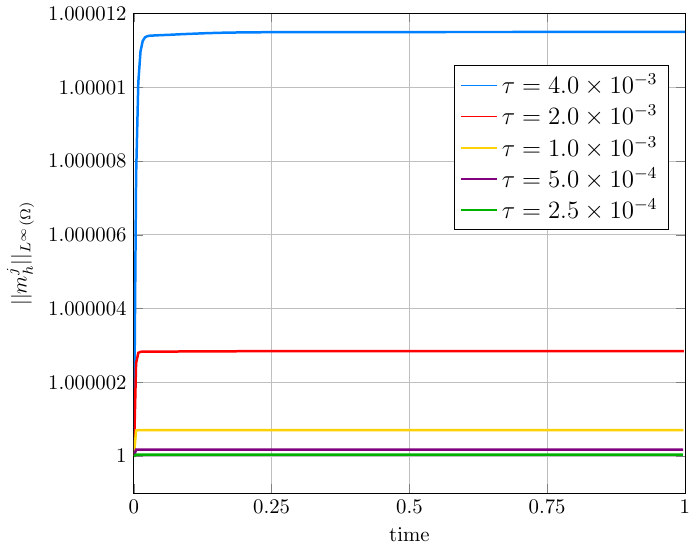}
        \caption{Alg.~\ref{alg: full discr} (BDF2), Experiment from Section~\ref{subsec: first}.}
    \end{subfigure}
    \centering
    \begin{subfigure}[b]{0.52\linewidth}
        \centering
        \includegraphics[width=0.95\linewidth]{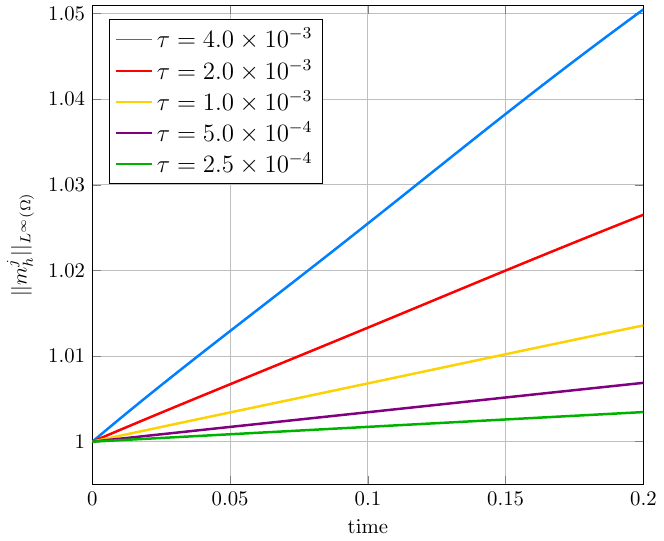}
        \caption{TPS, Experiment from Section~\ref{subsec: second}.}
    \end{subfigure}
    \begin{subfigure}[b]{0.53\linewidth}  
        \centering
        \includegraphics[width=0.95\linewidth]{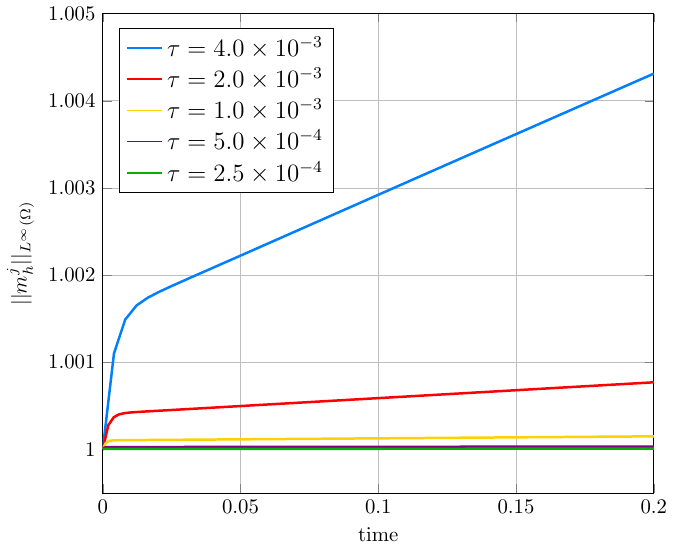}
    \caption{Alg.~\ref{alg: full discr} (BDF2), Experiment from Section~\ref{subsec: second}.}
    \end{subfigure}
    \caption{Empirical deviation from the unit-length constraint for TPS (left) and BDF2 (right) for the solution from Section \ref{subsec: first} (top) and the solution from Section \ref{subsec: second} (bottom).}
\label{fig:example_1_unit_modulus} 
\end{figure}

The first row of Figure \ref{fig:example_1_unit_modulus} shows the evolution for the first test case, while the second row considers the second one. It is clearly visible (right column) that in either experiment BDF2 is much superior to TPS; in particular, note the scaling of the $y$-axis.

Moreover, for both examples, we examine the convergence rate of the deviation of the discrete solutions from the unit sphere. According to Remark~\ref{remark: second-order convergence rate}, this rate is expected to be of second-order for BDF2, i.e.,
$$
\norm{|\mm_h^n|^2-1}_{L^1(\Omega)}\le C\tau^2,
$$ 
provided that the regularity condition
$$
\normO{\vv_h^0}^2+\tau^2\sum_{j=2}^n\normO{\d _t^2\mm_h^j}^2\lesssim 1
$$
is satisfied.
In both test cases, we use the same physical parameters of before. For the first test case, we discretize in time using a uniform partition of the time interval with an initial time-step size of $\tau_0 = 4 \times 10^{-3}$, and we set $\tau_\ell = 2^{-\ell}\tau_0$ for $\ell = 0, \dots, 6$. For the second test case, we start with an initial time-step size of $\tau_0 = 0.064$, which is successively halved nine times, i.e., $\tau_\ell = \tau_0 \cdot 2^{-\ell}$ for $\ell = 0, \dots, 9$.

\begin{figure}[ht!]
\centering
\begin{subfigure}[b]{0.47\linewidth}  
  \centering
  \includegraphics[width=1\linewidth]{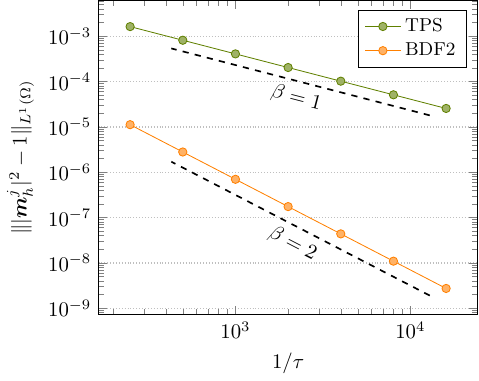}
  \caption{Experiment from Section \ref{subsec: first}.}
\end{subfigure} 
\begin{subfigure}[b]{0.47\linewidth}  
  \centering
  \includegraphics[width=1\linewidth]{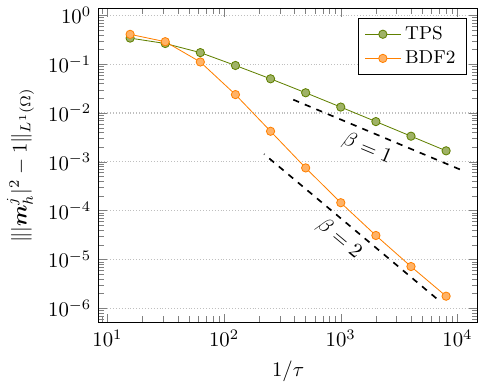}
  \caption{Experiment from Section \ref{subsec: second}.}
\end{subfigure}
\caption{First-order convergence (TPS) and second-order convergence (BDF2) for the deviation of the discrete solutions from the modulus one constraint for the solutions from Section \ref{subsec: first}--\ref{subsec: second}.}
\label{fig:unit_modulus_conv_rates}
\end{figure}

Figure~\ref{fig:unit_modulus_conv_rates} presents the convergence rates obtained for both examples. We observe first-order convergence for TPS and second-order convergence for BDF2. In Tables~\ref{tab:bounds1}-\ref{tab:bounds2}, we see that the quantities $\normO{\vv_h^0}^2$ and $\tau^2\sum_{j=2}^n\normO{\d _t^2\mm_h^j}^2$ are indeed bounded.

\begin{figure}[ht!]
    \centering
    \begin{subfigure}[b]{0.49\linewidth}
        \centering
        \includegraphics[width=1\linewidth]{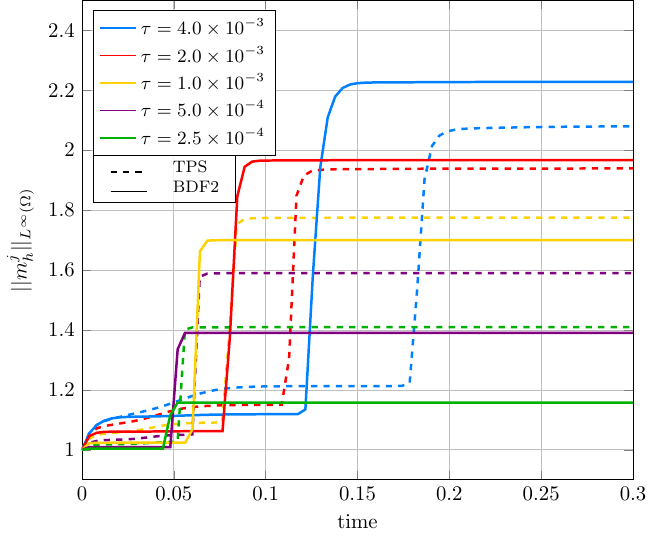}
        \caption{Growth in time of $\norm{\mm_{h\tau}(t)}_{\bL^\infty(\Omega)}$.}
    \end{subfigure}
    \begin{subfigure}[b]{0.51\linewidth}  
        \centering
        \includegraphics[width=1\linewidth]{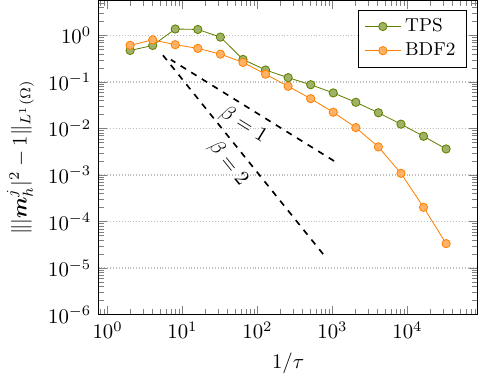}
        \caption{Convergence rate of $\norm{|\mm_h^j|^2-1}_{L^1(\Omega)}$}
    \end{subfigure}
    \caption{Deviation from the unit sphere for the singular solution from Section \ref{subsec: third}.}
    \label{fig:example_blowup_2} 
\end{figure}

\begin{table}[ht!]
\centering
\begin{tabular}{c|c|c|c}
\toprule
$\tau$ & $\normO{\vv_h^0}^2$ & $\tau^2\sum\limits_{j=2}^n\normO{\d _t^2\mm_h^j}^2$ (TPS) & $\tau^2\sum\limits_{j=2}^n\normO{\d _t^2\mm_h^j}^2$ (BDF2) \\
\midrule
$4\times10^{-3}$   & $0.455689$ & $1.63774\times10^{-3}$ & $1.65129\times10^{-3}$ \\
$2\times10^{-3}$   & $0.455708$ & $8.22262\times10^{-4}$ & $8.25698\times10^{-4}$ \\
$1\times10^{-3}$   & $0.455718$ & $4.11988\times10^{-4}$ & $4.12853\times10^{-4}$ \\
$5\times10^{-4}$   & $0.455723$ & $2.06209\times10^{-4}$ & $2.06426\times10^{-4}$ \\
$2.5\times10^{-4}$ & $0.455726$ & $1.03159\times10^{-4}$ & $1.03213\times10^{-4}$ \\
$1.25\times10^{-4}$& $0.455727$ & $5.15928\times10^{-5}$ & $5.16064\times10^{-5}$ \\
$6.25\times10^{-5}$& $0.455727$ & $2.57998\times10^{-5}$ & $2.58032\times10^{-5}$ \\
\bottomrule
\end{tabular}
\caption{Reported values for the experiment from Section \ref{subsec: first}.}\label{tab:bounds1}
\begin{tabular}{c|c|c|c}
\toprule
$\tau$ & $\normO{\vv_h^0}^2$ & $\tau^2\sum\limits_{j=2}^n\normO{\d _t^2\mm_h^j}^2$ (TPS) & $\tau^2\sum\limits_{j=2}^n\normO{\d _t^2\mm_h^j}^2$ (BDF2) \\
\midrule
0.064   & 12.04981 & 228.36674 & 65.90581 \\
0.032   & 42.57946 & 402.77878 & 157.35153 \\
0.016   & 59.76109 & 337.40996 & 236.61697 \\
0.008   & 65.07227 & 197.07895 & 185.57028 \\
0.004   & 66.48977 & 108.21309 & 109.69379 \\
0.002   & 66.85476 & 56.71912  & 57.78277 \\
0.001   & 66.95866 & 29.05572  & 29.41360 \\
0.0005  & 67.02836 & 14.75316  & 14.85658 \\
0.00025 & 67.18207 & 7.57171   & 7.62985  \\
0.000125& 67.51339 & 4.10776   & 4.29970  \\
\bottomrule
\end{tabular}
\caption{Reported values for the experiment from Section \ref{subsec: second}.}\label{tab:bounds2}
\end{table}
Finally, we investigate the conservation of length also in case of the singular solution from Section~\ref{subsec: third}, showing the deviation from the sphere when the singularity formation happens. We choose $T=1$ and $\tau_\ell \coloneq 2^{-\ell}\tau_0$ with $\tau_0=5\times10^{-1}$ and $\ell=0,\dots,15$. Figure \ref{fig:example_blowup_2} shows a significantly worse behavior than for the previous experiments, for the growth in time of $\norm{\mm_{h\tau}(t)}_{\bL^\infty(\Omega)}$, the convergence rate of the deviation from the unit-length constraint and the boundedness of the quantities $\normO{\vv_h^0}^2$ and $\tau^2\sum_{j=2}^n\normO{\d _t^2\mm_h^j}^2$. We observe that the unit modulus constraint is not respected (see the values on the $y$-axis in Figure \ref{fig:example_blowup_2}, left). Moreover, also the convergence rate of the deviation from the unit-length constraint deteriorates (see Figure \ref{fig:example_blowup_2}, right) and the reported values in Table~\ref{tab:bounds3} are significantly bigger than before. 

\begin{table}[ht!]
\centering
\begin{tabular}{c|c|c|c}
\toprule
$\tau$ & $\normO{\vv_h^0}^2$ & $\tau^2\sum\limits_{j=2}^n\normO{\d_t^2\mm_h^j}^2$ (TPS) & $\tau^2\sum\limits_{j=2}^n\normO{\d_t^2\mm_h^j}^2$ (BDF2) \\
\midrule 
$5.0\times10^{-1}$ & 1.42801 & 0.55275 & 0.49728 \\ 
$2.5\times10^{-1}$ & 5.23960 & 2.31982 & 3.81792 \\ 
$1.25\times10^{-1}$ & 16.97473 & 14.81064 & 9.60509 \\ 
$6.25\times10^{-2}$ & 43.03307 & 35.97531 & 44.50748 \\ 
$3.125\times10^{-2}$ & 85.88250 & 52.71513 & 159.06115 \\ 
$1.5625\times10^{-2}$ & 164.78944 & 192.54834 & 459.63490 \\ 
$7.8125\times10^{-3}$ & 347.73076 & 432.46000 & 967.93781 \\ 
$3.9063\times10^{-3}$ & 777.98794 & 1013.63901 & 1921.49996 \\ 
$1.9531\times10^{-3}$ & 1525.54306 & 2128.73352 & 3770.10938 \\ 
$9.7656\times10^{-4}$ & 2323.21045 & 4222.58285 & 7126.33462 \\ 
$4.8828\times10^{-4}$ & 2937.72530 & 7093.48840 & 12234.96584 \\ 
$2.4414\times10^{-4}$ & 3447.27162 & 8924.46862 & 17014.57798 \\ 
$1.2207\times10^{-4}$ & 3937.01609 & 8497.24954 & 16846.74545 \\ 
$6.1035\times10^{-5}$ & 4318.47730 & 6598.78479 & 11240.96887 \\ 
$3.0518\times10^{-5}$ & 4533.15425 & 4497.07936 & 6258.06810 \\ 
\bottomrule
\end{tabular}
\caption{Reported values for the singular experiment from Section \ref{subsec: third}.}\label{tab:bounds3}
\end{table}

\subsection{Role of the CFL condition $\boldsymbol{\tau = o(h^2)}$}\label{subsec: CFL}
In this section, we investigate the role of the CFL condition $\tau = o(h^2)$, which is used as a theoretical requirement to verify the energy estimate of Definition~\ref{def: weak solution}\eqref{def: weak solution iv} in the proof of Theorem \ref{thm: main result}.

\begin{figure}[ht!]
    \centering
    \begin{subfigure}[b]{0.51\linewidth}
        \centering
        \includegraphics[width=1\linewidth]{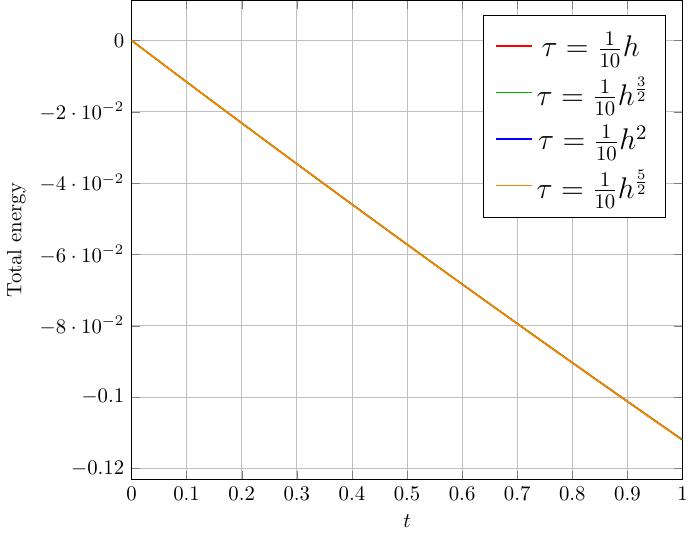}
        \caption{Experiment from Section \ref{subsec: first}.}
    \end{subfigure}
    \begin{subfigure}[b]{0.51\linewidth}  
        \centering
        \includegraphics[width=0.93\linewidth]{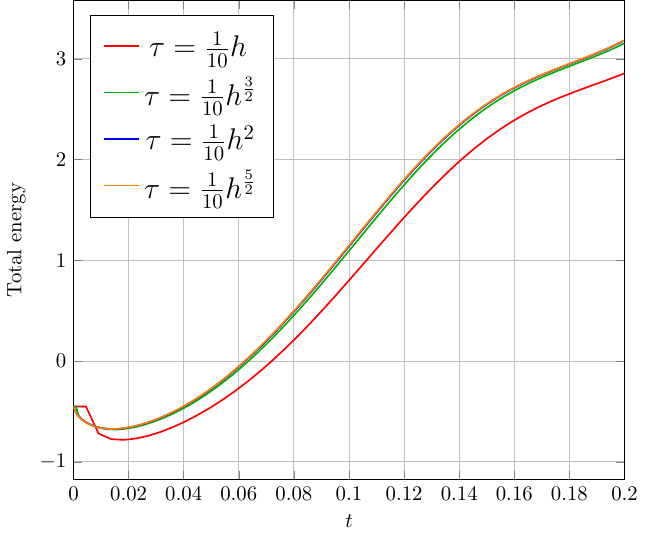}
        \caption{Experiment from Section \ref{subsec: second}.}
    \end{subfigure}
    \caption{Evolution of the total Gibbs energy for fixed mesh size $h = 2^{-5}\sqrt{2}$ and different time-step sizes $\tau=\frac{1}{10}h^\alpha$ and $\alpha~\in~\{1; 3/2; 2; 5/2\}$.}
    \label{fig:energy_evolution_smooth_examples} 
\end{figure}

First, we plot the decay in time of the total Gibbs energy 
$$
    \EE(\mm,\ff)
    =
    \frac{\lamex^2}{2} \int_\Omega|\nabla\mm|^2\d \xx - \int_\Omega\ff\cdot\mm \d \xx
$$ 
for the smooth solutions from Section~\ref{subsec: first} and \ref{subsec: second}. For both cases we fix a mesh size $h = 2^{-5}\sqrt{2}$ and we investigate the energy decay for different time-step sizes $\tau=\frac{1}{10}h^\alpha$ with $\alpha\in\left\{1; \frac{3}{2}; 2; \frac{5}{2}\right\}$. In Figure \ref{fig:energy_evolution_smooth_examples} we plot the evolution of the total energy for both examples. 
In order to avoid an overloading of plots, only the results obtained with BDF2 are reported. 
Notice that, for the test case from Section~\ref{subsec: second}, the total energy $\EE(\mm(t),\ff(t))$ is not guaranteed to be decreasing in time, since the external field $\ff$ is not constant in time.

\begin{figure}[ht!]
  \centering
  \begin{subfigure}[b]{0.49\linewidth}
    \centering
    \includegraphics[width=1\linewidth]{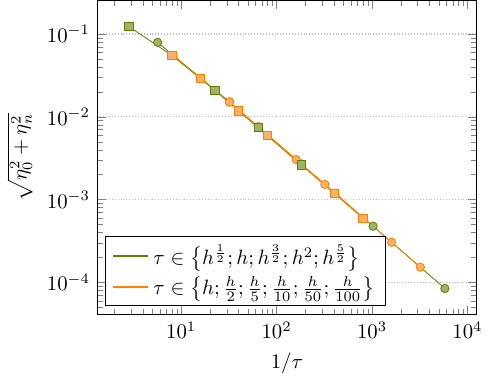}
    \caption*{\small Experiment from Section \ref{subsec: first}.}
  \end{subfigure}
  \begin{subfigure}[b]{0.49\linewidth}
    \centering
    \includegraphics[width=1\linewidth]{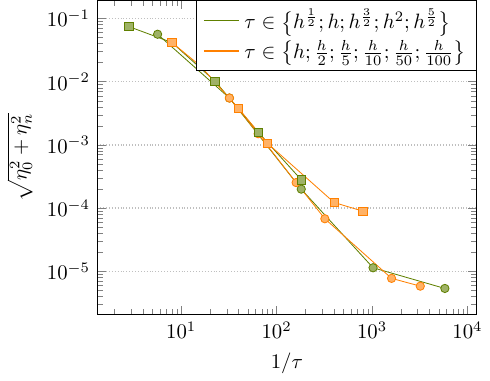}
    \caption*{\small Experiment from Section \ref{subsec: second}.}
  \end{subfigure}

  \vspace{6pt}
  \begin{subfigure}[b]{0.49\linewidth}
    \centering
    \includegraphics[width=1\linewidth]{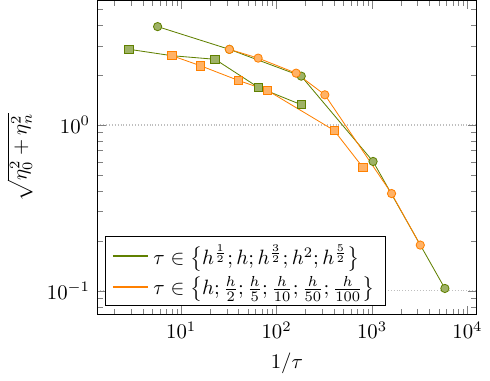}
    \caption*{\small Experiment from Section \ref{subsec: third}.}
  \end{subfigure}

  \caption{Decay of $C(\tau) = \sqrt{\eta_0^2 + \eta_n^2}$ for $h = 0.125$ (\protect\tikz[baseline=-0.5ex]{
  \protect\draw[line width=1pt] (-0.4,0) -- (-0.15,0); 
  \protect\draw[line width=1pt] (-0.15,-0.15) rectangle (0.15,0.15); 
  \protect\draw[line width=1pt] (0.15,0) -- (0.4,0); 
 })
 and $h = 0.03125$ 
 (\protect\tikz[baseline=-0.5ex]{
  \protect\draw[line width=1pt] (-0.4,0) -- (-0.15,0); 
  \protect\draw[line width=1pt] (0,0) circle [radius=0.15]; 
  \protect\draw[line width=1pt] (0.15,0) -- (0.4,0); 
 }) and different choices of $\tau$ decaying as fractions or powers of $h$.}
  \label{fig:eta_decay}
\end{figure}

Second, we fix a mesh size $h$ and investigate the decay of the quantity
$$
C(\tau) \coloneqq \sqrt{\eta_0^2 + \eta_n^2} 
$$
as $\tau\to 0$, for all the three experiments from Sections \ref{subsec: first}--\ref{subsec: third}. This is motivated by the fact that the CFL condition is introduced only in order to guarantee that $C(\tau)\to 0$ as $\tau h^{-2}\to 0$. We choose $\tau$ decaying to zero as fractions of $h$, i.e., $\tau \in \big\{h; \frac{h}{2}; \frac{h}{5}; \frac{h}{10}; \frac{h}{50}; \frac{h}{100}\big\}$ or as powers of $h$, i.e., $\tau \in \{h^{\frac{1}{2}}; h; h^{\frac{3}{2}}; h^{2}; h^{\frac{5}{2}} \}.$ For each case, we choose two fixed mesh sizes, corresponding to a coarse mesh ($256$ elements; mesh size $h = 0.125$) and to a finer one ($4096$ elements; mesh size $h = 0.03125$). Figure \ref{fig:eta_decay} illustrates the decay of $C(\tau)$, which happens always and independently of the choice of $\tau$ as fractions or powers of $h$, and independently of the chosen mesh size. This underlines that the condition $C(\tau)\to 0$ appears to be satisfied in the present examples, so that the CFL condition $\tau = o(h^2)$ might only be a technical requirement in the proof of Theorem \ref{thm: main result} and is indeed not needed.
\appendix
\section*{Appendix}
\begin{proof}[{\bfseries Proof of Lemma~\ref{lem: norm equivalence 2}}]
From Algorithm~\ref{alg: full discr}\eqref{item: def mm1} and Lemma~\ref{lem: identities}, it follows that
\begin{equation*}
 \begin{split}
 &\tau\sum_{j=0}^{n-1}\normO{\vv_h^j}^2
 \eqreff{eq: identity1}=
 \tau\normO{\d_t\mm_h^1}^2+\tau\sum_{j=1}^{n-1}\normO[\Big]{\frac{3}{2}\d_t\mm_h^{j+1}-\frac{1}{2}\d _t\mm_h^{j}}^2
 \\& \qquad
 \le
 \tau\normO{\d_t\mm_h^1}^2+\frac{9}{2}\tau\sum_{j=1}^{n-1}\normO{\d_t\mm_h^{j+1}}^2+\frac{1}{2}\tau\sum_{j=1}^{n-1}\normO{\d_t\mm_h^j}^2
 \le
 5\tau\sum_{j=0}^{n-1}\normO{\d_t\mm_h^{j+1}}^2.
 \end{split}
\end{equation*}
This proves the lower estimate in~\eqref{eq: norm equivalence 2}.
Analogously, it also holds that
\begin{equation*}
 \begin{split}
 \tau\sum_{j=0}^{n-1}\normO{\d_t\mm_h^{j+1}}^2
 &\eqreff{eq: identity1}=
 \tau\normO{\vv_h^0}^2+\tau\sum_{j=1}^{n-1}\normO[\Big]{\frac{2}{3}\vv_h^j-\frac{1}{3}\d_t\mm_h^j}^2
 \\&
 \le
 \tau\normO{\vv_h^0}^2+\frac{8}{9}\tau\sum_{j=1}^{n-1}\normO{\vv_h^j}^2+\frac{2}{9}\tau\sum_{j=1}^{n-1}\normO{\d_t\mm_h^j}^2.
 \end{split}
\end{equation*}
Absorbing the last term of the right-hand side by the left-hand side, it follows that 
\begin{equation*}
 \begin{split}
 \frac{7}{9} \, \tau\sum_{j=0}^{n-1}\normO{\d_t\mm_h^{j+1}}^2
 \le
 \tau\sum_{j=0}^{n-1}\normO{\vv_h^j}^2.
 \end{split}
\end{equation*}
This proves the upper estimate in~\eqref{eq: norm equivalence 2}.
The same argument proves~\eqref{eq: norm equivalence grad 2}.
\end{proof}

\begin{proof}[{\bfseries Proof of Lemma~\ref{lem: inverse estimate}}]
Note that 
\begin{align*}
 \normO{\d_t^2\mm_h^j}^2
 = \frac{1}{\tau^2} \, \normO{\d_t\mm_h^j - \d_t\mm_h^{j-1}}^2
 \le \frac{2}{\tau^2} \big( \normO{\d_t\mm_h^j}^2 + \normO{\d_t\mm_h^{j-1}}^2 \big).
\end{align*}
Together with~\eqref{eq: norm equivalence 2} from Lemma~\ref{lem: norm equivalence 2}, this yields that
\begin{align*}
 \tau\sum_{j=2}^n \normO{\d_t^2\mm_h^j}^2 
 \le \frac{2}{\tau} \sum_{j=1}^{n} \normO{\d_t\mm_h^j}^2 
 \eqreff{eq: norm equivalence 2}\le \frac{18/7}{\tau^2} \, \Big(\tau \sum_{j=0}^{n-1} \normO{\vv_h^j}^2 \Big).
\end{align*}
This proves~\eqref{eq: inverse estimate} and the same argument proves~\eqref{eq: inverse estimate grad}.
\end{proof}
\renewcommand{\thetheorem}{A.\arabic{theorem}}
\setcounter{theorem}{0} 

\printbibliography

\end{document}